\documentclass{article}

\usepackage{amsmath}
\usepackage{amssymb}
\usepackage{url}
\usepackage{graphicx,color}
\usepackage{geometry}
\input otex.sty

\def\sig{\sigma}
\def\i{\,\ro i\,}

\newdimen\fiiiskip \fiiiskip=-.4in
\newdimen\figwidth \figwidth=4cm
\newdimen\figboxl \figboxl=1.5in
\newdimen\figlabl \figlabl=1.65in
\newdimen\prefigskip \prefigskip=15pt
\newdimen\figunderskip \figunderskip=5pt

\def\Fig#1#2#3{\begin{figure}[htbp]#3\caption{#2.\label{#1}}\end{figure}}

\def\Figbox#1{\includegraphics[width=\figwidth]{#1.pdf}}

\def\Figboxl#1#2{\hbox to \figboxl{\vbox{\hbox to \figboxl{\hfil\includegraphics[width=\figwidth]{#1.pdf}\hfil}\vglue\figunderskip\hbox to \figlabl{\hfil #2\hfil}}\hfil}}

\def\Fiiil#1#2#3#4#5#6{\hfill\break\hglue\prefigskip
\Figboxl{#1}{$t = #2$}\hfill\Figboxl{#3}{$t = #4$}\hfill\Figboxl{#5}{$t = #6$}\hfill\break\hfill}

\begin{document}

\title{Dispersive Lamb Systems}
\date{\today}
\author{Peter J. Olver$^1$ and Natalie E. Sheils$^2$ \\
\footnotesize 1. School of Mathematics, University of Minnesota, Minneapolis, MN 55455\\ \footnotesize{\texttt{olver@umn.edu}}\qquad  \footnotesize{\url{http://www.math.umn.edu/~olver}}\\
\footnotesize 2. School of Mathematics, University of Minnesota, Minneapolis, MN 55455\\ \footnotesize{\texttt{nesheils@umn.edu}}\qquad \footnotesize{\url{http://www.math.umn.edu/~nesheils}}
}

\maketitle

\vskip-.2in

\Abstract
Under periodic boundary conditions, a one-dimensional dispersive medium driven by a Lamb oscillator exhibits a smooth response when the dispersion relation is asymptotically linear or superlinear at large wave numbers, but unusual fractal solution profiles emerge when the dispersion relation is asymptotically sublinear. Strikingly, this is exactly the opposite of the superlinear asymptotic regime required for fractalization and dispersive quantization, also known as the Talbot effect, of the unforced medium induced by discontinuous initial conditions.  


\Section i Introduction.

The original Lamb system was introduced by Horace Lamb in 1900, \rf{Lambp}, as a simple model for the phenomenon of \is{radiation damping} experienced by a vibrating body in an energy conducting medium. Examples include vibrations of an elastic sphere in a gaseous medium, electrical oscillations of a spherical conductor, a dielectric sphere with large inductance, relativistic radiation of energy from a concentrated mass via gravity waves, and quantum resonance of nuclei. To understand this phenomenon, Lamb proposed a simpler one-dimensional model consisting of an oscillatory point mass--spring system directly coupled to an infinite string modeled by the usual second order wave equation. The oscillator transfers energy to the string, generating waves that propagate outwards at the intrinsic wave speed. Meanwhile, the string tension induces a damping force on the oscillator, with the effect that the propagating waves are of progressively smaller and smaller amplitudes. It is worth pointing out a potential source of confusion: Lamb refers to the point mass as a ``nucleus'', although he clearly does not have atomic nuclei in mind as these were not discovered by Rutherford until 1911. 

More recently, Hagerty, Bloch, and Weinstein, \rf{HaBlWe}, investigated mechanical gyroscopic systems, also known as Chetaev  systems, such as the spherical pendulum and rigid bodies with internal rotors, that are coupled to both the classical non-dispersive wave equation and a dispersive equation of Klein--Gordon form.  In such systems, the gyroscopic Lamb oscillator can even induce instabilities through the effects of Rayleigh dissipation. However, as we will see, their analysis does not cover the full range of phenomena exhibited by dispersive Lamb models.

This paper analyzes the both the full line and periodic Lamb problems when the oscillatory mass is coupled to more general one-dimensional dispersive media. Our motivation is to see whether such systems exhibit the remarkable phenomenon of dispersive quantization, \rf{Odq}, also known as the Talbot effect, \rf{BMS}, in honor of a striking 1836 optical experiment of William Henry Fox Talbot, \rf{Talbot}. Fractalization and quantization effects arise in both linear and nonlinear unforced dispersive wave equations on (at least) periodic domains. In the linear regime, a solution to the periodic initial-boundary value problem produced by rough initial data, \eg a step function, is ``quantized'' as a discontinuous but piecewise constant profile or, more generally, piecewise smooth or even piecewise fractal, at times which are rational multiples of $L^n/\pi^{n-1}$, where $L$ is the length of the interval, but exhibits a continuous but non-differentiable fractal profile at all other times. (An interesting question is whether such fractal profiles enjoy self-similarity properties similar to the Riemann and Weierstrass non-differentiable functions, as explored in \rf{Duistermaat}.  However, the Fourier series used to construct the latter are quite different in character, and so the problem remains open and challenging.)

Dispersive quantization relies on the slow, conditional convergence of the Fourier series solutions, and requires the dispersion relation to be asymptotically polynomial at large wave number. 
The key references include the 1990's discovery of Michael Berry and collaborators, \rf{BMS}, in the context of optics and quantum mechanics, and the subsequent analytical work of Oskolkov, \rf{OskolkovV}. In particular, this effect underlies the experimentally observed phenomenon of quantum revival, in which an electron that is initially concentrated at a single location of its orbital shell is, at rational times, re-concentrated at a finite number of orbital locations. The subsequent rediscovery by the first author, \rf{COdisp, Odq}, showed that the phenomenon appears in a range of linear dispersive partial differential equations, while other models arising in fluid mechanics, plasma dynamics, elasticity, DNA dynamics, etc.  exhibit a fascinating range of as yet poorly understood behaviors, whose qualitative features are tied to the large wave number asymptotics of the dispersion relation. These studies were then extended, through careful numerical simulations, to  nonlinear dispersive equations, including both integrable models, such as the cubic nonlinear Schr\"odinger, Korteweg--deVries, and modified Korteweg--deVries equations, as well as non-integrable generalizations with higher degree nonlinearities. Some of these numerical observations were subsequently rigorously confirmed in papers of Chousionis, Erdo\utxt gan, and Tzirakis, \rf{ChErTz,ErTznls,ErTz}, and, further, in the very recent paper by Erdo\utxt gan and Shakan, \rf{ErSh}, which extends the analysis to non-polynomial dispersion relations, but much more work remains to be completed, including extensions to other types of boundary conditions and dispersive wave models in higher space dimensions. 

In the original case of a Lamb oscillator coupled to the wave equation, but now on a periodic domain, the solution is readily seen to remain piecewise smooth, the only derivative discontinuities occurring at the wave front and the location of the oscillator. Extending the analysis to a higher order linearly elastic string model reveals that the solution continues to exhibit a piecewise smooth profile. Indeed, smoothness will hold for both bidirectional and unidirectional dispersive linear partial differential equations with linear or super-linear dispersion asymptotics at large wave numbers, including the Klein--Gordon model considered in \rf{HaBlWe}. However, the periodic Lamb problem with dispersion asymptotics of the order $\O{\sqrt{\abs k}}\ustrut{10}$ as the wave number $k \to \infty $ does exhibit a fractal solution profile. Furthermore, recent results of Erdo\utxt gan and Shakan, \rf{ErSh} can be adapted to the present context, to rigorously establish convergence of the formal Fourier series solution, at each fixed time $t >0$, to a fractal profile whose maximal fractal dimension lies between $\fr54$ and $\fr74$.  On the other hand, when the dispersion relation is asymptotically constant, the resulting highly oscillatory plots of the partial sums suggest weak convergence of the Fourier series to some kind of distributional solution. 

%

%

Our present analysis of the Lamb systems relies on classical Fourier transform and Fourier series techniques. One can envision applying the more general and powerful Unified Transform Method (UTM), due to Fokas and collaborators, \rf{FokasUTM, Fokas}. However, the space-dependent coefficient places the system outside the class of equations currently solvable by the UTM. A second possible way to approach such problems is to view Lamb's original formulation as an interface problem, as in~\rf{DePeSh, ShDe}, and combine this with recent work applying the UTM to systems of equations,~\rf{DGSV}. 
However, the Lamb interface condition is more complicated than those considered to date, and extending current work on interface problems remains an interesting challenge.


\Rmk The paper includes still shots of a variety of solutions at selected times. \Mathematica\ code for generating the movies, which are even more enlightening, can be found on the first author's web site:\hfbg{100} {\www math.umn.edu/$\sim$olver/lamb }

\Section 2 The Bidirectional Lamb Model.

The original Lamb model, \rf{Lambp}, consists of an oscillating point mass that is connected to an infinite elastic string and constrained to move in the transverse direction. With the system starting at rest, the mass is subject to a sudden blow, and the continuous medium serves to dampen its ensuing vibrations of the mass. In the underdamped regime, the vibrating mass spawns an oscillatory traveling wave which propagates along the string at its innate wave speed, while the string tension acts as a damping force on the oscillating mass, which vibrates at an exponentially decreasing amplitude that in turn produces a similarly decaying propagating wave profile. 

Thus, in the small amplitude regime, the string displacement $u(t,x)$ satisfies the usual bidirectional wave equation
\Eq{wave}
$$\qqeq{u_{tt} = c^2 u_{xx}, &x \ne 0,}$$ 
away from the fixed location of the mass, which we take to be at the origin $x = 0$. 
Here $c = \sqrt{T/\rho }$ is the wave speed, with $\rho $ representing the string density and $T$ its tension, which are both assumed constant for simplicity.
As in \rf{Lambp}, force balance on the mass displacement $h(t) = u(t,0)$ at the origin yields
\Eq{fb}
$$M \pa{h'' + \sigma^2 h} = -\: T \,\bbk{ u_x}_0,$$
where $M$ is the mass and $\sigma$ its uncoupled oscillatory frequency. The forcing term on the right hand side is the negative of the product of the string tension $T$ and the jump in the spatial derivative $u_x$ at the location of the mass.  \is{Warning\/}: Lamb appears to include a factor of $2$ in his formulation of \eq{fb}, but he is using the fact that $u$ is even in $x$ and hence its jump at the origin is twice its limiting value. 

The first observation is that the Lamb model can be rewritten as a forced wave equation of the form
\Eq{fwave}
$$u_{tt} = c^2 u_{xx} - 2\:c\:h'(t) \, \delta (x),$$
where $\delta (x)$ is the usual Dirac delta function and $h(t)$ satisfies the damped oscillator equation
\Eq{mass}
$$\req{h'' + 2\:\beta \: h' + \sig ^2 h = 0, \\ h(0) = 0,}$$
whose damping coefficient is 
$$\beta 
 = \frac T{c\:M} = \frac{\sqrt{\rho \:T\,\ustrut9}} M \, 
\,.$$
To prove this, just integrate \eq{fwave} from $x = -\: \varepsilon $  to $x = +\: \varepsilon $ and let $\varepsilon \to 0$.

According to \rf{P; Theorem 2.18}, the solution to \eq{fwave} with zero initial conditions
\Eq{ic0}
$$\req{u(0,x) = u_t (0,x) = 0,}$$
 (\ie at rest initially) is
\Eq{Fwavesol}
$$\req{u(t,x) = -\;\Int s0t{\Int y{x-c\,(t-s)}{x+c\,(t-s)}{h'(s) \,\delta (y)} } ,\\t > 0.}$$
Integrating twice, we deduce that
\Eq{sol}
$$\req{u(t,x) = -\,\htilde\bpa{c\:t - \abs x} ,\\t > 0,}$$
where
\Eq{htilde} 
$$\htilde(t) = \mcases{h(t),& t > 0,\\0&t < 0.}$$ 

On the other hand, the rescaled function $f(t) = h(t/c)$
satisfies Lamb's oscillator equation, \crf{Lambp; equation (8)}:
\Eq f
$$\req{f'' + \fra b\, f' + \frac{\sig ^2}{c^2} f = 0, \\ f(0) = 0,}$$
where $b = c/(2\:\beta )$.
Thus,
\Eq{Lambf}
$$\req{f(t) = C\, e^{-\:t/(2\:b)} \sin \kappa \:t,}$$
where $C$ is the integration constant, while
\Eq{Lambkappa}
$$\req{\kappa = \sqrt{\frac{\sig ^2}{c^2} - \fra{4\:b^2} } = \frac \varsigma c,\\ \varsigma = \sqrt{\sig ^2 - \beta ^2} \,,}$$
where we assume, as does Lamb, that the mass oscillator \eq f is underdamped, so
$$\sig > \beta > 0, \roq{or, equivalently,} c < 2\: b\: \sig.$$
Thus, setting 
\Eq h
$$h(t) = f(c\:t) = C\, e^{-\:\beta \: t} \sin \varsigma \, t,$$
reduces \eqr{sol}{htilde} to Lamb's solution, \rf{Lambp; equation (11)}:
\Eq{Lambsol}
$$u(t,x) = \mcases{-\:C\, e^{-\:(c\:t - \abs x)/(2\:b)} \sin \kappa \, (c\:t - \abs x),& \abs x < c\:t,\\0&\abs x > c\:t.}$$ 

The resulting solution profile is continuous and piecewise smooth, and consists of a symmetric pair of successively damped oscillatory disturbances generated by the mass that propagate into the undisturbed region at velocities $\pm\:c$. The only discontinuities in its spatial derivative occur at the front of the disturbance and at the origin.  The following plots graphs the solution spatial profile at a few times on the interval $-6\pii\leq x \leq 6\pii$ with vertical scale $-.45 \leq u \leq .45 $, for the specific values $\meq{c = 1,\\C=-1/2,\\b=5,\\\kappa=\sqrt{.99}\>\ustrut{11}}$. The damping effect of the Lamb coupling is already evident in the waves generated by the oscillator.

\Fig1{The Lamb Oscillator on the Line}{\vskip15pt
\figwidth=1.7in\figboxl=3in \figlabl=3.1in \hglue-40pt
\hbox{\Figboxl{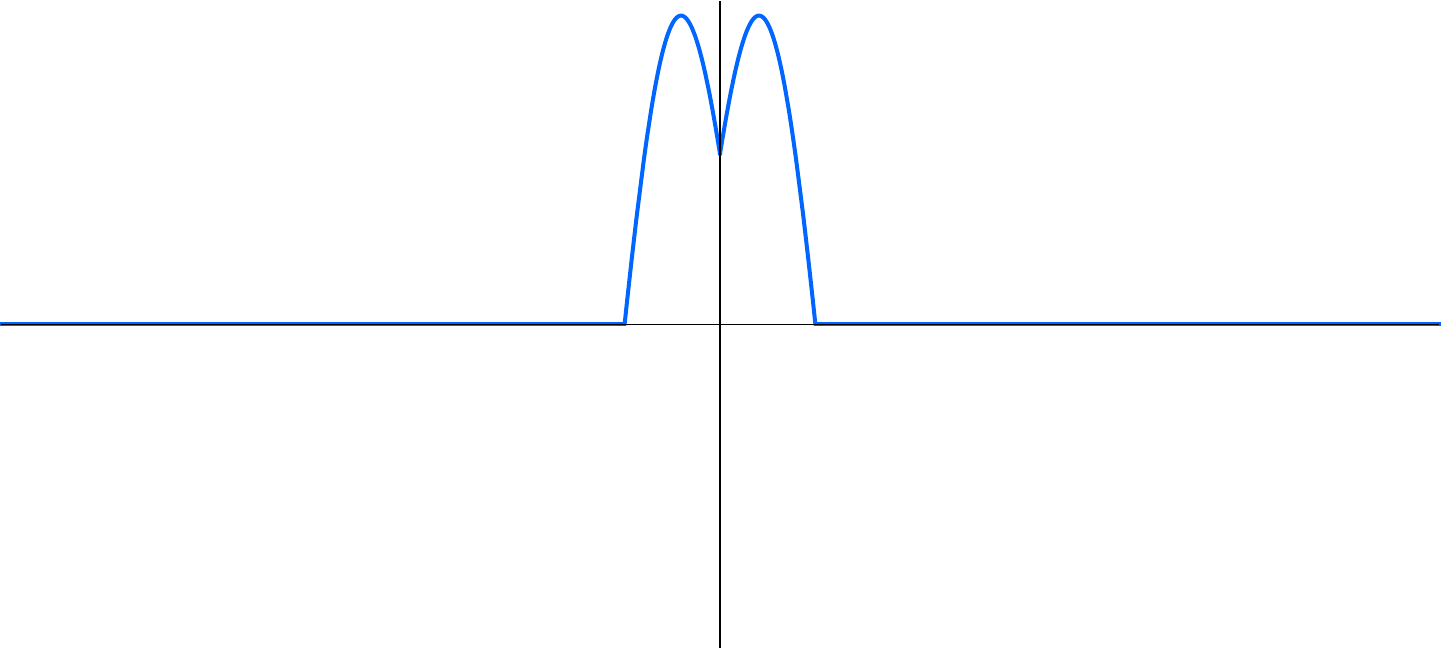}{$t=2.5$} \hskip -80pt\Figboxl{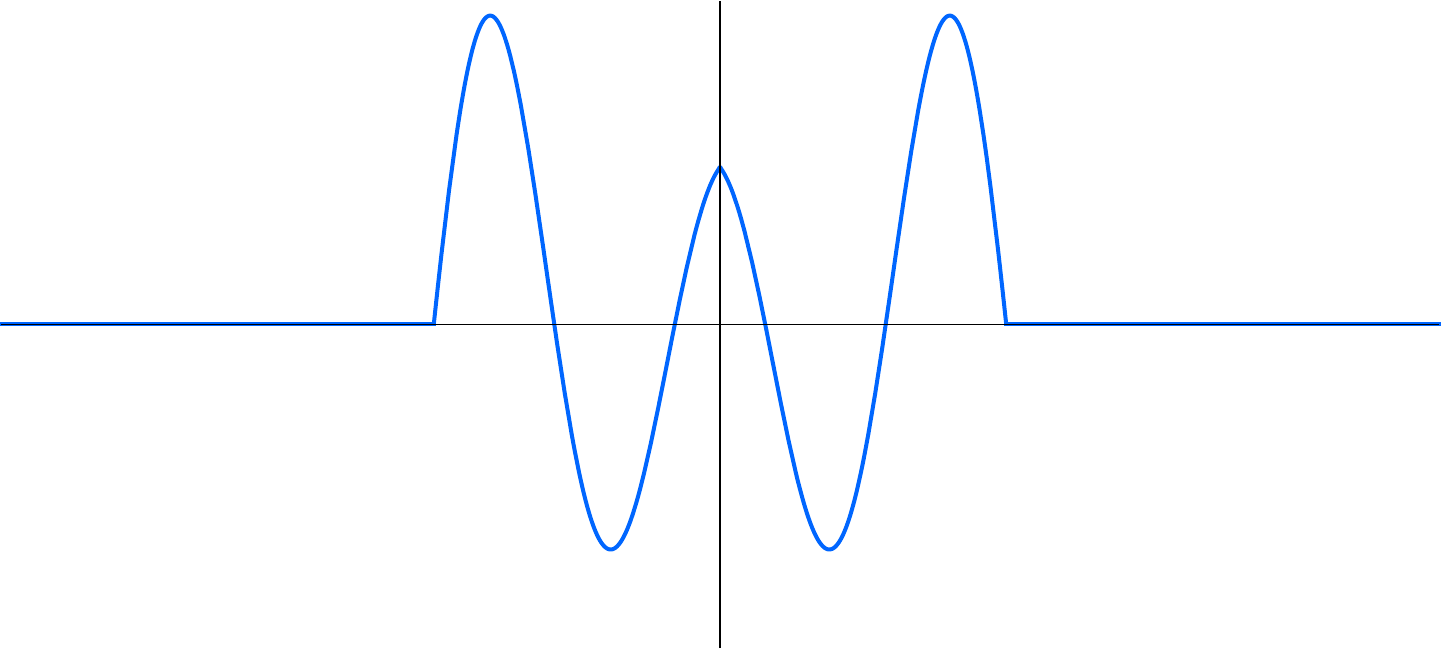}{$t=7.5$}\hskip -80pt\Figboxl{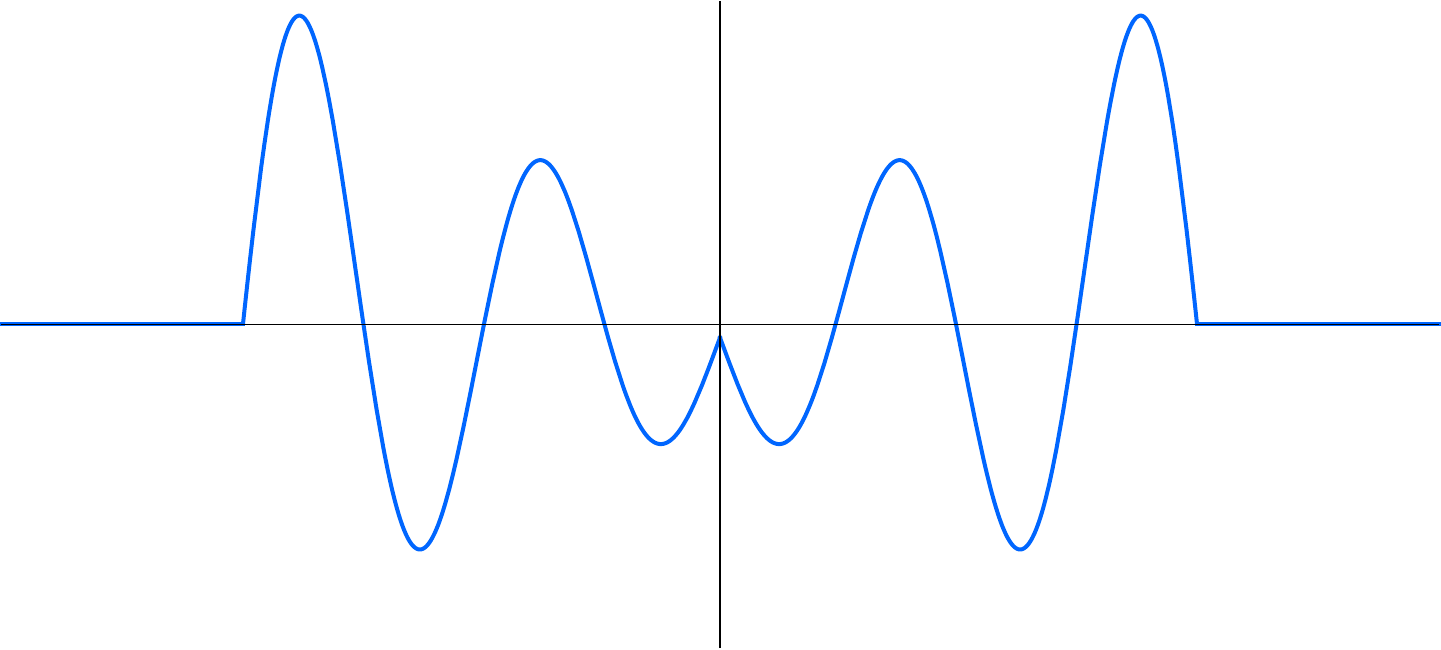}{$t=12.5$}}}
\vskip5pt

\noindent As time progresses, the initial disturbance continues to propagate outwards, while the damping completely dominates on an increasingly large interval centered around the origin, as shown in the following graph, which is plotted with the same vertical scale on the  larger interval $-50\pii\leq x \leq 50\pii$.

\vskip5pt

\Fig2{The Lamb Oscillator on the Line at Large Time}{\figwidth=5.4in\vskip15pt\hglue10pt\hbox{\figlabl=6in\Figboxl{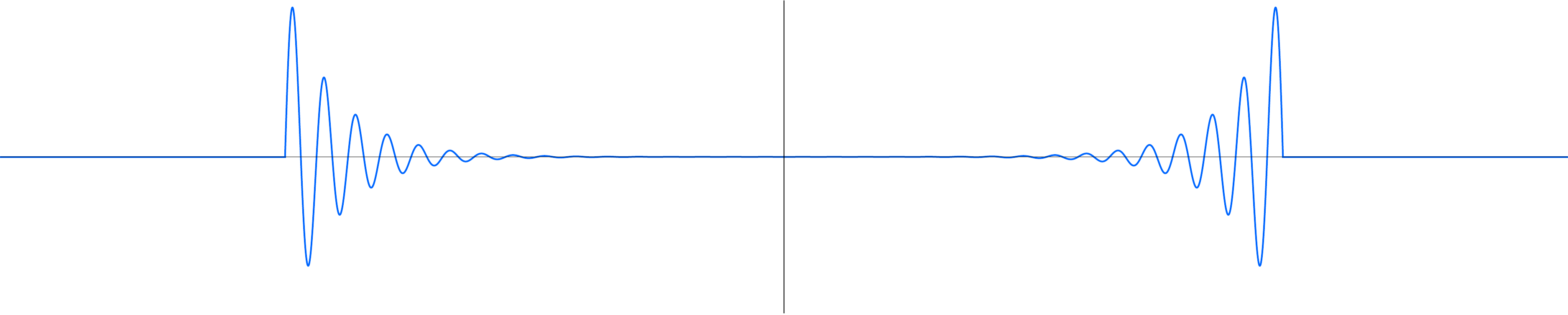}{$t=100$}}}
\vskip5pt


\Section p The Periodic Lamb Model.

Let us now focus our attention on the periodic version of the classical Lamb oscillator, that is equations \eqr{fwave}{mass} 
on the interval $-\pii < x < \pi$ with periodic boundary conditions. (Equivalently, we can regard $x$ as an angular coordinate on the unit circle.) As before, we assume the system is initially at rest, and so the initial conditions are \eq{ic0}. Thus, if we view the forcing term in \eq{fwave} as a $2\pii$ periodically extended delta function, the solution can be written in d'Alembert form as a formally infinite sum
\Eq{perLambsol}
$$\req{u(t,x) = \Sumii n \htilde(t - \abs {x- 2\:n\pii}/c),\\t > 0.}$$
Since $\htilde(t) = 0$ whenever $t \leq 0$, at any given point in space-time, only finitely many summands in \eq{perLambsol} are non-zero, and hence, by the preceding remarks, the solution is continuous and piecewise smooth.

Alternatively, one can solve the periodic Lamb problem via a Fourier series approach. Since the solution is even in $x$, we can expand it into a cosine series: 
\Eq{uF}
$$u(t,x) = \f2\:a_0(t) +\ \Sumi k a_k(t) \cos k\:x.$$
Substituting~\eq{uF} into \eqc{fwave}{ic0}, and using the fact that the delta function has the (weakly convergent) Fourier expansion
$$\delta(x)\sim \frac{1}{2\pi}+\frac{1}{\pi}\sum_{k=1}^\infty \cos kx,$$
we deduce that the Fourier coefficients $a_k(t)$ must satisfy the following decoupled system of initial value problems:
\Eq{aeq} 
$$\req{a_k'' + c^2k^2 a_k = -\frac{2\:c}\pi\:h'(t),\\a_k(0) = a_k'(0) = 0.}$$
Given the particular Lamb forcing function \eq h, the solution is straightforwardly constructed modulo some tedious algebra:
\Eq{LambFak}
$$a_k(t) = -\frac{2\: c \:C}{\pi}\bbk{p_k \cos c\:k \:t + q_k \sin c\:k \:t -e^{-\:\beta \: t} (p_k \cos \varsigma \:t + r_k \sin \varsigma \:t )} ,$$
where 
\Eq{pqr}
$$\ceq{\req{p_k = \frac{\varsigma\:(\sig ^2 -c ^2k^2)}{(\sigma^2 - c^2k^2)^2 + 4\:c^2k^2\beta ^2}\,,&
q_k = \frac{2\: c\:k\:\beta\:\varsigma}{(\sigma^2 - c^2k^2)^2 + 4\:c^2k^2\beta ^2}\,,}\\
r_k = \frac{\beta\:(\sig ^2 + c^2k^2)}{(\sigma^2 - c^2k^2)^2 + 4\:c^2k^2\beta ^2}\,,}$$
where, by \eq{Lambkappa}, $\sigma ^2 = \varsigma^2 + \beta ^2$.

Some representative solution profiles over a single period appear in the following plots. Keep in mind that the solution is both even and periodic.

\Fig p{The Periodic Lamb Oscillator}
{\prefigskip=5pt\figlabl=4.75cm\figwidth=4.5cm
\vskip-10pt
\Fiiil{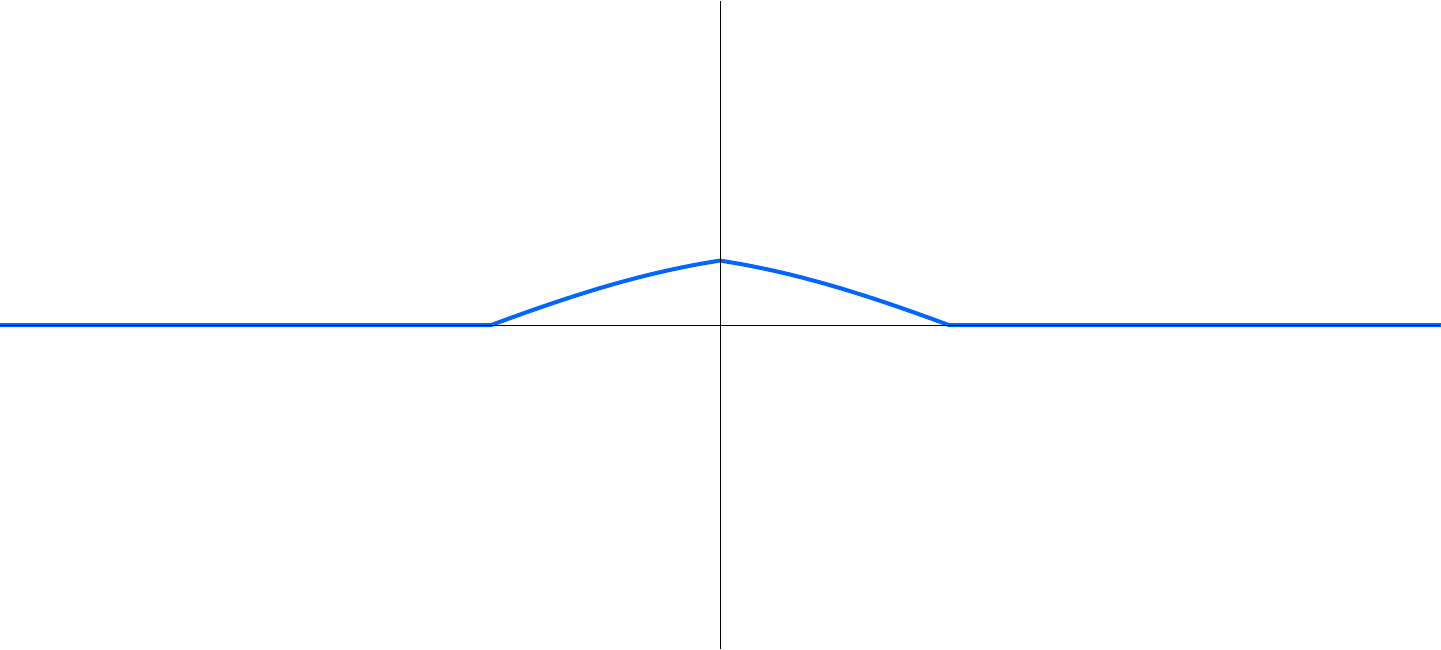}{1.}{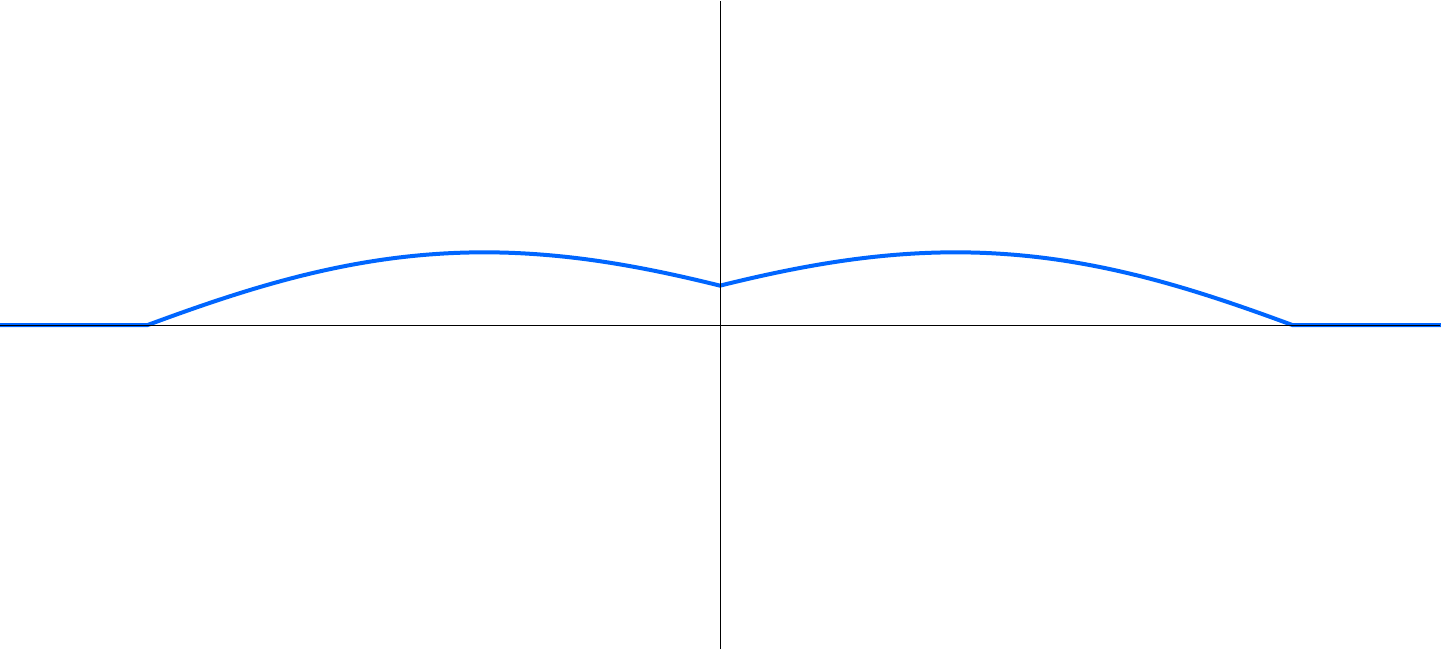}{2.5}{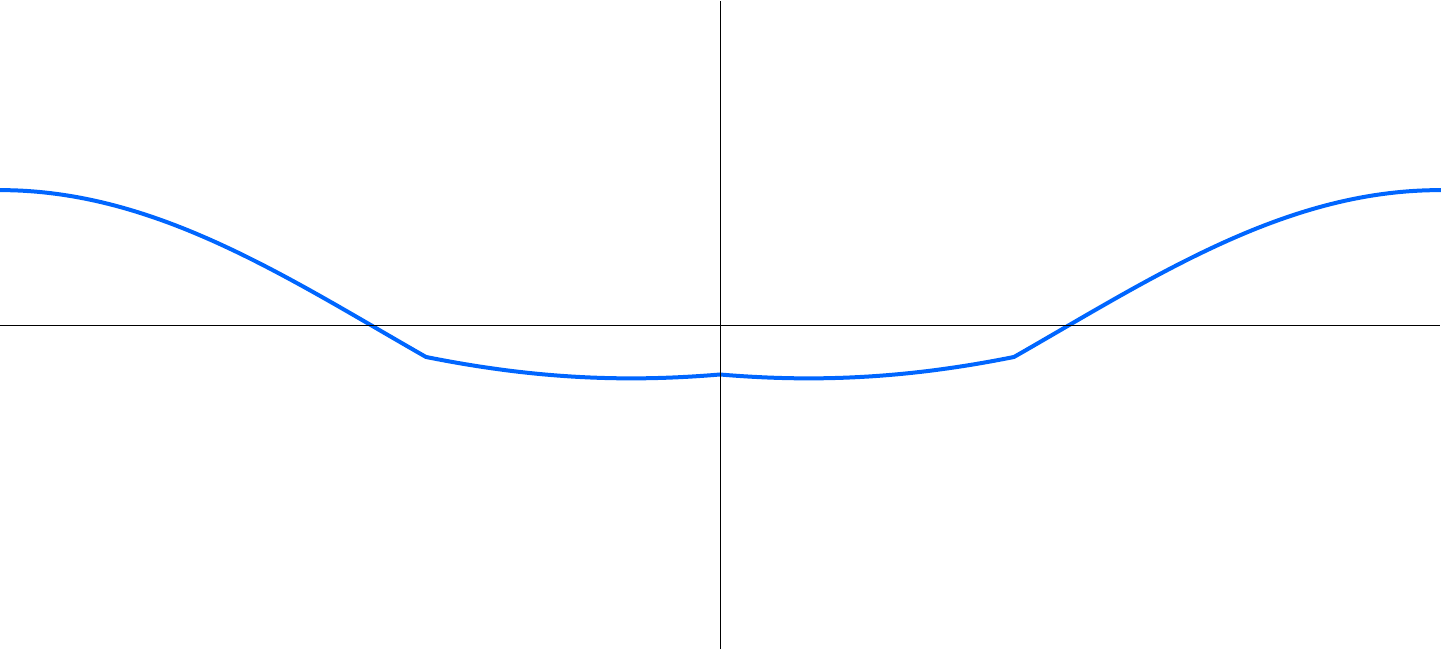}{5.}\vskip-10pt
\Fiiil{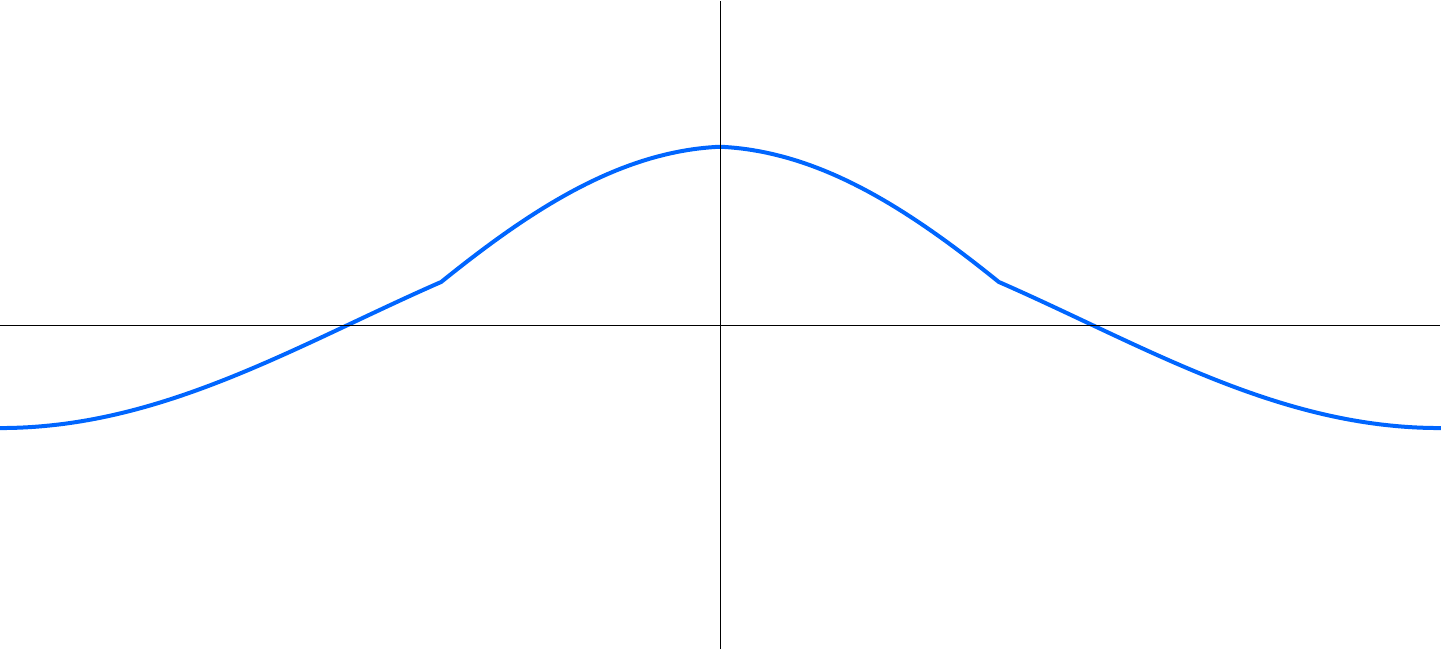}{7.5}{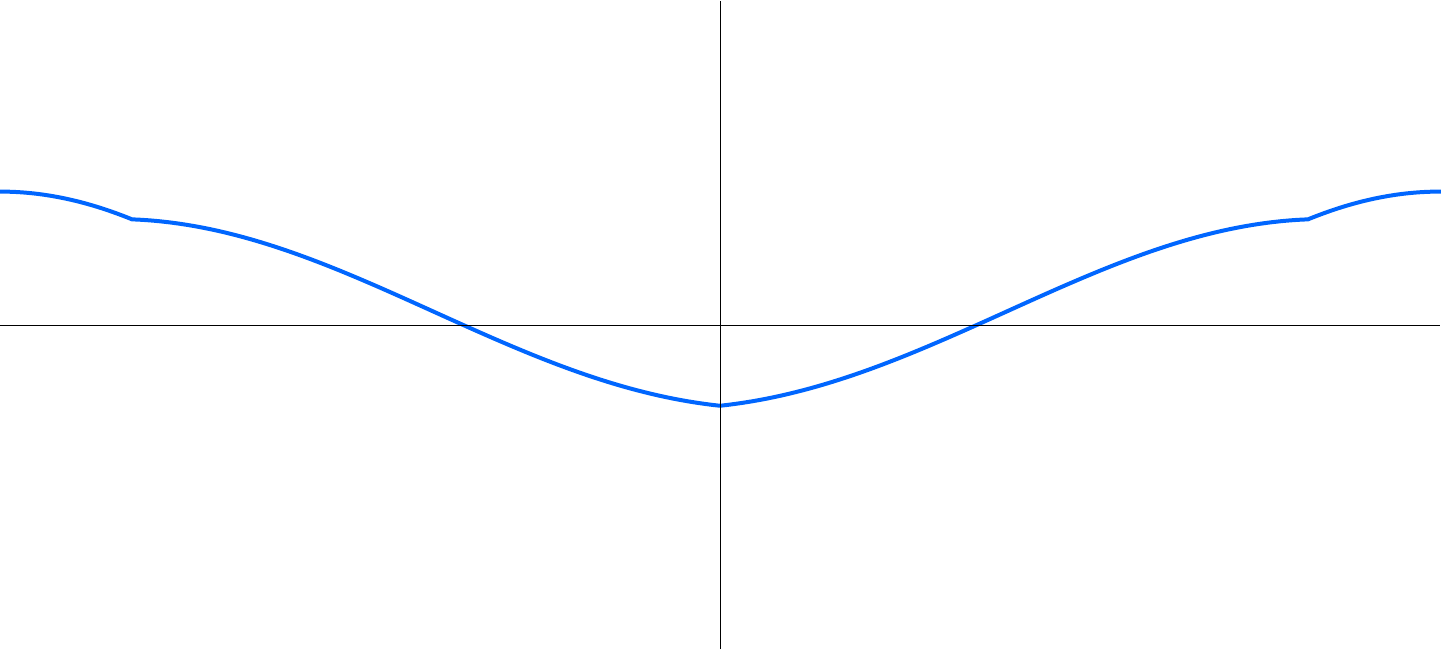}{10.}{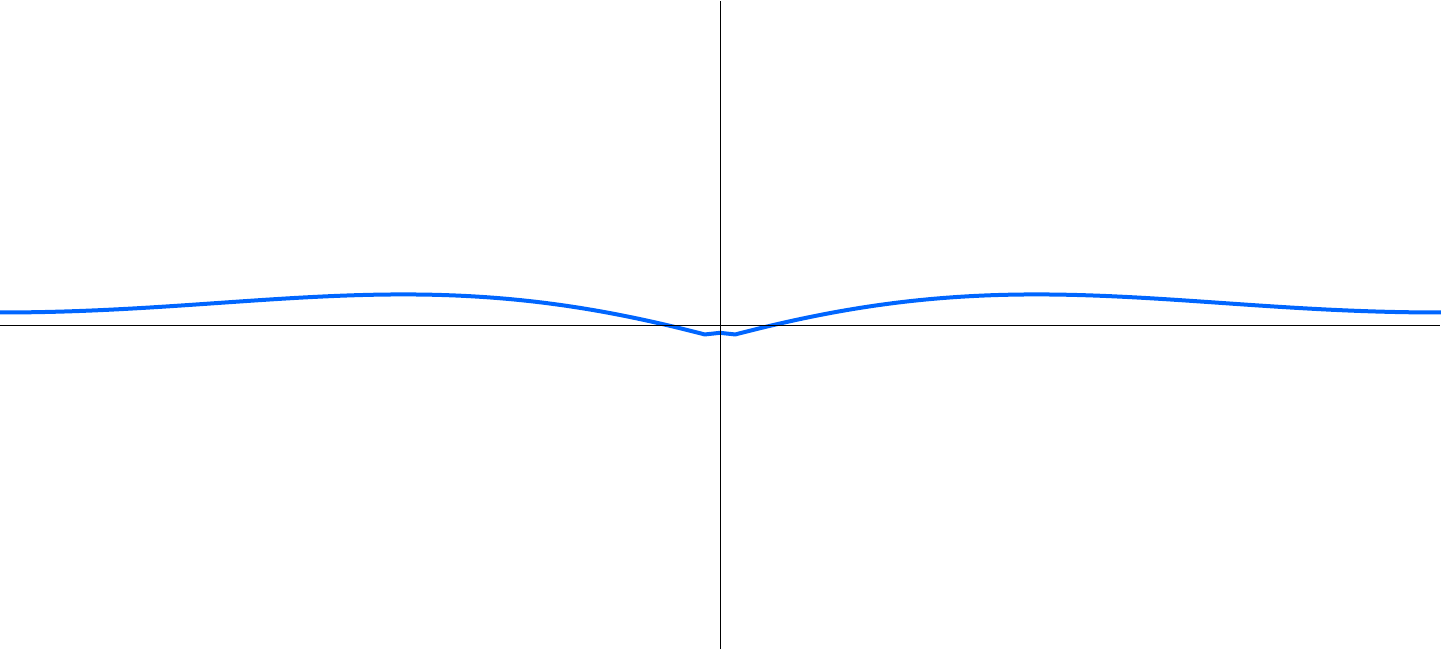}{12.5}}


\Section d Bidirectional Dispersive Lamb Models.

Let us now generalize the dispersion relation for the linear system being driven by a Lamb oscillator. We replace the second order wave equation by a second order (in time) dispersive linear (integro-)differential equation coupled to an oscillator:
\Eq{gdl}
$$u_{tt} = L\br u + h'(t) \, \delta (x),$$
where $h(t)$ is as above, \cf \eq h. (For simplicity, we have absorbed the prefactor $-2\:c$ in \eq{fwave} into the forcing function: $h(t) \longmapsto -\:h(t)/(2\:c)$.)
The linear operator $L$ is assumed to have dispersion relation $\omega = \omega (k)$ relating temporal frequency $\omega $ and wave number $k$. As usual, \rf{Whitham}, $\omega (k)$ is found by substituting the exponential ansatz $e^{\i(k \:x - \omega \:t) }$ into the unforced equation.  We assume that the system is purely dispersive, which is equivalent to requiring that $\omega (k)$ be a real-valued function of the wave number $k \in \R$. The operator $L$ is, in general, of integro-differential form, and is a differential operator if and only if $\omega (k)$ is a polynomial. 

In particular, the wave equation \eq{fwave} has the linear dispersion relation $\omega = c\:k$. A higher-order correction modeling a linearly elastic string takes the form, \rf{Kunin, Weinberger},
\Eq{fB}
$$u_{tt} = c^2 u_{xx} - \varepsilon \:u_{xxxx} + h'(t) \, \delta (x),$$
whose dispersion relation
\Eq{fBd}
$$\omega = \sqrt{c^2 k^2 + \varepsilon\: k^4}\,,$$
is asymptotically quadratic for large wave number: $\omega (k) \sim k^2$ as $\abs k \to \infty $. The same dispersion relation arises in the Boussinesq approximations to the free boundary problem for water waves, \rf{Whitham}. A regularized version
\Eq{rB}
$$u_{tt} = c^2 u_{xx} + \varepsilon \:u_{xxtt} + h'(t) \, \delta (x),$$
in which two of the $x$ derivatives are replace by $t$ derivatives, has the same order of approximation to the full physical model, and has also been proposed as the linearization of a model for DNA dynamics, \rf{Scott}. \eqE{rB} has dispersion relation
\Eq{rBd}
$$\omega = \frac{c\,\abs k}{\sqrt{1 + \varepsilon\: k^2}}\,,$$
which is asymptotically constant at large wave number.
In their work on gyroscopic Lamb systems, \rf{HaBlWe}, Hagerty, Bloch, and Weinstein investigate coupling the Lamb oscillator to the linearly dispersive Klein--Gordon equation:
\Eq{KGL}
$$u_{tt} = c^2 u_{xx} - m^2 u + h'(t) \, \delta (x),$$
with ``mass'' $m$. In this case, the dispersion relation
\Eq{KGLd}
$$\omega = \sqrt{c^2 k^2 + m^2}\,,$$
is asymptotically linear: $\omega (k) \sim \abs k$. 

We begin by investigating the solution to \eq{gdl} on the line with zero initial condition $u(0,x) = 0$. 
While there is no longer a d'Alembert style formula for the solution, the Fourier transform will enable us to express the solution in the following form:
\Eq{biD}
$$ u(t,x)= \frac{C}{4\: c \pii}\int_{-\infty}^\infty e^{\i k\:x}\bbk{p_k \cos \omega(k)\: t+q_k \sin \omega(k)\: t-e^{-\beta t}(p_k \cos \varsigma\: t+r_k \sin \varsigma \:t ) }~dk,$$
where
\Eq{pqrw}
$$\xeq{p_k = \frac{\varsigma\:(\sig ^2 -\omega^2)}{(\sigma^2 - \omega ^2)^2 + 4\:\omega^2\beta ^2}\,,\\
~q_k = \frac{2\: \omega\:\beta\:\varsigma}{(\sigma^2 - \omega ^2)^2 + 4\:\omega^2\beta ^2}\,,\\
~r_k = \frac{\beta\:(\sig ^2 + \omega^2)}{(\sigma^2 - \omega ^2)^2 + 4\:\omega^2\beta ^2}\,.}$$
Note that \eq{pqrw} coincides with \eq{pqr} upon setting $\omega=c\:k$, which is the dispersion relation for the classical wave equation \eq{wave}. Thus, \eq{biD} and the Fourier transform of~\eq{Lambsol} are equivalent, modulo the factor $-2\:c$ which was absorbed into $h(t)$. 


Turning to the corresponding periodic problem for a dispersive Lamb system \eq{gdl} on the interval $-\pii < x < \pi$, we work with the Fourier series representation \eq{uF}.
The \eqse{aeq} for the Fourier coefficients become
\Eq{daeq} 
$$\sqeq{a_0'' = h'(t)/\pi\,,& a_0(0) = a_0'(0) = 0\,,\\
a_k'' + \omega(k)^2 a_k = h'(t)/\pi\,,& a_k(0) = a_k'(0) = 0.}$$
Thus, the solution to the periodic dispersive Lamb problem has 
\Eq{LambdFak}
$$\eeq{a_0(t) = \frac{C}{\pii(\beta^2+\varsigma^2)} \bbk{ \varsigma-e^{-\:\beta \: t}(\varsigma\cos \varsigma \:t+\beta \sin \varsigma \:t)},\\
a_k(t) = \frac{C}{\pii}\bbk{p_k \cos \omega(k) \:t + q_k \sin \omega(k) \:t - e^{-\:\beta \: t} (p_k \cos \sig \:t + r_k \sin \sig \:t )},}$$
for $k\geq1$ where $p_k,q_k,r_k$ are as in~\eq{pqrw}.

\Fig4{The Dispersive Periodic Lamb Oscillator with $\omega (k) = k^2$}
{\prefigskip=5pt\figlabl=4.75cm\figwidth=4.5cm
\vskip-10pt\hglue 30pt
\Fiiil{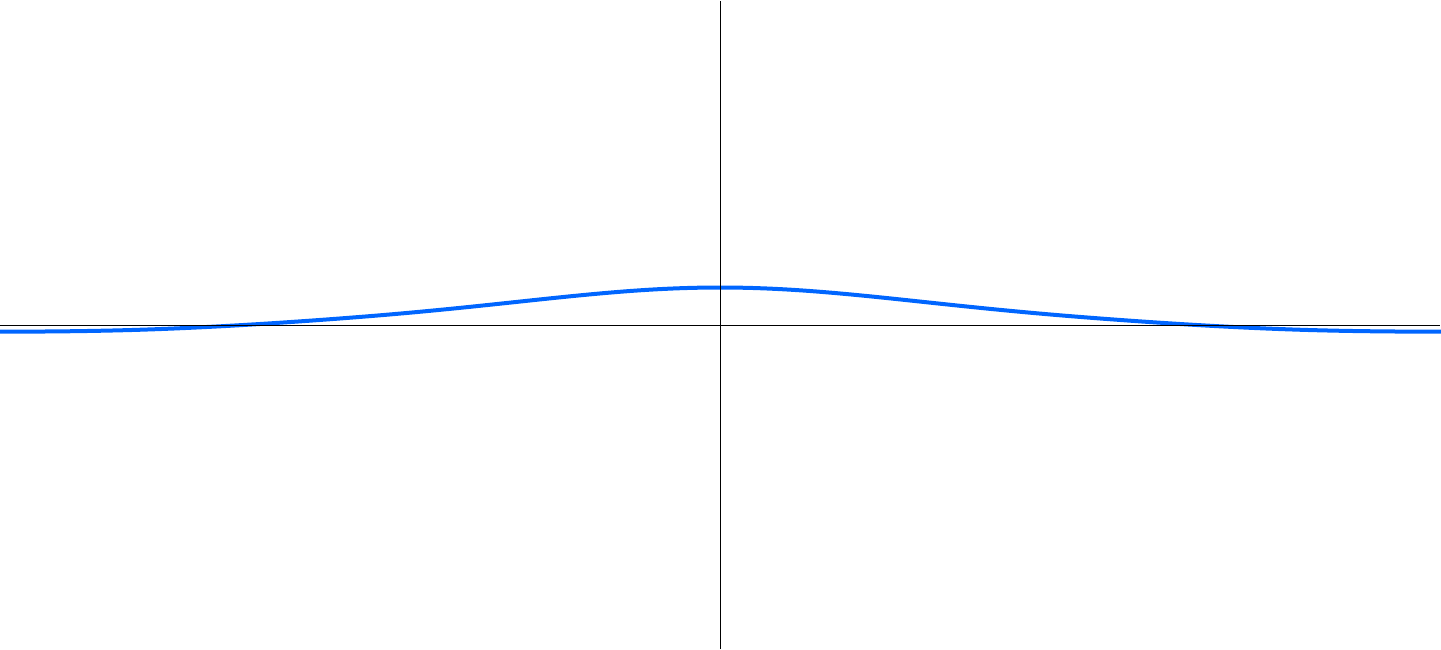}{1.}{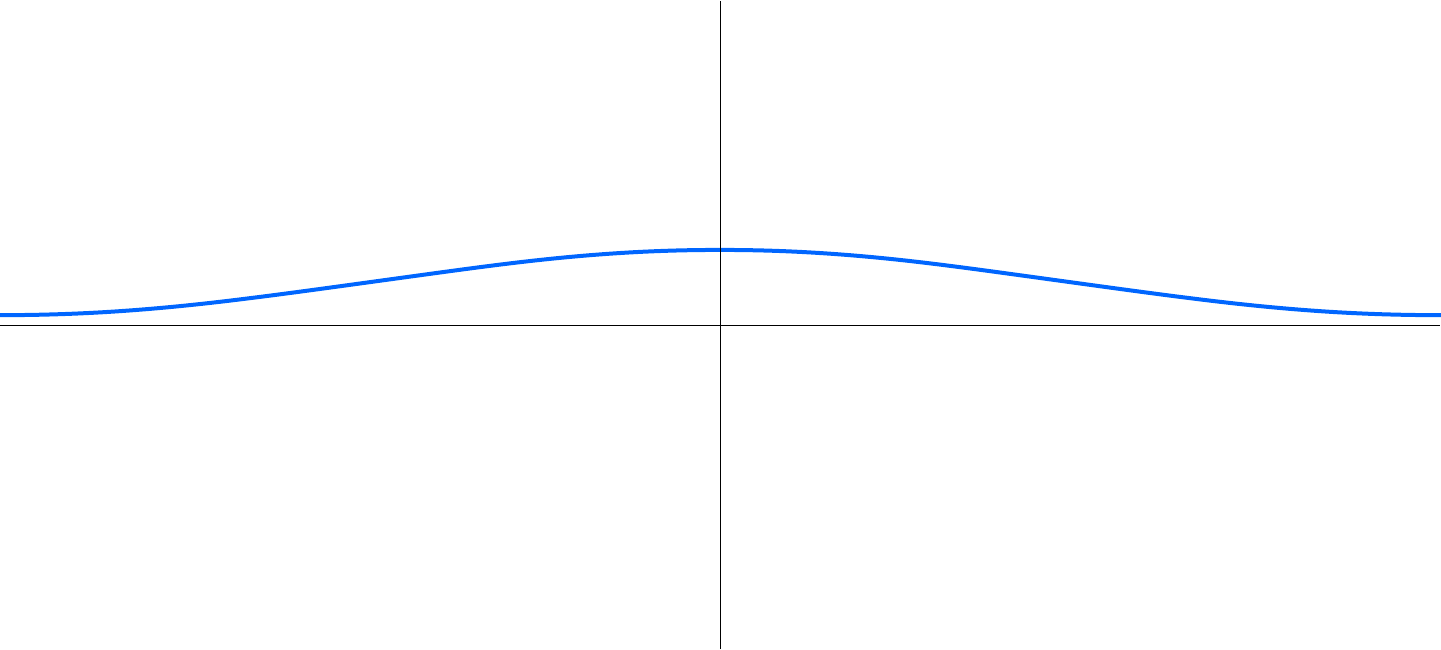}{2.5}{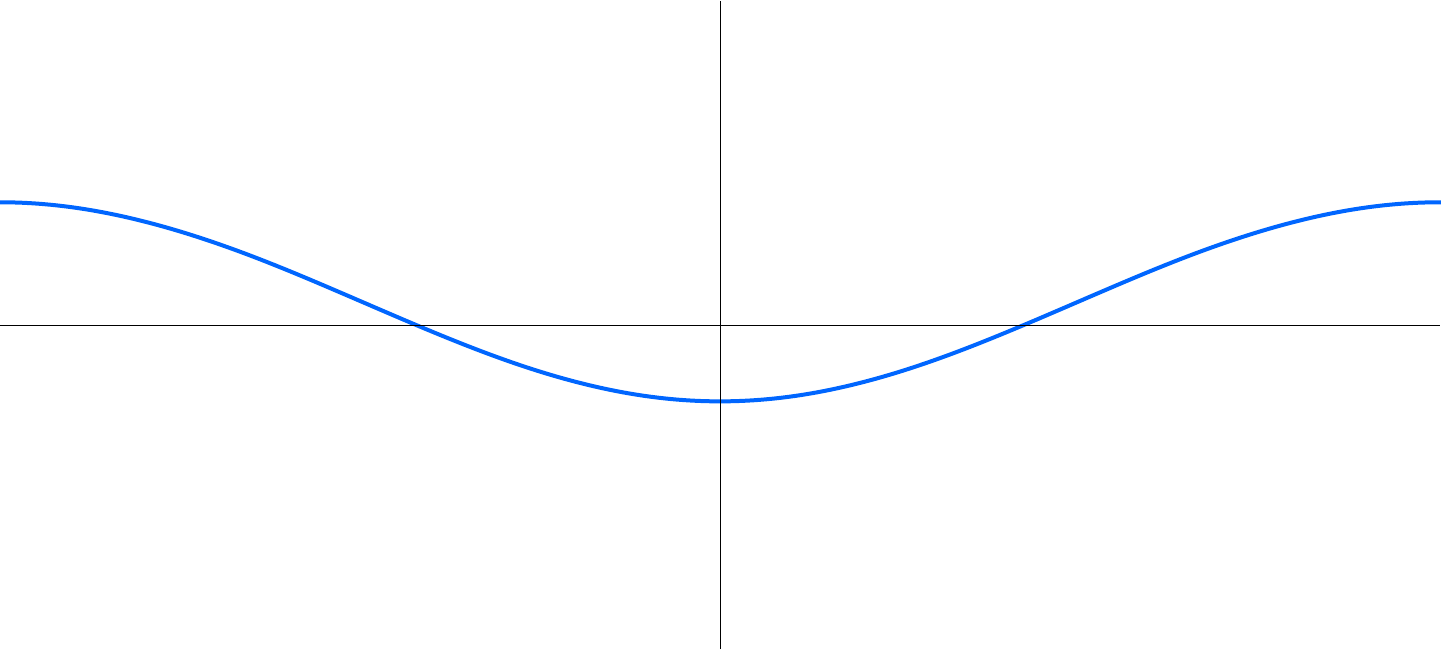}{5.}\vskip-10pt
\Fiiil{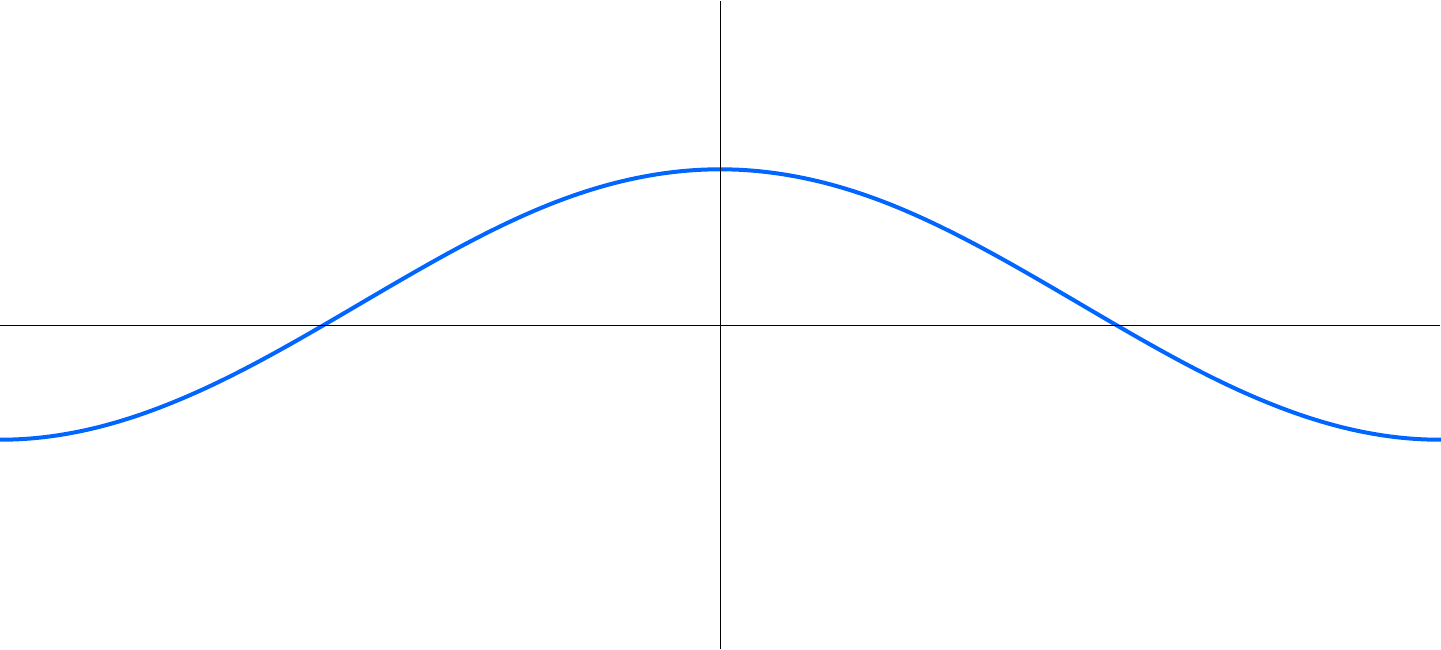}{7.5}{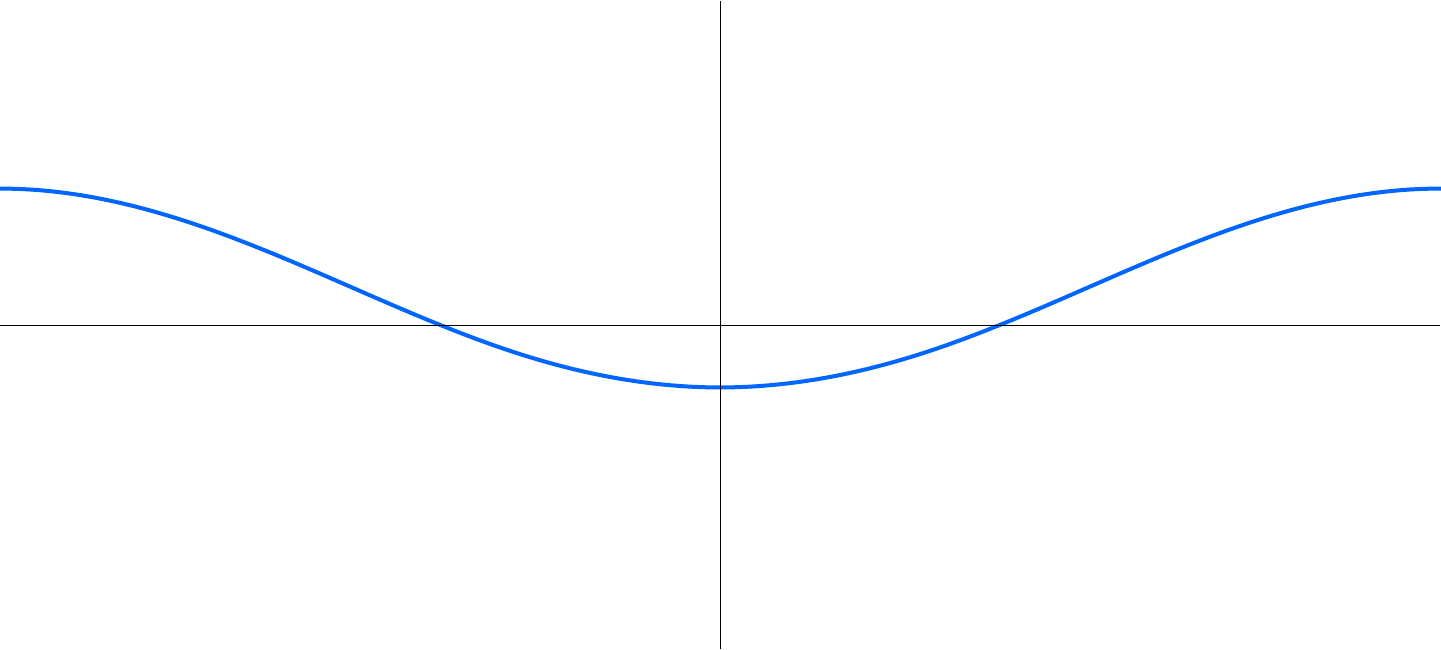}{10.}{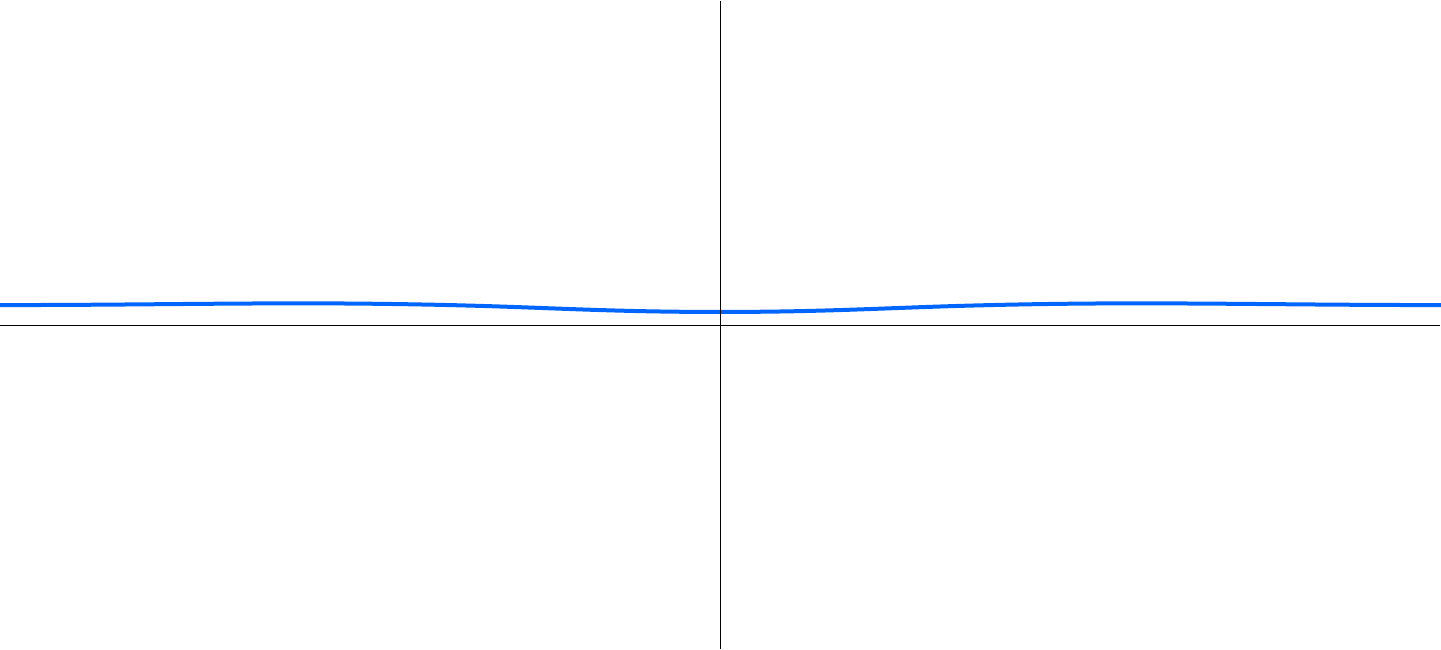}{12.5}\vskip10pt
\vskip-5pt}

\Fig5{The Dispersive Periodic Lamb Oscillator for the Klein--Gordon Model}
{\prefigskip=5pt\figlabl=4.75cm\figwidth=4.5cm
\Fiiil{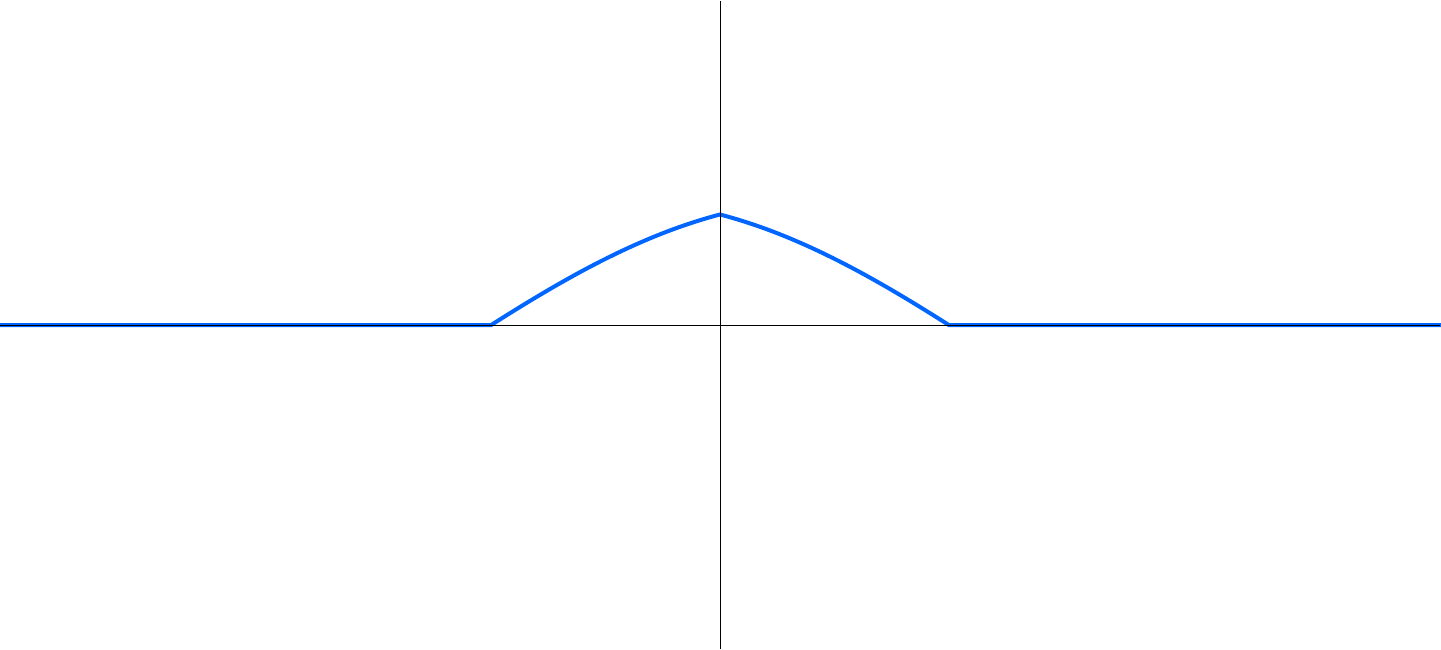}{1.}{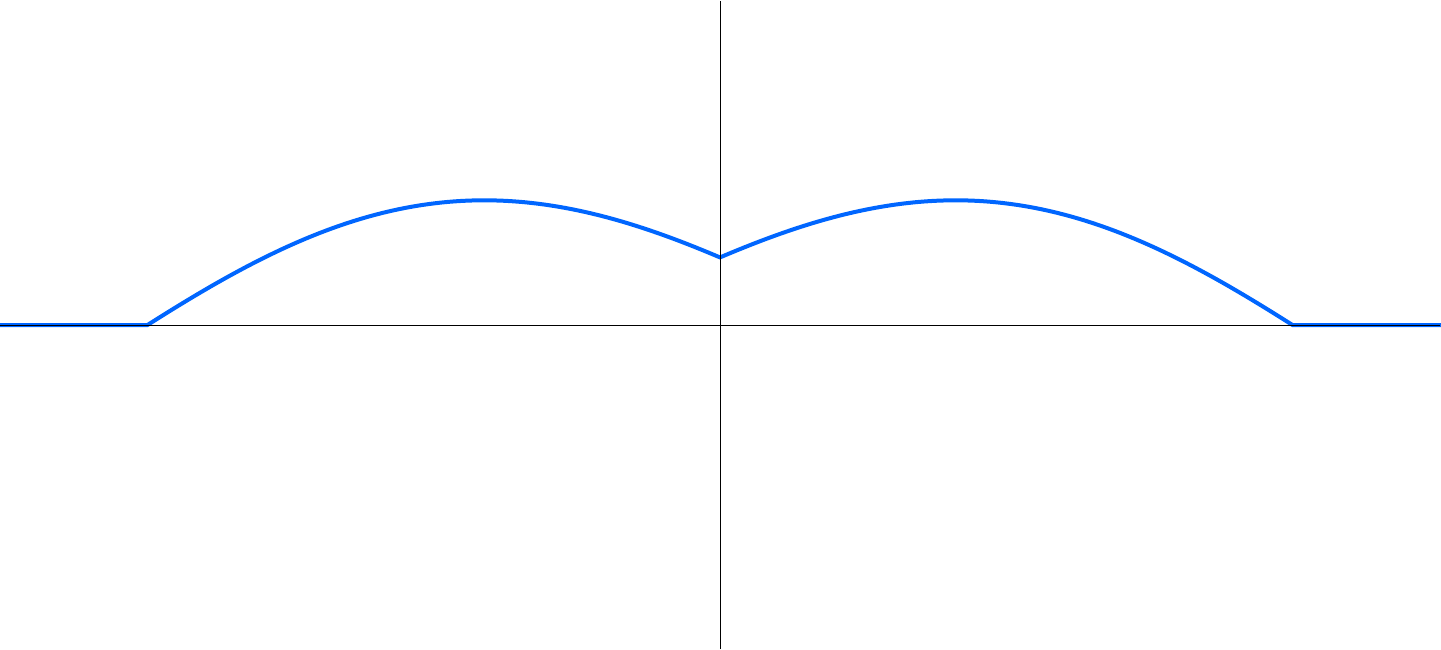}{2.5}{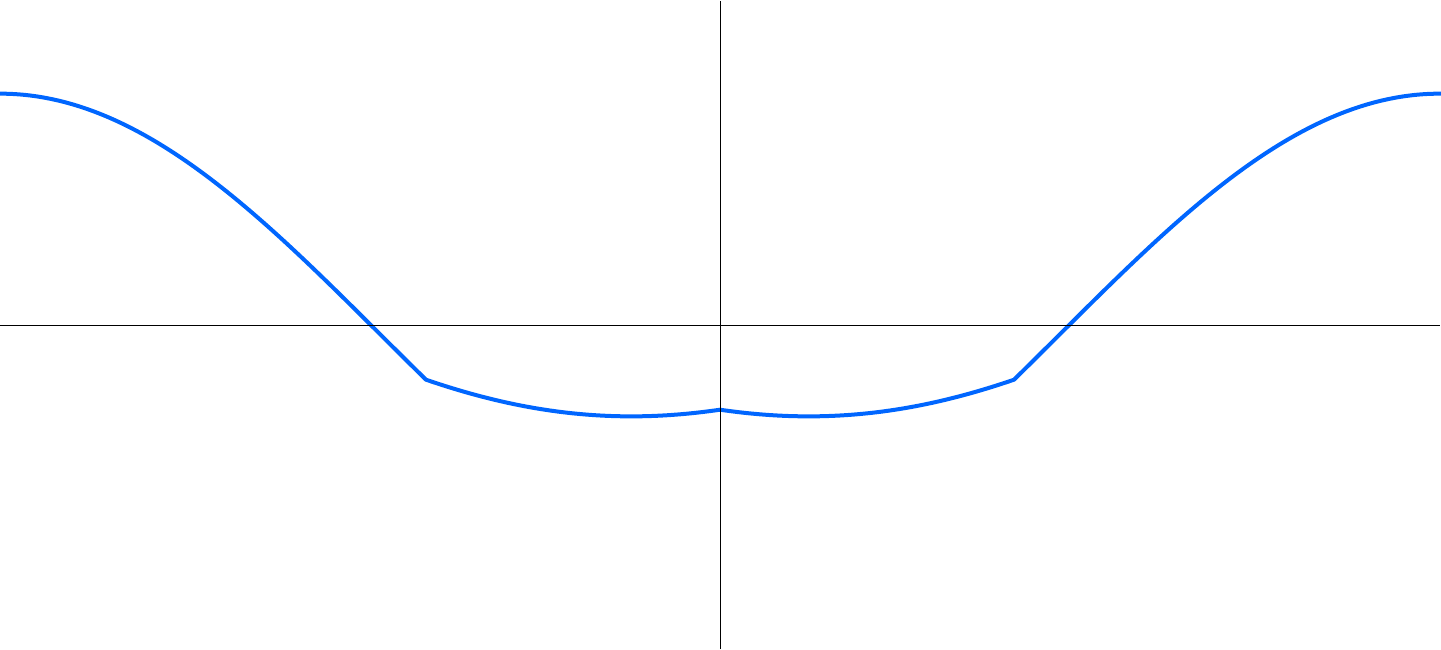}{5.}\vskip-10pt
\Fiiil{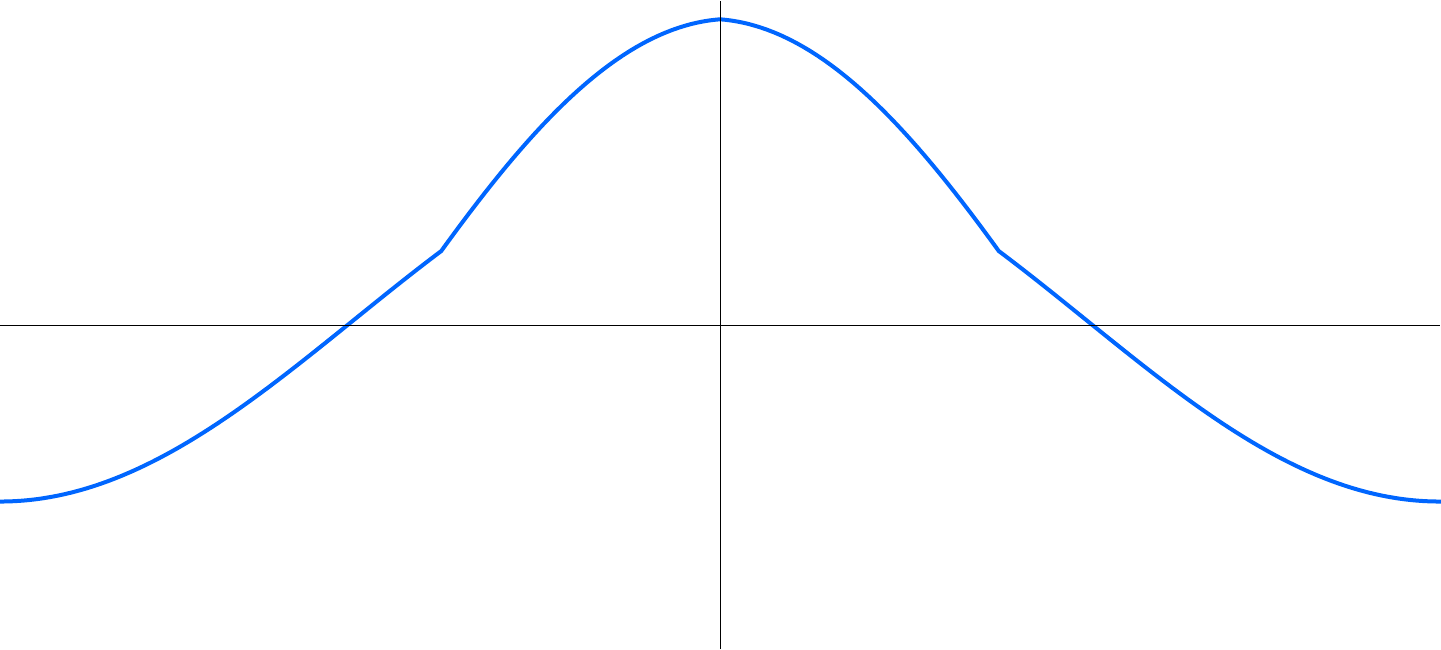}{7.5}{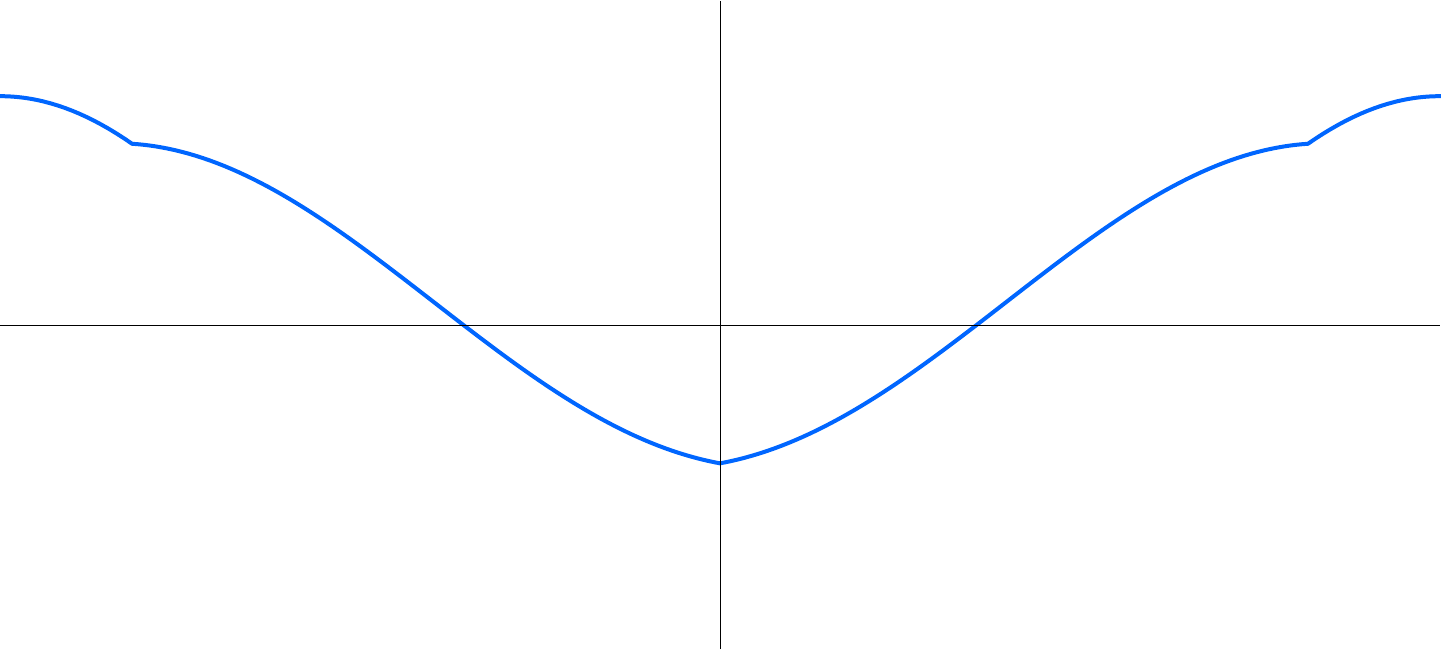}{10.}{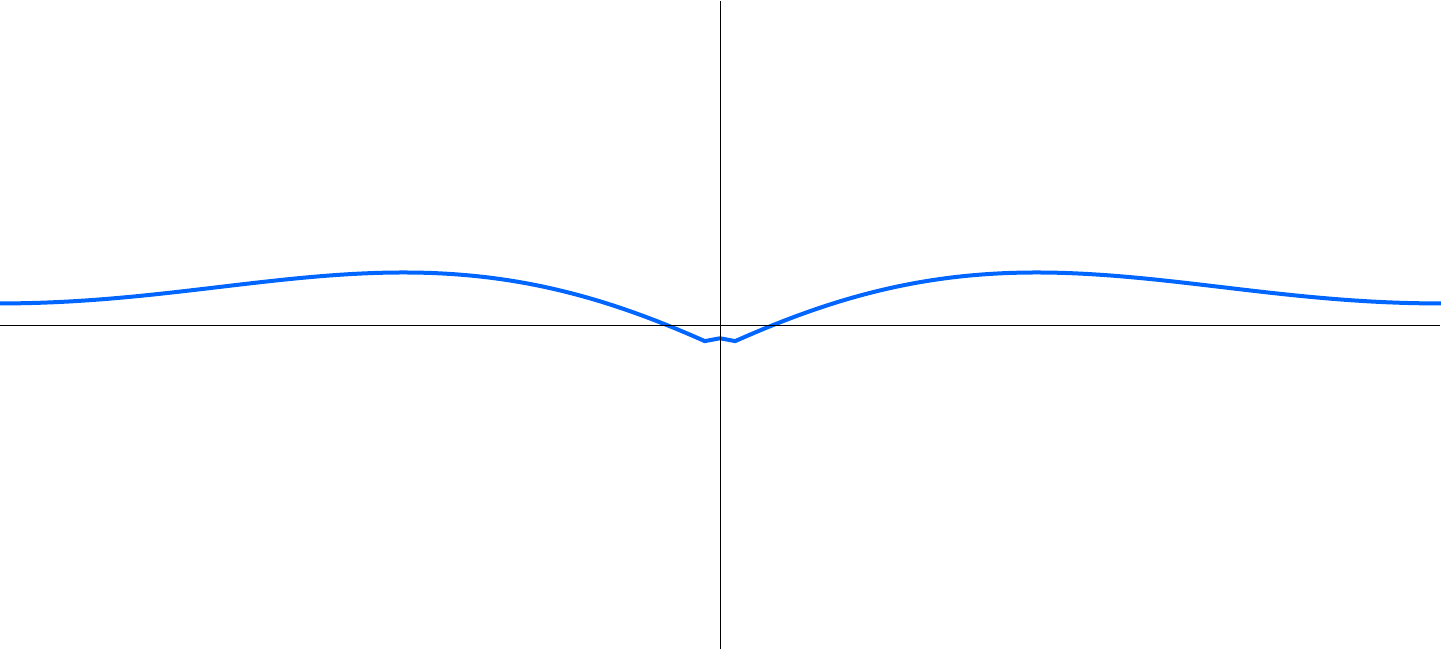}{12.5}\vskip10pt
\vskip-5pt}

In particular, suppose that the dispersion relation has large wave number asymptotics given by a power law:
\Eq{drm} 
$$\omega (k) \simx \abs k^m \roq{as} \abs k \longrightarrowx \infty ,$$
for some $m \in \R$.
Then, assuming $m > 0$, according to \eq{LambdFak},
\Eq{dak} 
$$a_k(t) \simx \omega(k)^{-2} \simx \abs k^{-\:2\:m} \roq{as} \abs k \longrightarrowx \infty .$$
The resulting Fourier coefficients therefore decay as $\abs k \to \infty$, which, depending upon the magnitude of $m$, implies smoothness, meaning differentiability, of the solution. Namely, when $m > 1$, \eq{dak} implies that $u(t,\cdot ) \in \Cn$ provided $2\:m-1 > n \in \N$, by a standard result in Fourier analysis, \rf P. Thus, in this situation, the solution to the periodic dispersive Lamb problem is continuously differentiable up to order $n$ and \iz{there is no fractalization}.
This observation is borne out by \Mathematica\ calculations, obtained by summing the first $1000$ terms of the Fourier series \eq{uF}; examples appear in the above plots for quadratic dispersion and the Klein--Gordon model \eq{KGL}. 


When the exponent $m > 1$, in direct contrast to the smoothness of the Lamb solution, starting with rough initial conditions, \eg a step function, fractalization and, when $m \in \N$, dispersive quantization at rational times will be manifested in the unforced periodic solution profiles, \rf{COdisp}.  The effect of combining such rough initial data with the Lamb oscillator is simply a linear superposition of the two solutions.



On the other hand, if the dispersion relation is sublinear, $m < 1$, then the resulting Lamb solution coefficients exhibit the slow decay that underlies the fractalization effects observed is dispersive quantization. The most interesting examples occur when $\omega(k) \sim \sqrt{\abs k}$, whereby \eq{dak} implies that $a_k(t) \sim 1/\abs k$, which is exactly the decay rate that produces the dispersive quantization effects.  In particular, $\omega(k) = \sqrt{\abs k}$ corresponds to the (complex) linear system
$$ u_{tt} = \i u_x,$$
\ie the linear Schr\"odinger equation with the roles of space and time reversed. Interestingly, the dispersion relation for the free boundary problem for water waves, $\ustrut{11}\omega (k) = \sqrt{k \tanh k}$ (ignoring physical constants), \rf{Whitham}, is also asymptotically of this form. 



\Fig6{The Dispersive Periodic Lamb Oscillator with $\omega(k) = \sqrt{\abs k}$}{\prefigskip=5pt\figlabl=4.75cm\figwidth=4.5cm
\Fiiil{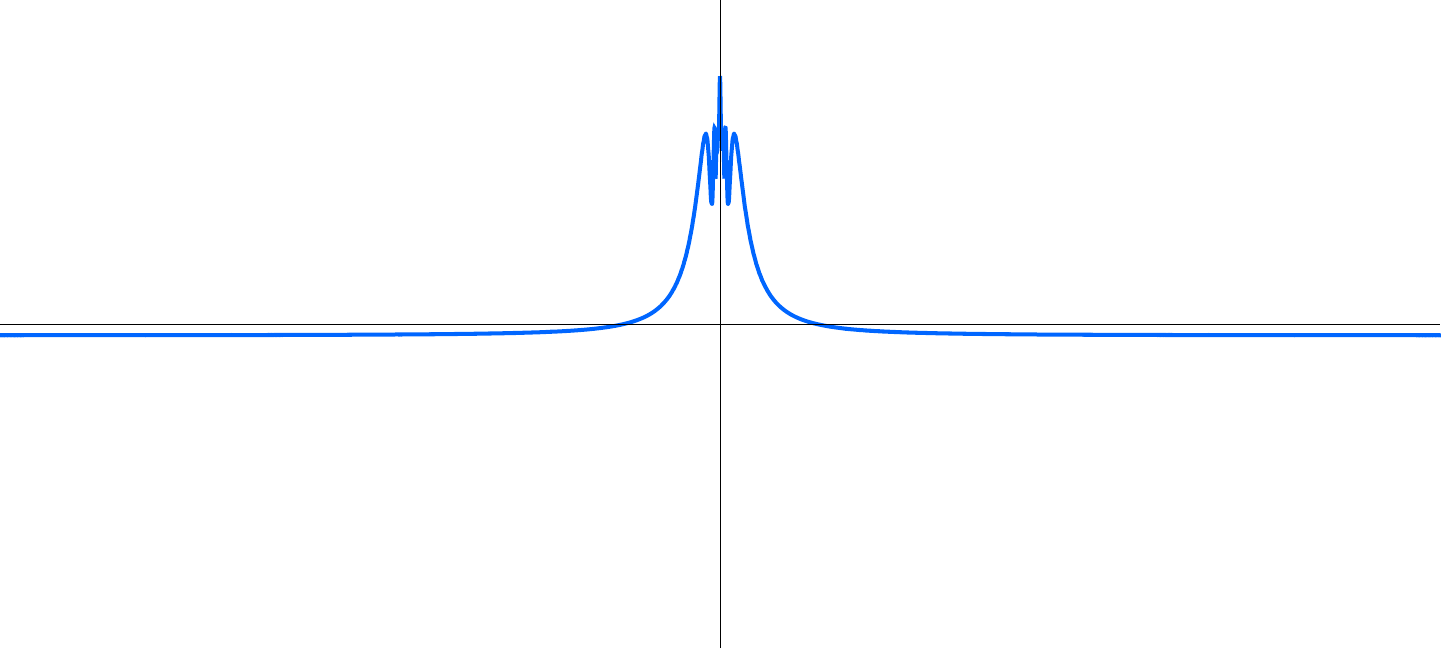}{1.}{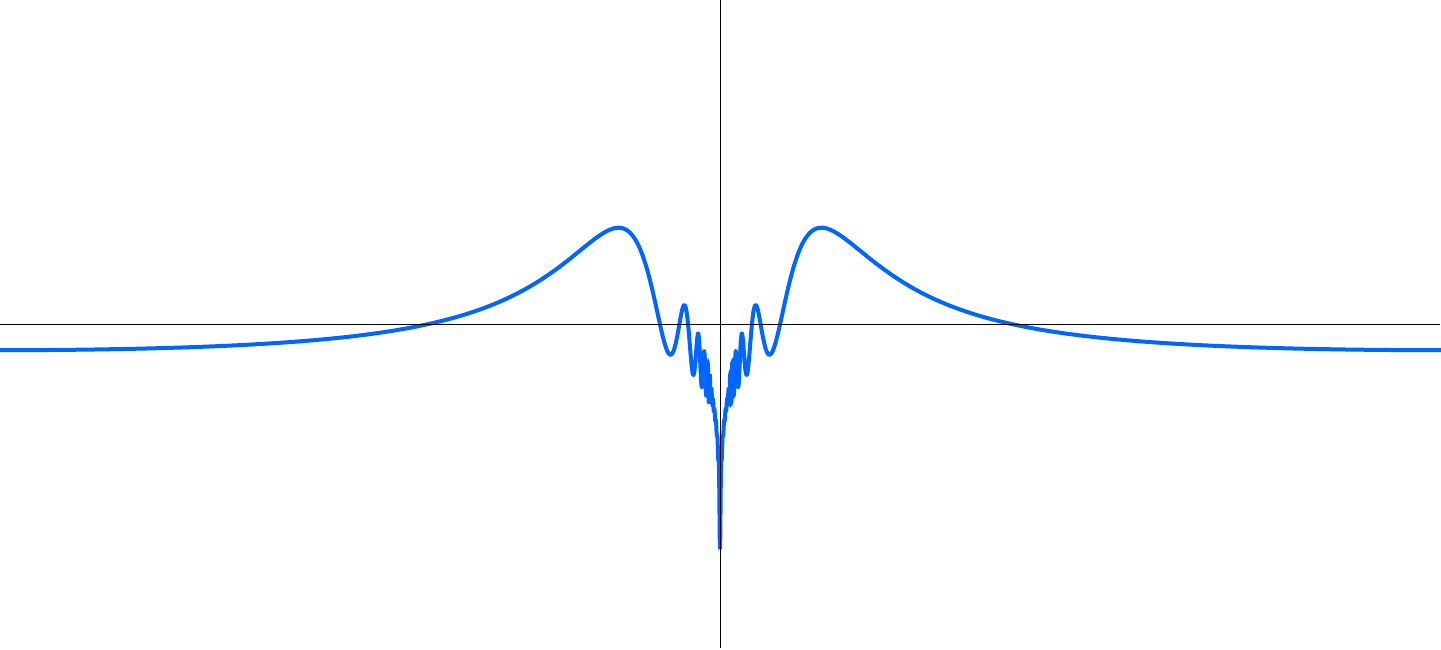}{2.5}{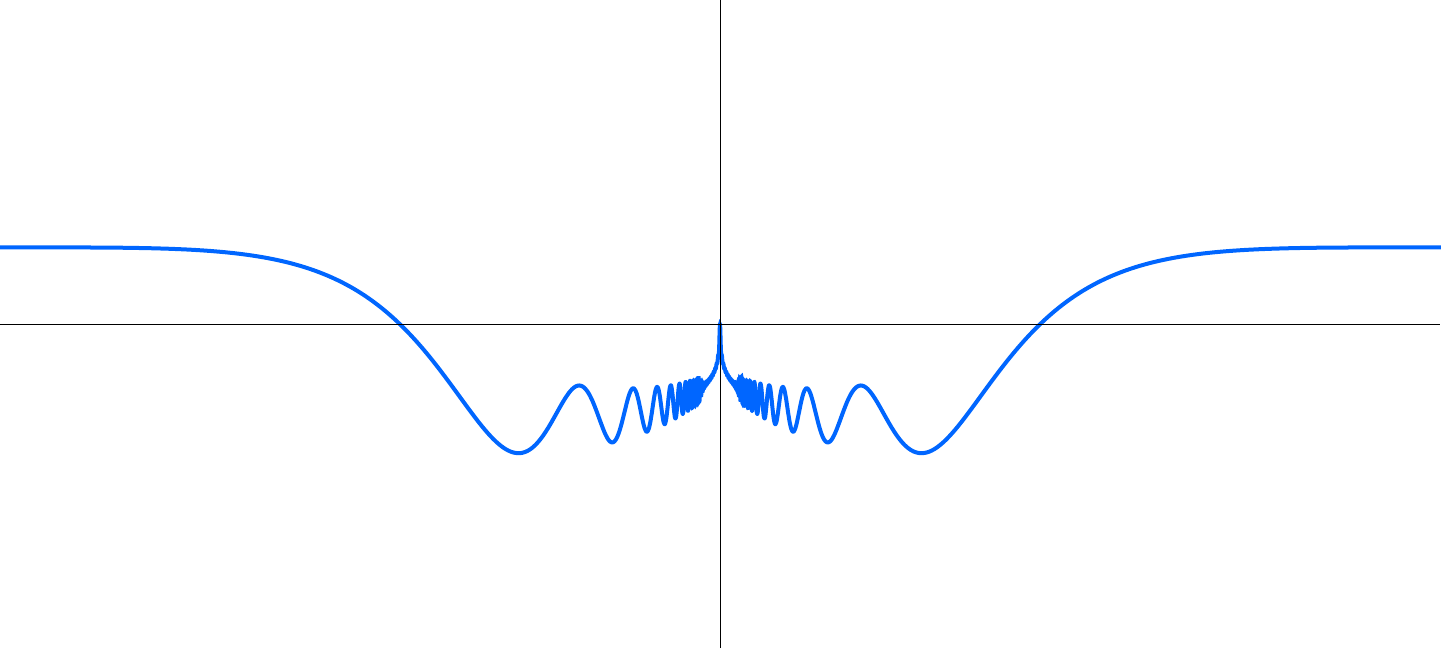}{5.}\vskip-20pt
\Fiiil{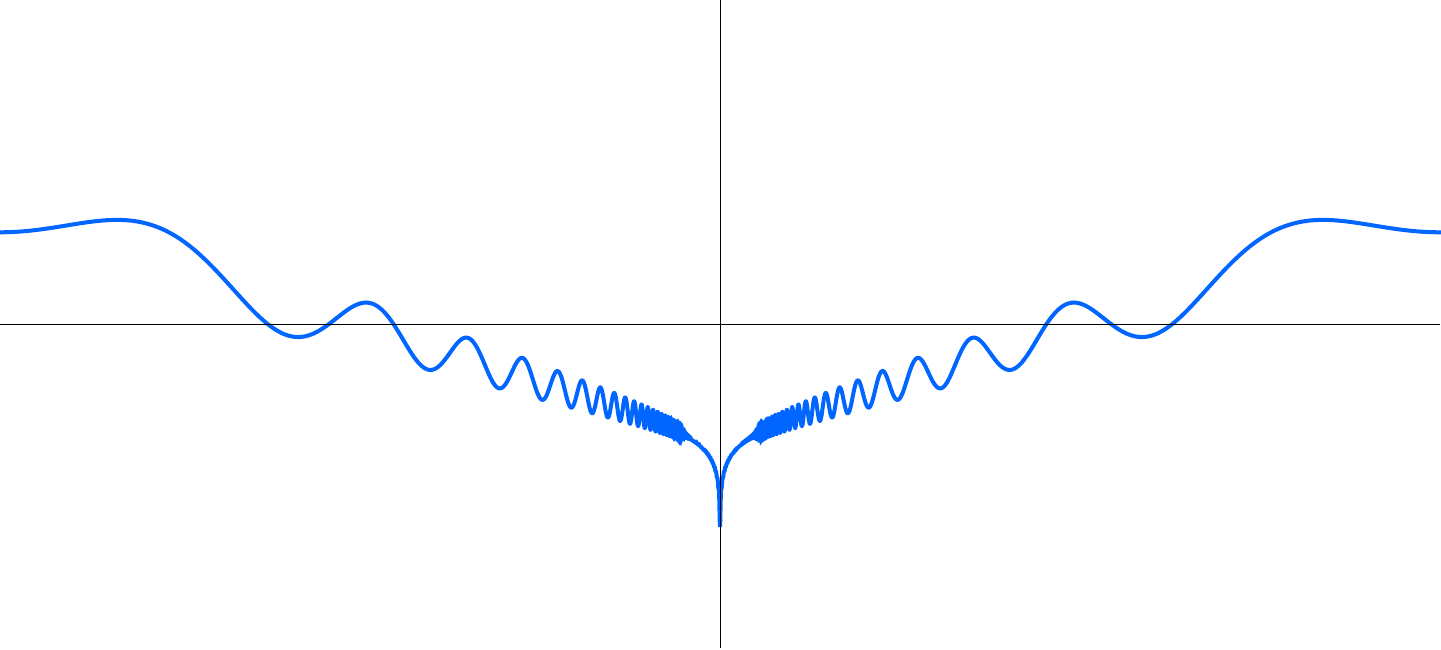}{10.}{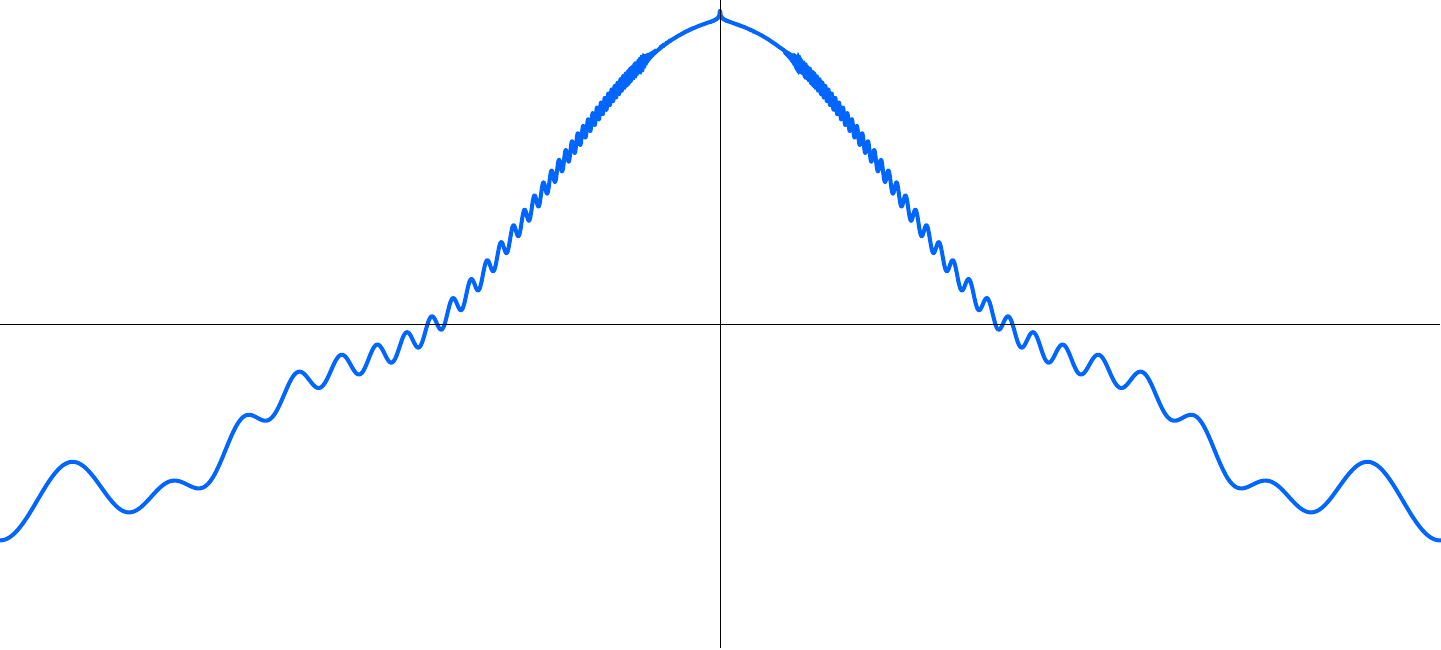}{20.}{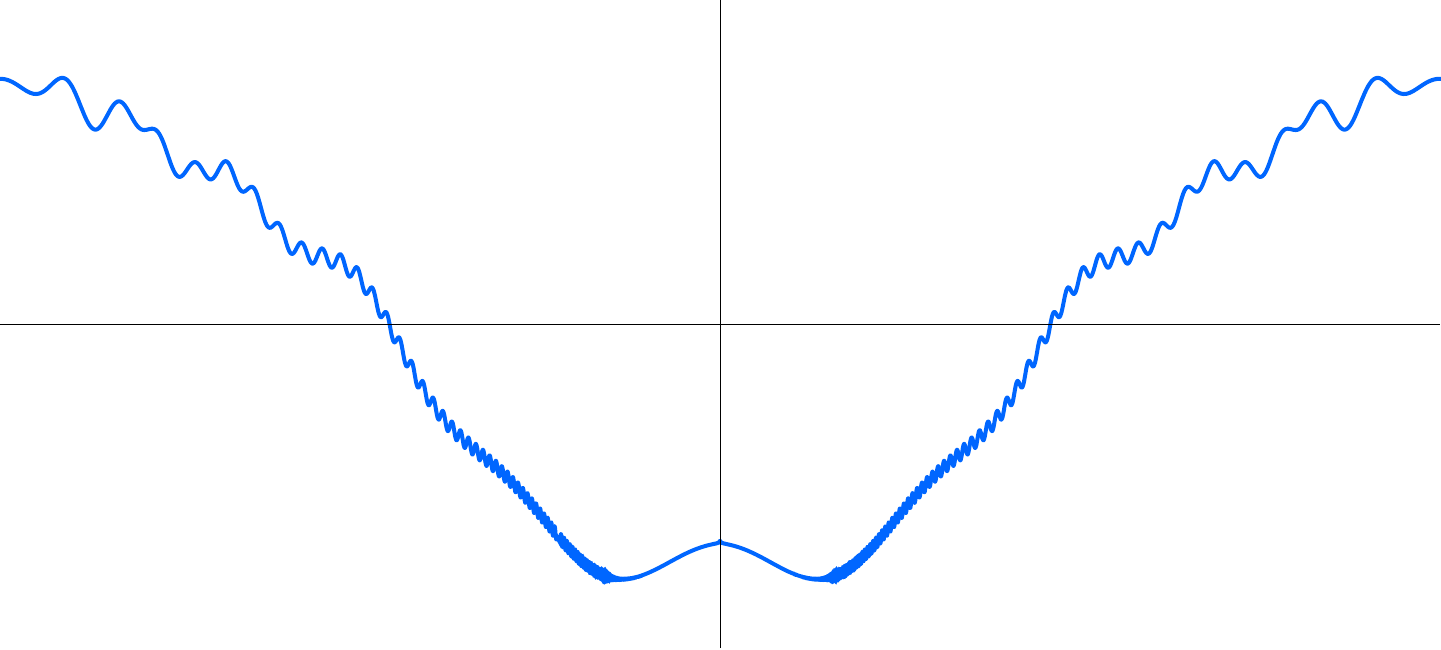}{30.}\vskip-20pt
\Fiiil{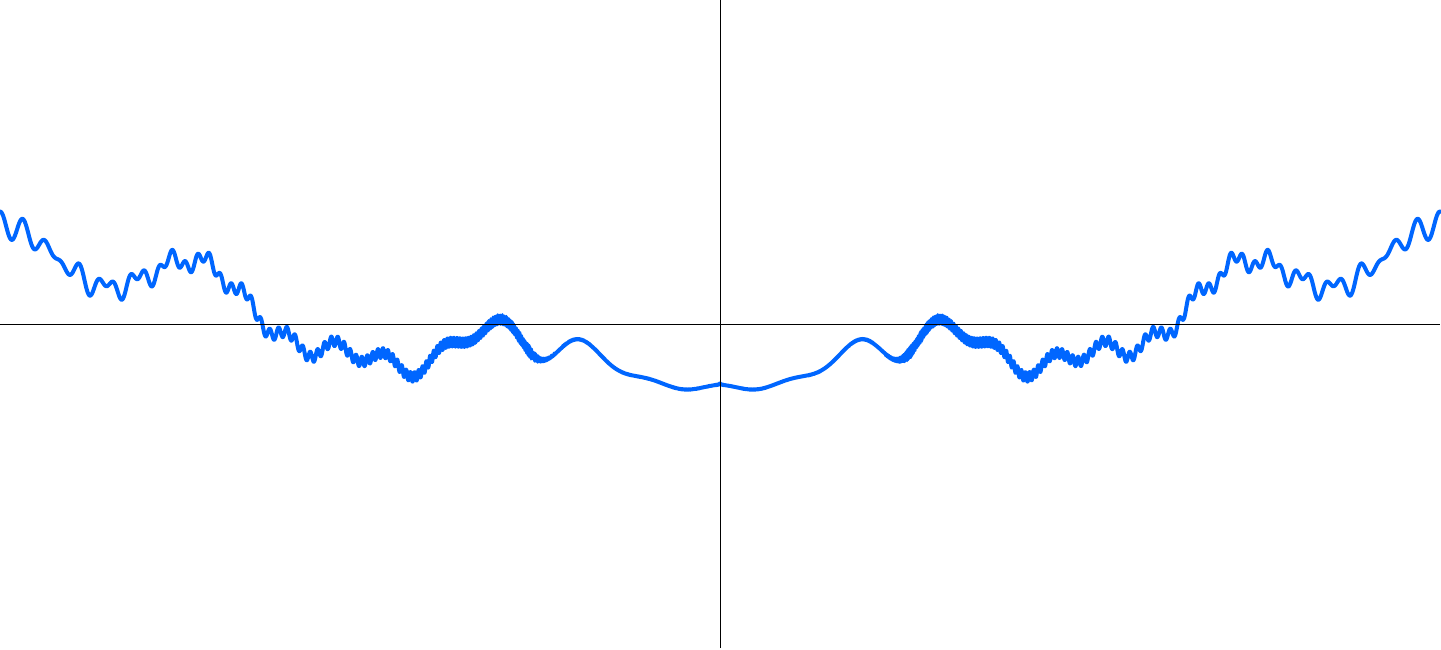}{50.}{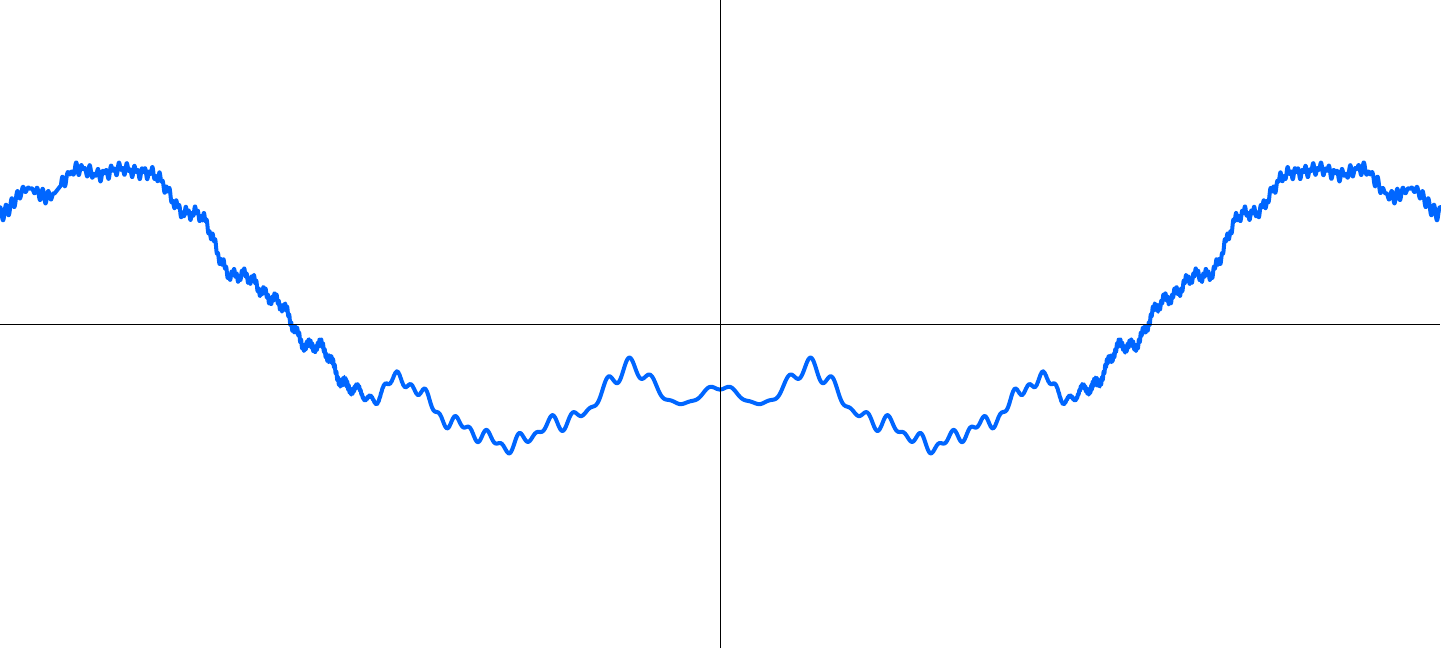}{100.}{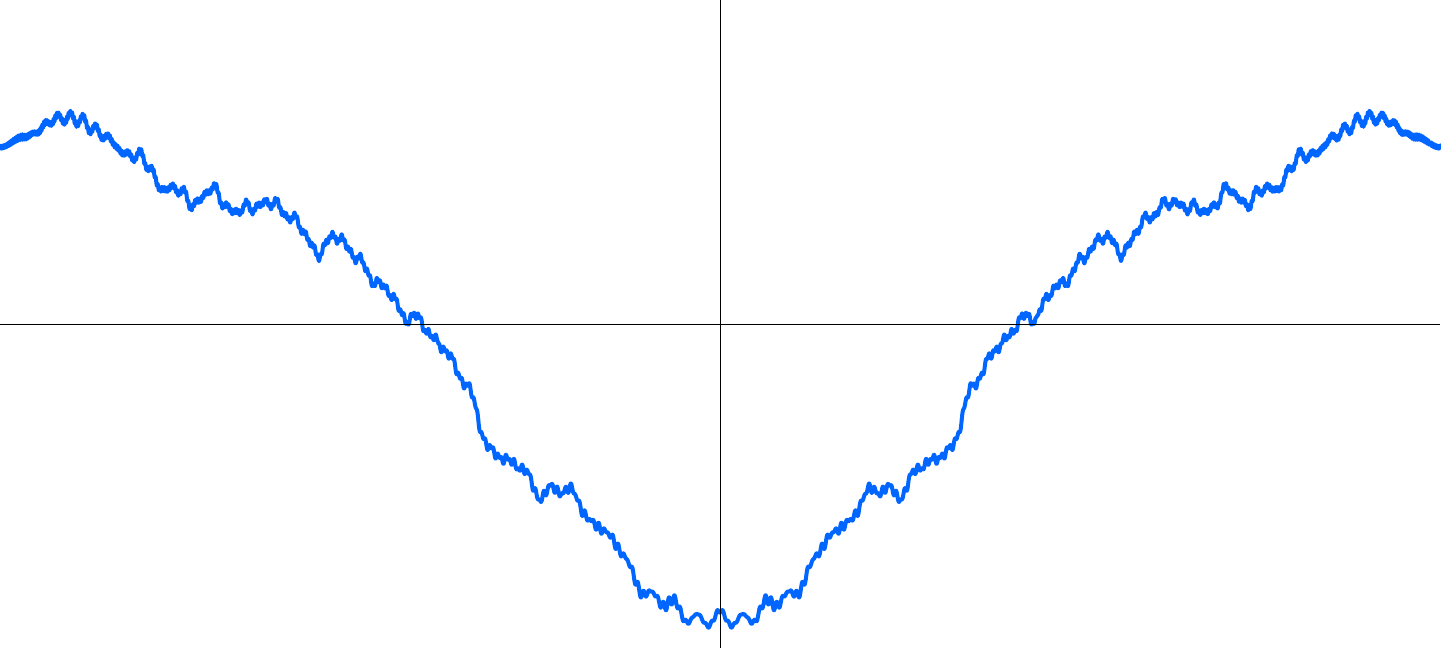}{200.}}
\vskip5pt

As above, the plotted profiles are the order $N=1000$ partial sums of the Fourier series \eqs{uF}{LambdFak}.  Plots of both smaller and larger order partial sums, \eg $N=500$ and $1500$, exhibit no noticeable change in the overall profile at the resolution provided by the figures, and thus strongly suggest convergence of the series solution.  Moreover, it appears that the graphs are fractal, at least on a subinterval centered at the location of the oscillator, but perhaps smooth on the remainder.  The width of the fractal region appears to grow as $t$ increases, so that by $t\geq 50$ or so the entire interval appears to be fractal.

After making our numerical observations, we were fortunate to receive a new paper by Erdo\utxt gan and Shakan, \rf{ErSh}, that contains a proof the following result, which can be used to guarantee convergence of our series solutions and, moreover, confirms the fractal nature by providing an estimate of its (maximal) fractal dimension.

\newpage

\Th{ErSh} Suppose $c_k$ are the complex Fourier coefficients of a function of bounded variation. Let $\omega (k) \simx \abs k^{1/2} $ as $k \to \infty $, then, for any $t \ne 0$, the ``dispersive'' Fourier series
\Eq{fs}
$$v(t,x) \simx \Sumkii c_k \,e^{\i(k \:x - \omega(k) \:t) }$$
converges to a function whose real and imaginary parts have graphs whose maximal fractal dimension $D_t$ satisfies the following estimate\/\ro: $\fr54\leq D_t \leq \fr74$.

In our case, referring back to \eq{LambdFak}, the initial Fourier coefficients  $c_k$ are proportional to $p_k \pm \i q_k$.  
A simple partial fraction decomposition reduces this to a linear combination of Fourier series with coefficients of the form $c_k = 1/\bbk{\omega (k)^2 +  \alpha }$ for real $\alpha $.  In the special case $\omega (k) = \sqrt{\abs k}$, we have
\Eq{ck}
$$\req{c_k = \fra{\abs k +  \alpha} = \fra{\abs k} - \frac \alpha {k^2} + \frac{\alpha ^2}{k^2(\abs k + \alpha) },\\k \ne 0.}$$
The final term is $\O{1/\abs k^{3}}$ and hence represents a Fourier series of a continuously differentiable function. 
The first two terms are the Fourier coefficients of explicitly summable series, namely
\Eq{ckx}
$$\Sumiio k \Pa{\fra{\abs k} - \frac \alpha {k^2}} e^{\i k\:x} = -2 \log\Abs{2 \sin \f2\:x} - \alpha \Pa{\f2\:x^2 - \pii \abs x + \f3\:\pi^2}.$$
However the initial term has a logarithmic singularity at $x=0$, and so is \is{not} of boundfed variation. On the other hand, at $t=0$, the exponentially decaying terms in \eq{LambdFak} cancel out the singularity, since the initial value is identically zero.  Our numerical experiments indicate that there is also no singularity in the solution profile once $t>0$.  Furthermore, the methods of proof based on ``summation by parts'' used in \rf{ErSh} go through even without verifying bounded variation if there is a differentiable function $h(k)$ for $\abs k$ sufficiently large such that $c_k = h(k) = \O{1/k}\sstrut4$ and $d_k = h'(k) = \O{1/k^2}$, and this holds provided $\omega (k) = \O{\abs k^{1/2}}\sstrut4$ and  $\omega '(k)  = \O{1}$. This establishes convergence of the Lamb solution Fourier series and the preceding bound on the maximal fractal dimension of its graph at each $t >0$.  (The authors thank Burak Erdo\utxt gan, \rf{Erdogan}, for essential help with the above arguments.)



Finally, consider the case of the linear regularized Boussinesq equation with Lamb forcing, \eq{rB}, with asymptotically constant dispersion relation \eq{rBd}.  In this case, the solution profiles exhibit much more noticeable high frequency oscillations, whose overall amplitude, as indicated by the ``width'' of the fattened graphs, is time-varying. At earlier times, \eg $t\leq 20$, the graph is noticeably thicker, while later on, \eg at $t=57$ and $t=95$, it has thinned out, with barely any noticeable superimposed oscillations. A while later, at $t=205$, the larger scale oscillations thickening the graph have re-emerged. As yet, we have no explanation for this observed phenomenon.  In this case, the persistence of high frequency oscillations in the graphs of the partial sums remind us of the weak convergence of Fourier series to distributions, the prototypical example being the weak convergence of the Dirichlet kernel to the delta function, \rf P.


\Fig7{The Dispersive Periodic Lamb Oscillator for the Regularized Boussinesq Model}
{\prefigskip=5pt\figlabl=4.75cm\figwidth=4.5cm
\Fiiil{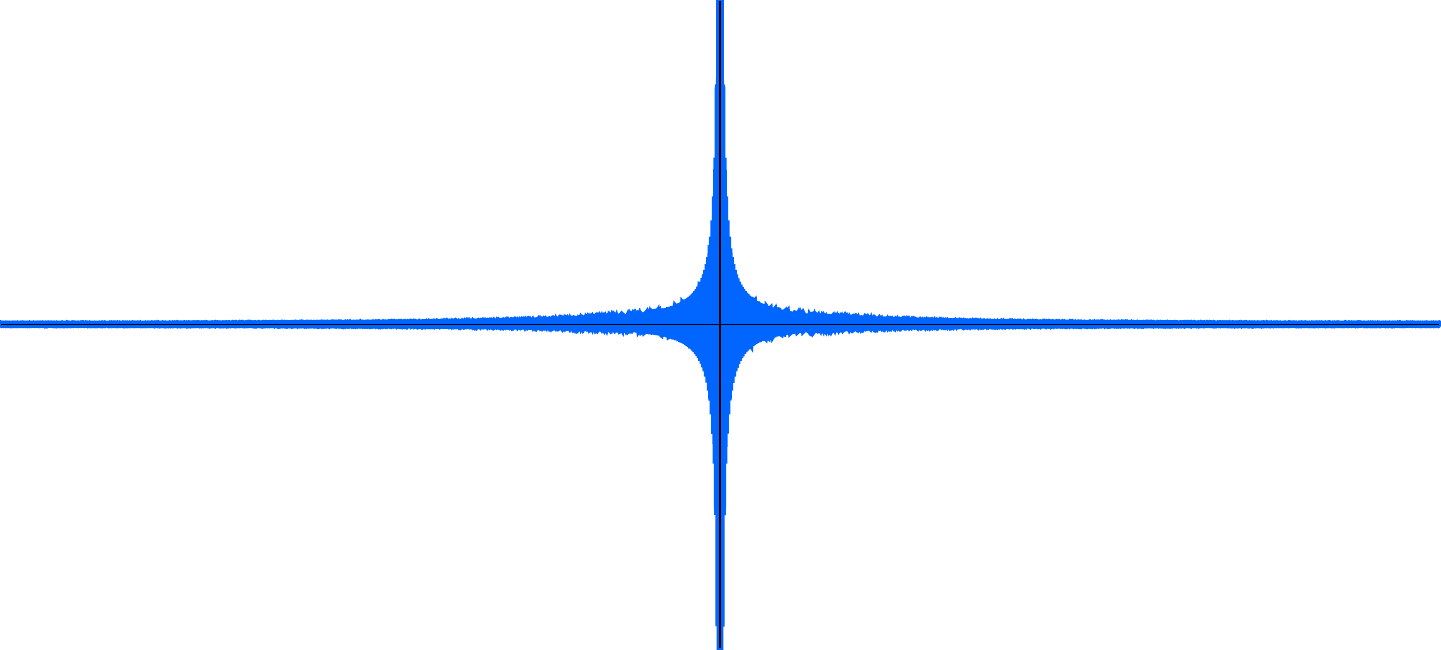}{1.}{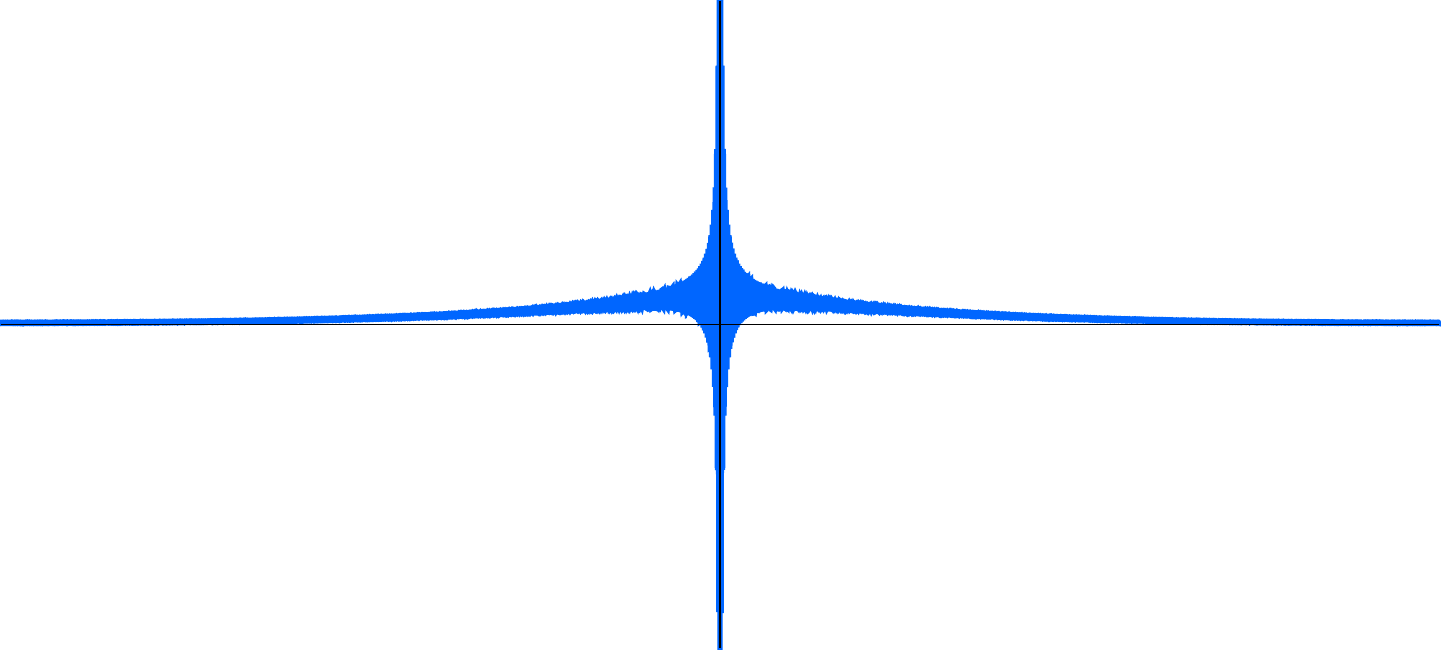}{2.5}{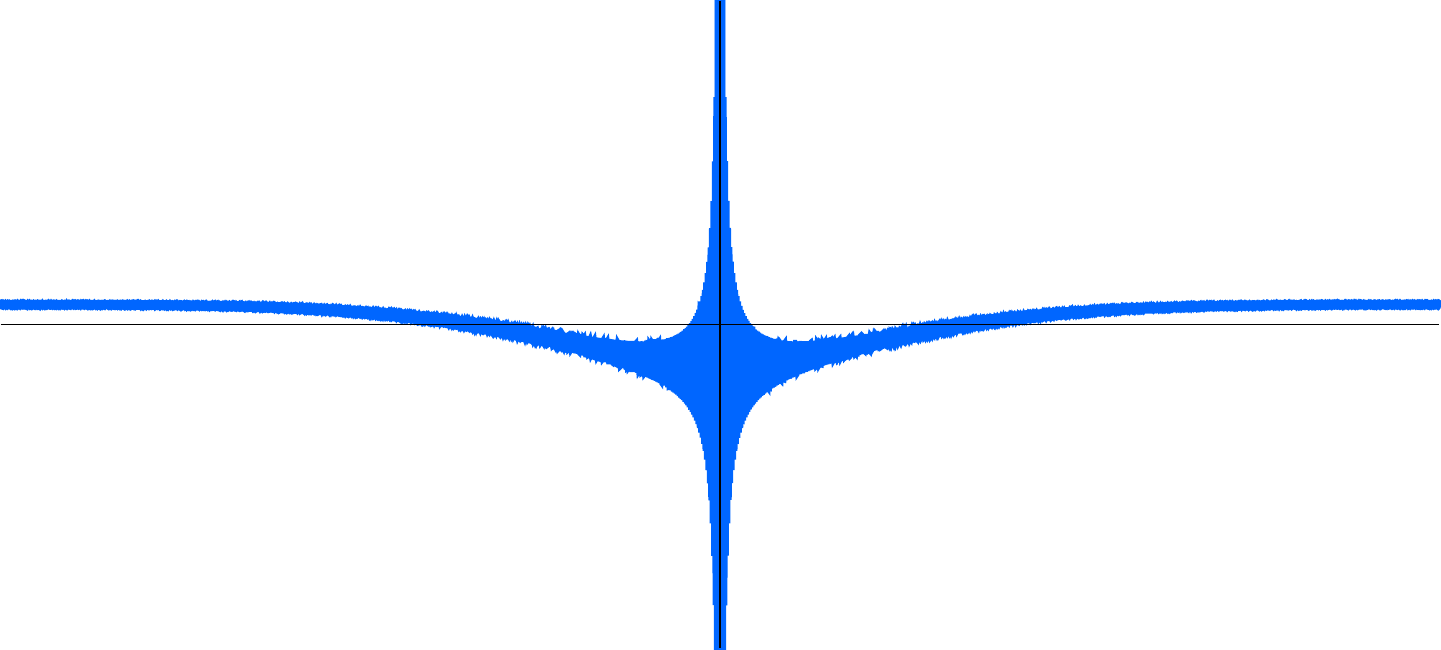}{5.}\vskip-20pt
\Fiiil{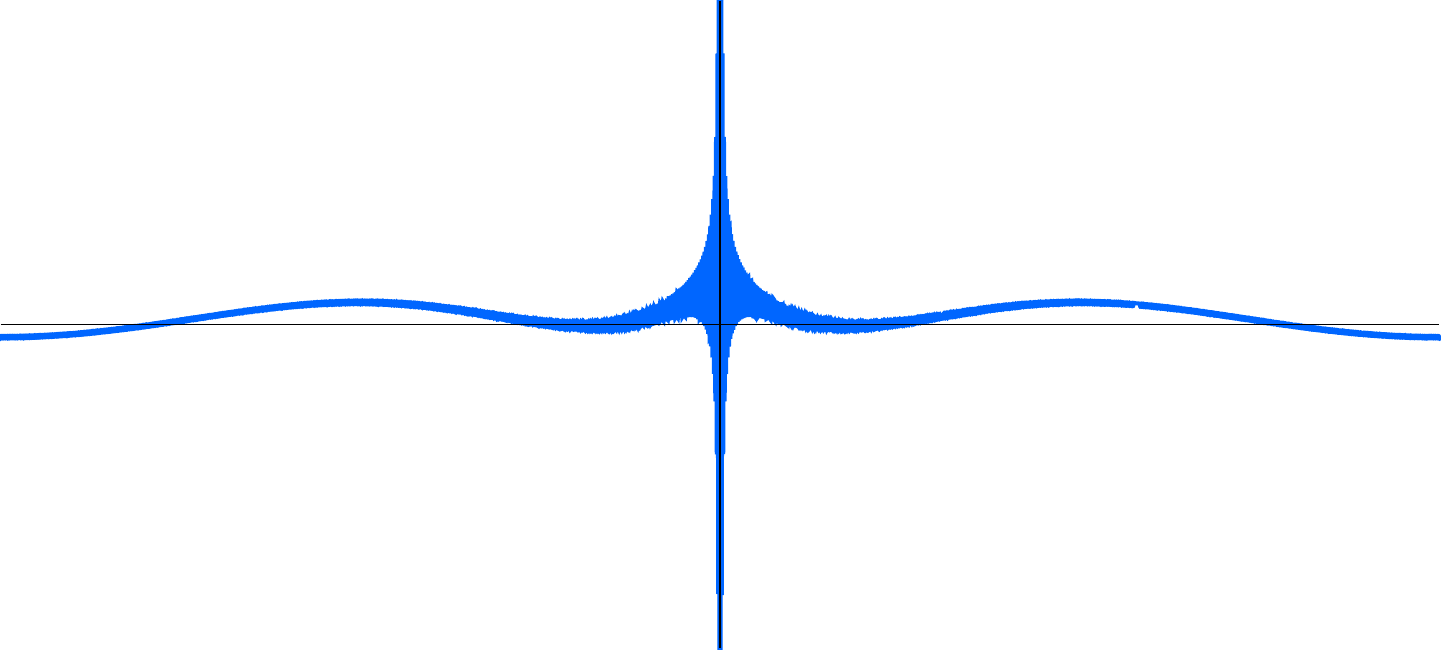}{10.}{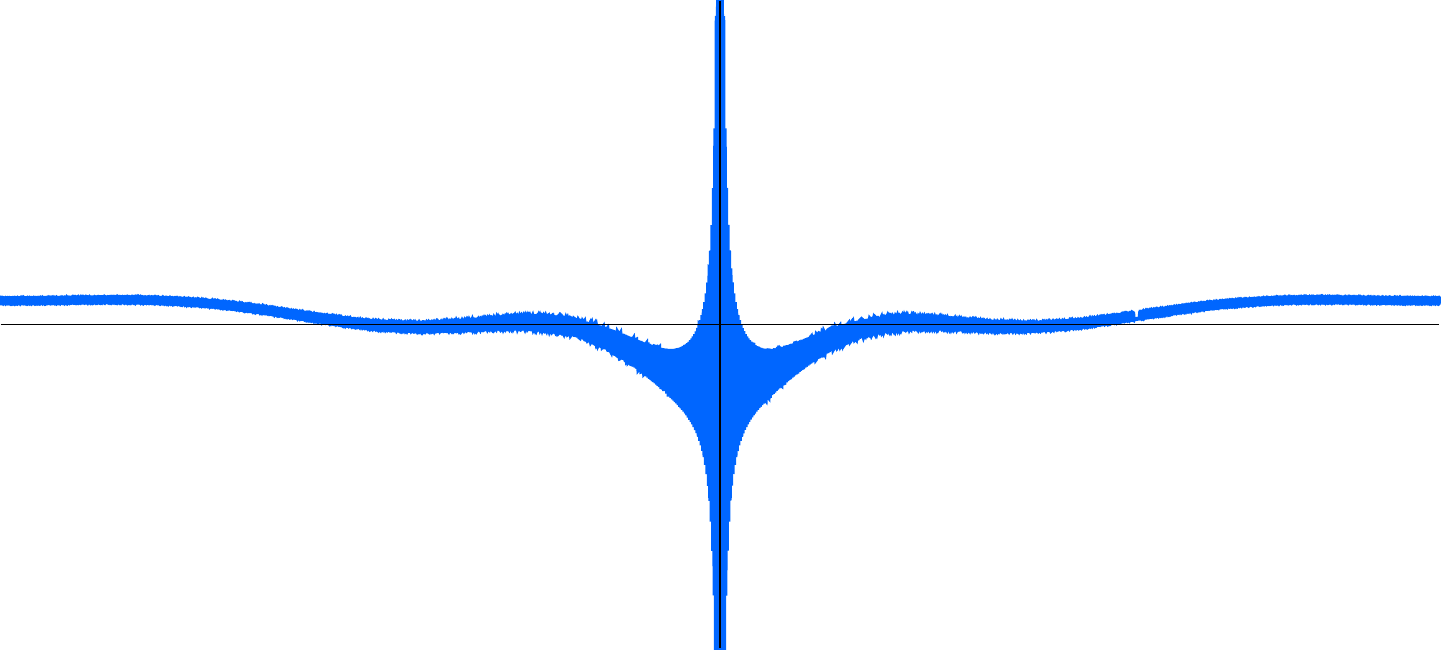}{20.}{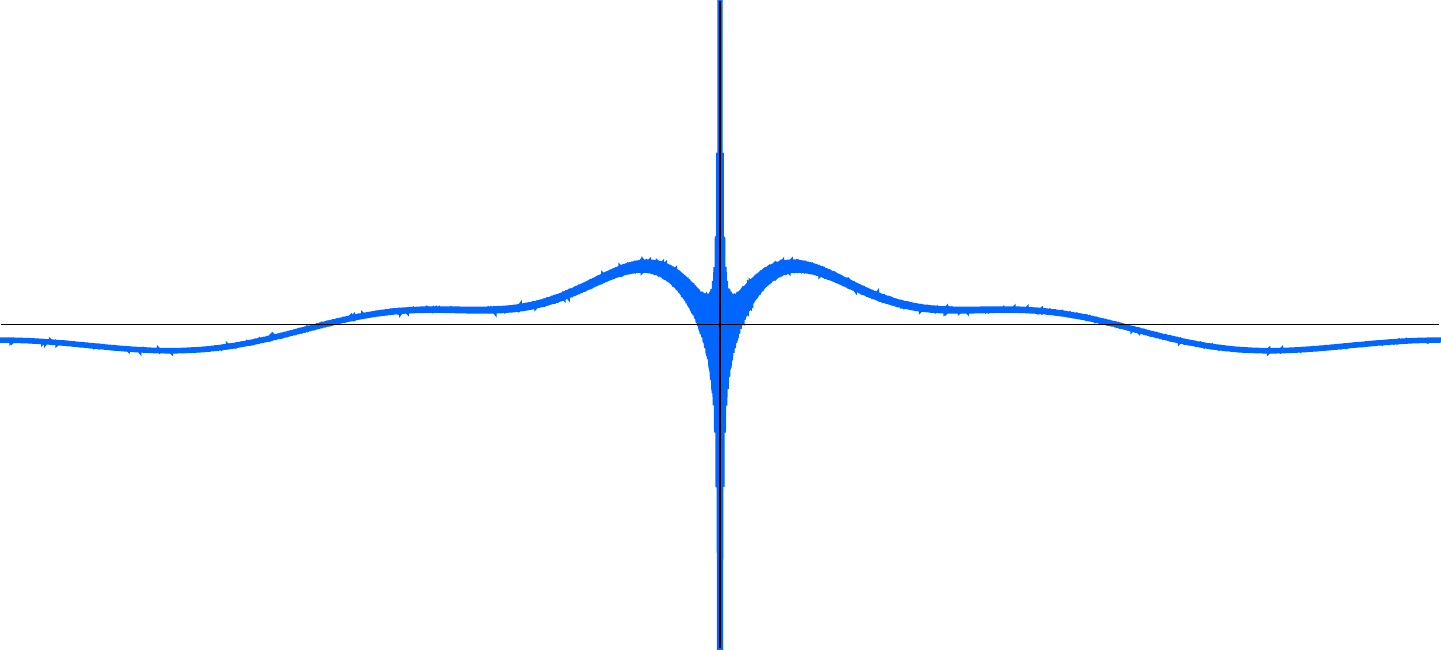}{30.}\vskip-20pt
\Fiiil{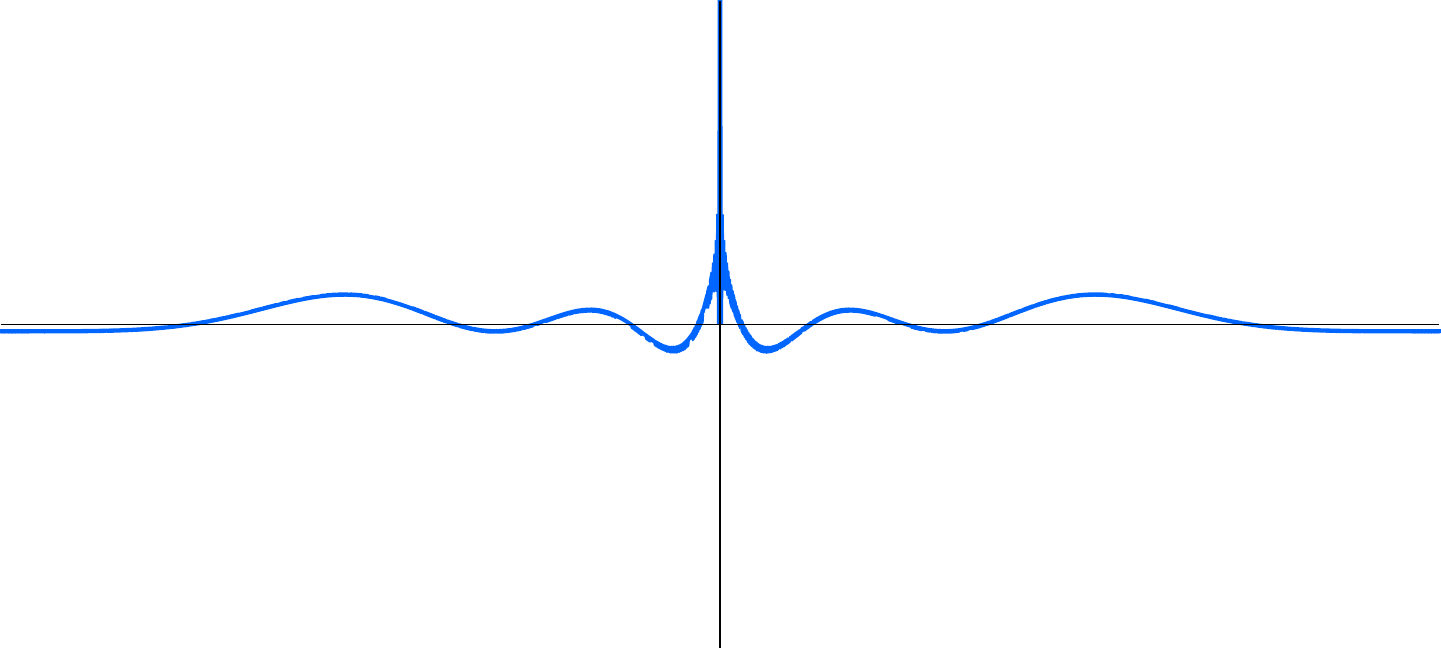}{57.}{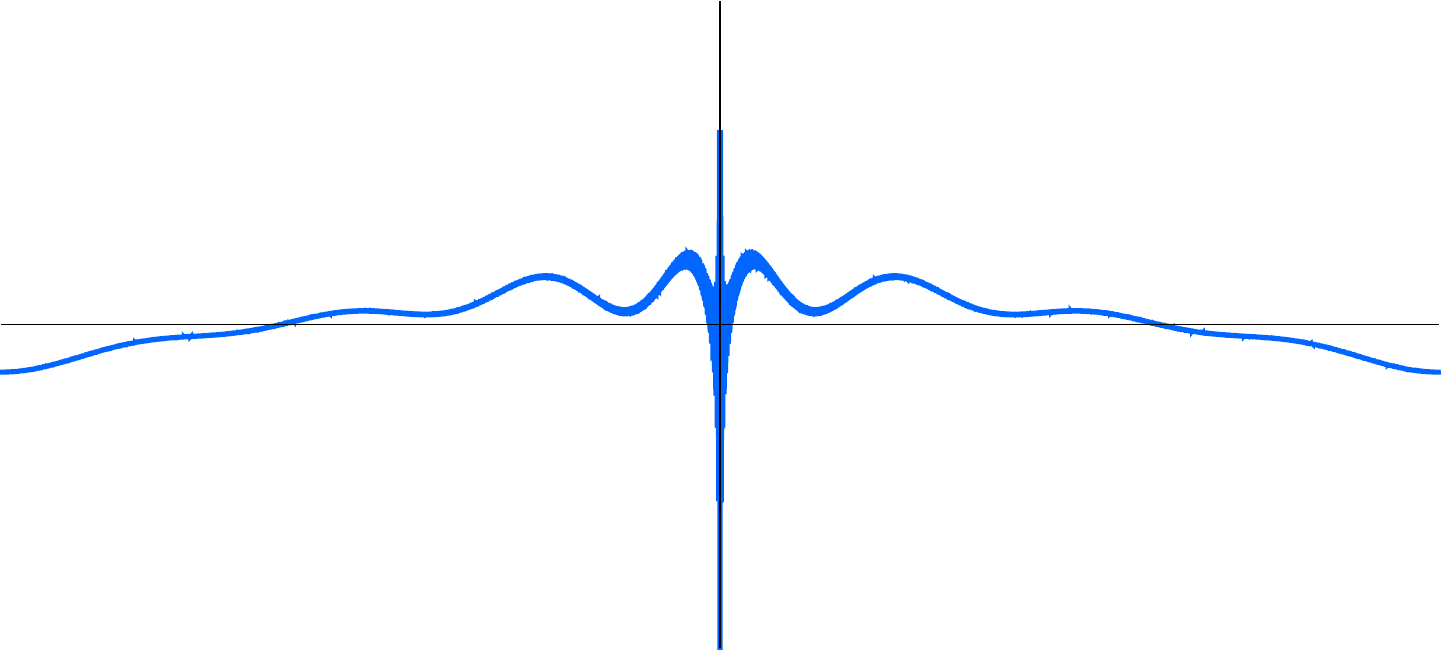}{95.}{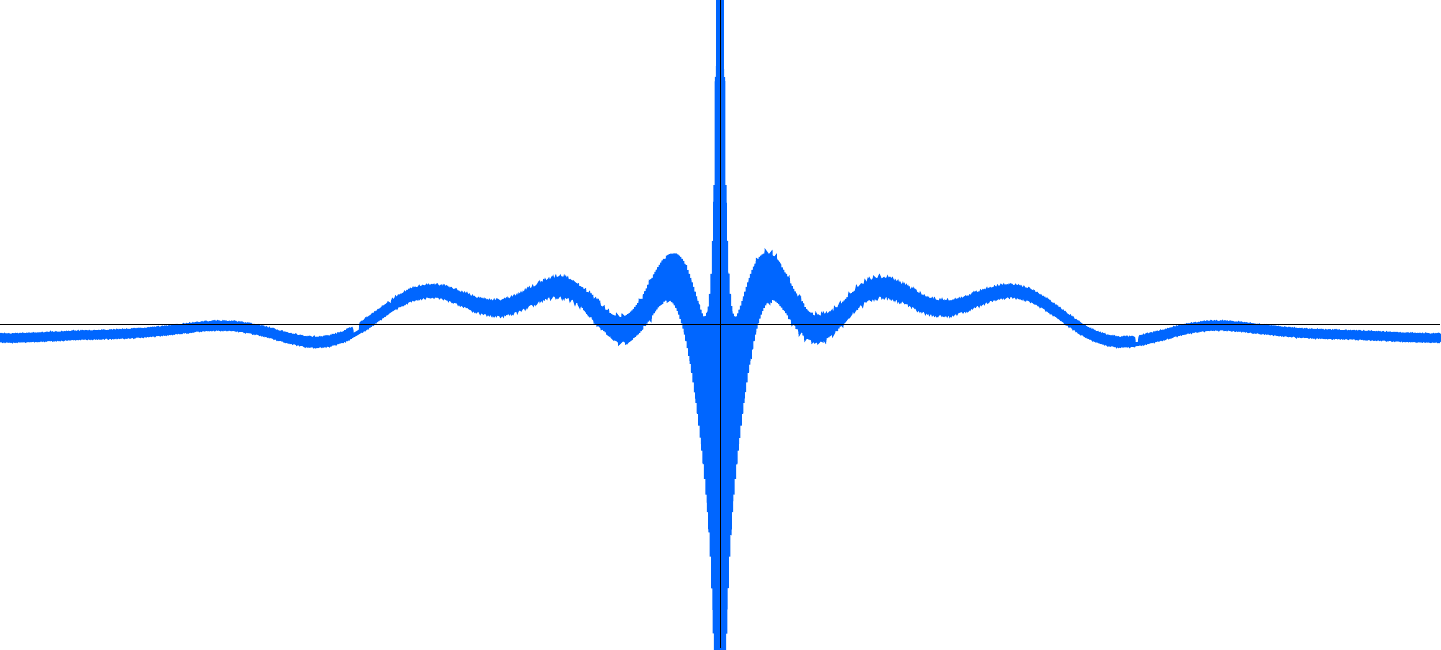}{205.}}

\Section 1 Unidirectional Models.

Now let us turn our attention to unidirectional dispersive wave models forced by a Lamb oscillator. We begin with the linear wave equation and factor the differential operator as usual, \rf P:
$$\partial _t^2 - c^2\partial _x^2 = (\partial _t - c\:\partial _x)(\partial _t + c\:\partial _x),$$
where the factors govern the unidirectional waves $u(t,x) = f(x \pm c\:t)$, that move off in opposite directions, each of which is constant on the characteristic lines associated with its annihilating differential operator.

The bidirectional Lamb solution \eq{Lambsol} is a superposition of left and right moving waves, which do not overlap since the forcing only occurs at the origin:
\Eq{ubi}
$$u(t,x) =\mcases{ v(t,x), & x > 0,\\ v(t,-\:x),&x < 0,}\where
v(t,x) = \mcases{f(c\:t - x), & x > 0,\\0,&x < 0,}$$
and $f(t)$ is as above, \eq{Lambf}.  Thus, $v(t,x)$ solves the quarter plane initial-boundary value problem
\Eq{p4}
$$\req{v_t + c\: v_x = 0,& v(0,x) = 0,\\v(t,0) = h(t) = f(c\:t),& x,t > 0,}$$
or, equivalently, the forced unidirectional wave equation
\Eq{f1}
$$\req{v_t + c\: v_x = h(t) \, \delta (x),&v(0,x) = 0,\\ t > 0,} $$
since we assume that $v(t,x) = 0 $ for all $x < 0$.
As in \sect d, we have absorbed the factor of $-2\:c$ into the forcing function $h(t)$.
More generally, we can replace the spatial transport term to construct a general unidirectional linear dispersive Lamb model
\Eq{ul}
$$\req{v_t + L\br v = h(t) \, \delta (x),&v(0,x) = 0,\\ t > 0,} $$
in which the linear integro-differential operator $L\br v$ has real dispersion relation $\omega (k)$.

On the full line $x \in \R$, we can solve the forced initial-boundary value problem \eq{p4} using the Fourier transform:
\Eq{uniD}
$$ v(t,x)= \frac{C}{4\: c \pii}\int_{-\infty}^\infty \frac{e^{-\i\omega(k) t} \varsigma +e^{-\beta t}\bbk{(\i\omega(k)-\beta)\sin \varsigma t-\varsigma\cos \varsigma t }}{\sigma^2-2\i\beta\,\omega(k)-\omega(k)^2}\;e^{\i k\:x}~dk.$$
The convergence of the integral depends upon the form of the dispersion relation $\omega (k)$.

Here we focus our attention on the periodic initial-boundary value problem corresponding to \eq{ul} on the interval $-\pii < x < \pi$. We expand the solution in a Fourier series:
\Eq{vF}
$$v(t,x) = \f2\:a_0(t) +\ \Sumi k \bbk{a_k(t) \cos k\:x + b_k(t) \sin k\:x}.$$
Substituting \eq{vF} into the initial value problem \eq{ul} produces the Fourier coefficients 
\Eq{vFcf}
$$\eeq{
a_0(t)=\frac{C}{2\:c\:\pi\: \sigma^2}\bbk{\varsigma-e^{-\:\beta \: t}(\varsigma \cos \varsigma \:t+\beta \sin \varsigma \:t) },\\
a_k(t) = \frac{C}{2\:c\:\pi}\bbk{p_k \cos \omega \:t + q_k \sin \omega \:t - e^{-\:\beta \: t} (p_k \cos \varsigma \:t + r_k \sin \varsigma \:t )},\\
b_k(t) = \frac{C}{2\:c\:\pi}\bbk{-\:q_k \cos \omega \:t + p_k \sin \omega \:t + e^{-\:\beta \: t} (q_k \cos \varsigma \:t + s_k \sin \varsigma \:t )},}$$
for all $k \geq 1$, where $p_k,q_k,r_k$ are given in \eq{pqrw} and 
\Eq{biDs}
$$s_k=\frac{\omega\:(2\beta^2-\sigma^2+\omega^2)}{(\sigma^2 - \omega ^2)^2 + 4\:\omega^2\beta ^2}.$$ 
In contrast to the bidirectional equations, the Fourier coefficients in~\eq{vF} decay like $1/{\omega(k)}$ as $k\to\infty$. As above, we deduce differentiability of the solution profiles provided the dispersion relation grows sufficiently rapidly at large wave number; specifically, \eq{drm} with $m \geq 2$ implies $u(t,\cdot ) \in \Cn$ provided $m-1 > n \in \N$. 


We now examine some representative solutions.  Again, we plot the order $N = 1000$ partial sum of the Fourier series~\eq{vF}. 
First, for a linear dispersion relation $\omega(k)=c\:k$, corresponding to the unidirectional transport equation with a Lamb oscillator \eq{f1}, the solution remains piecewise smooth.

\Fig9{The Unidirectional Periodic Lamb Oscillator for the Transport Model}
{\prefigskip=5pt\figlabl=4.75cm\figwidth=4.5cm
\Fiiil{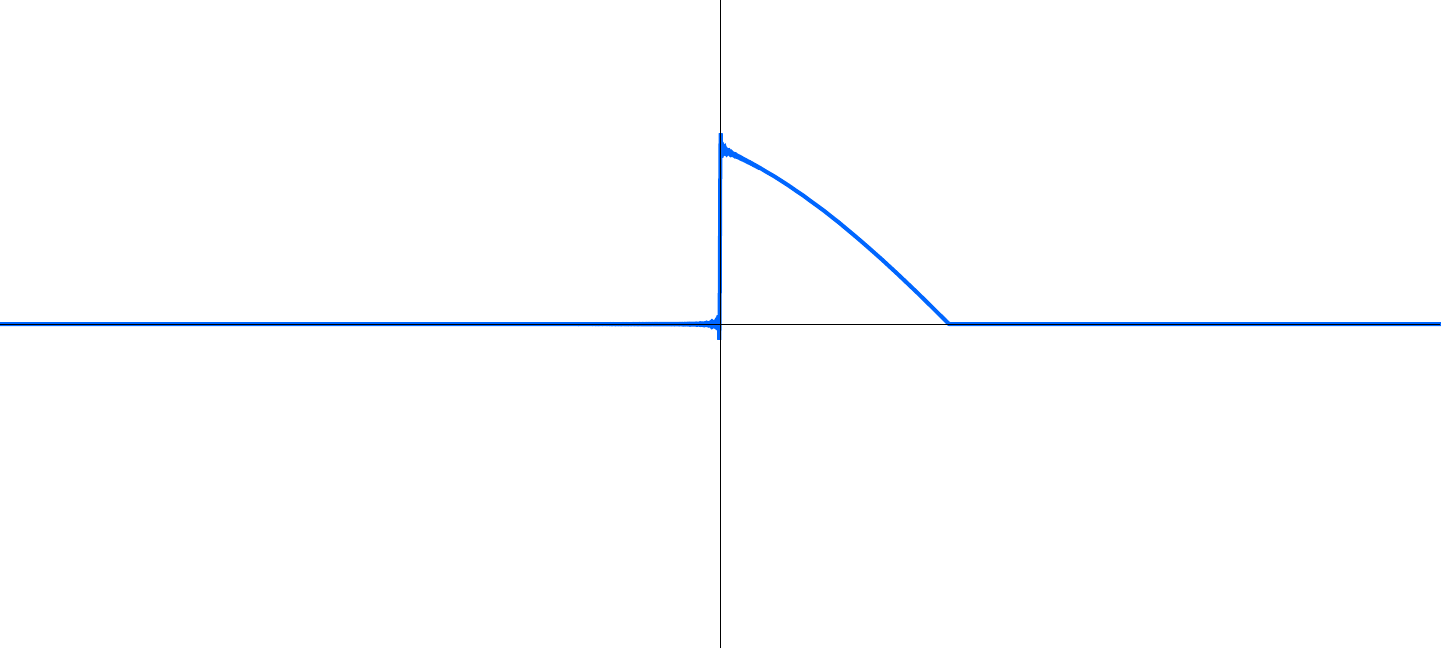}{1.}{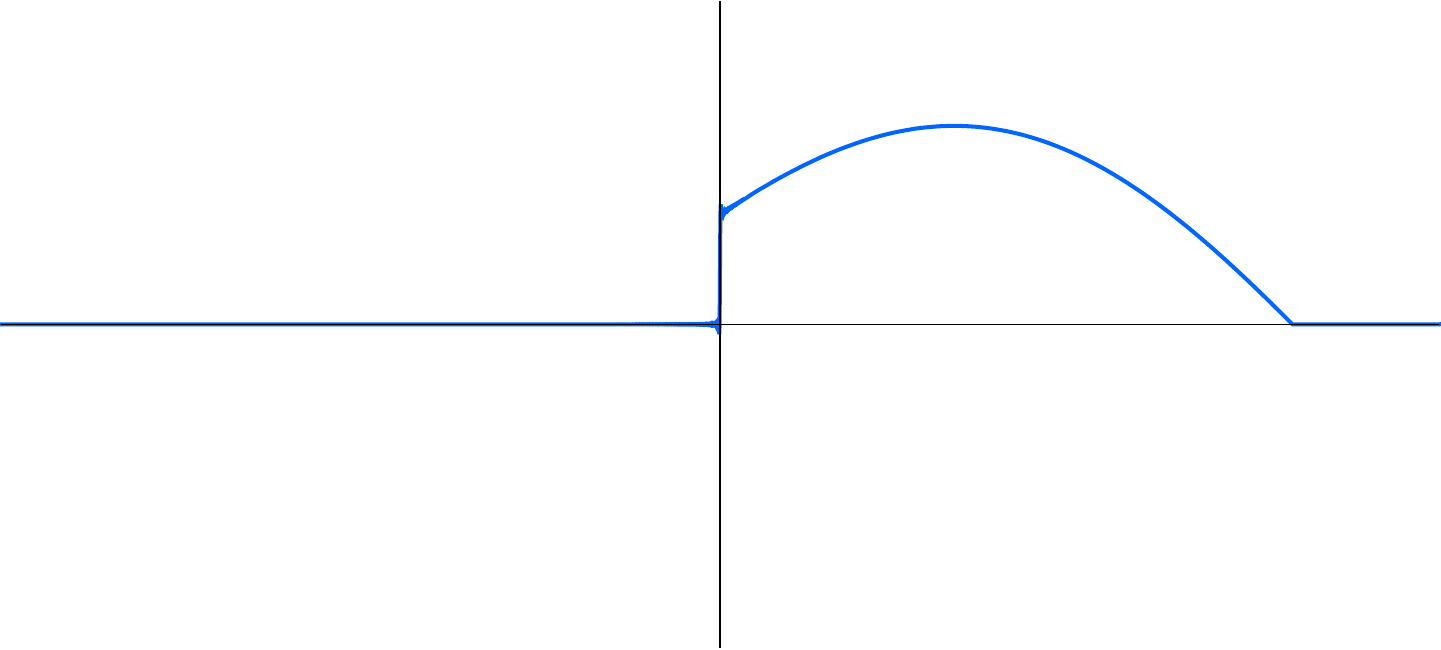}{2.5}{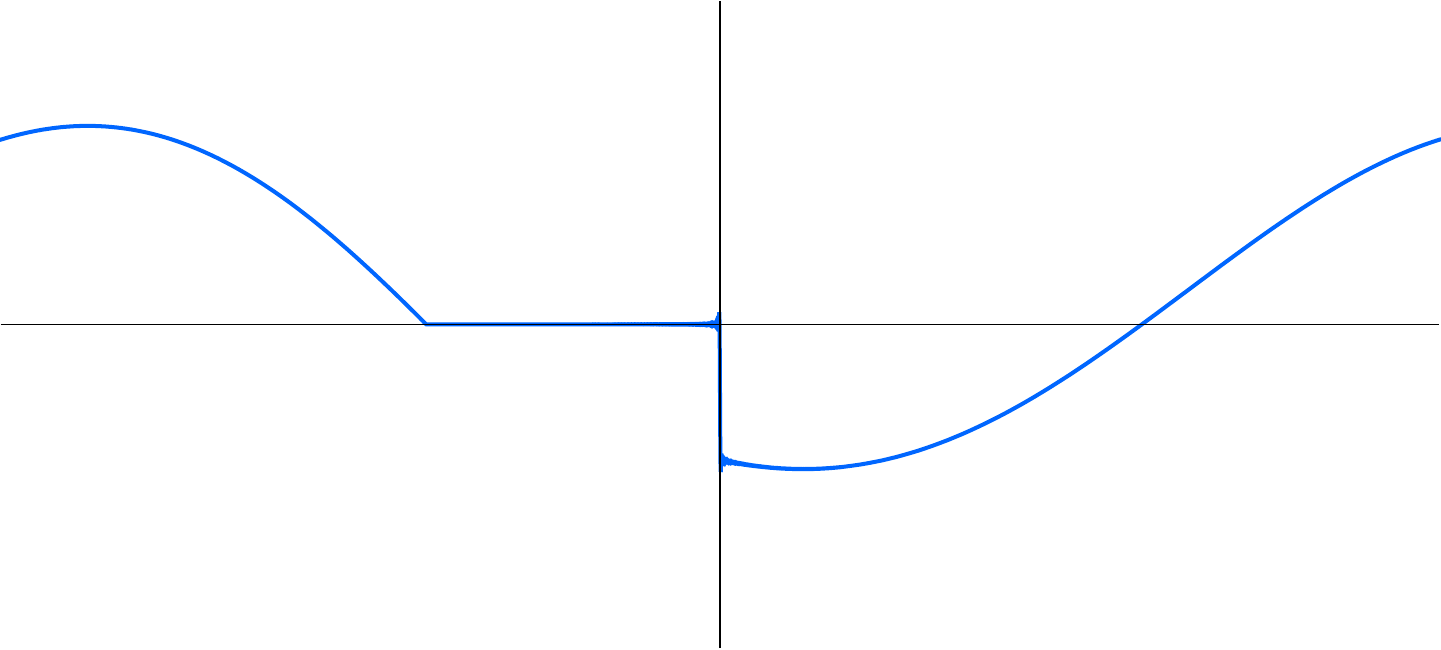}{5.}\vskip-15pt
\Fiiil{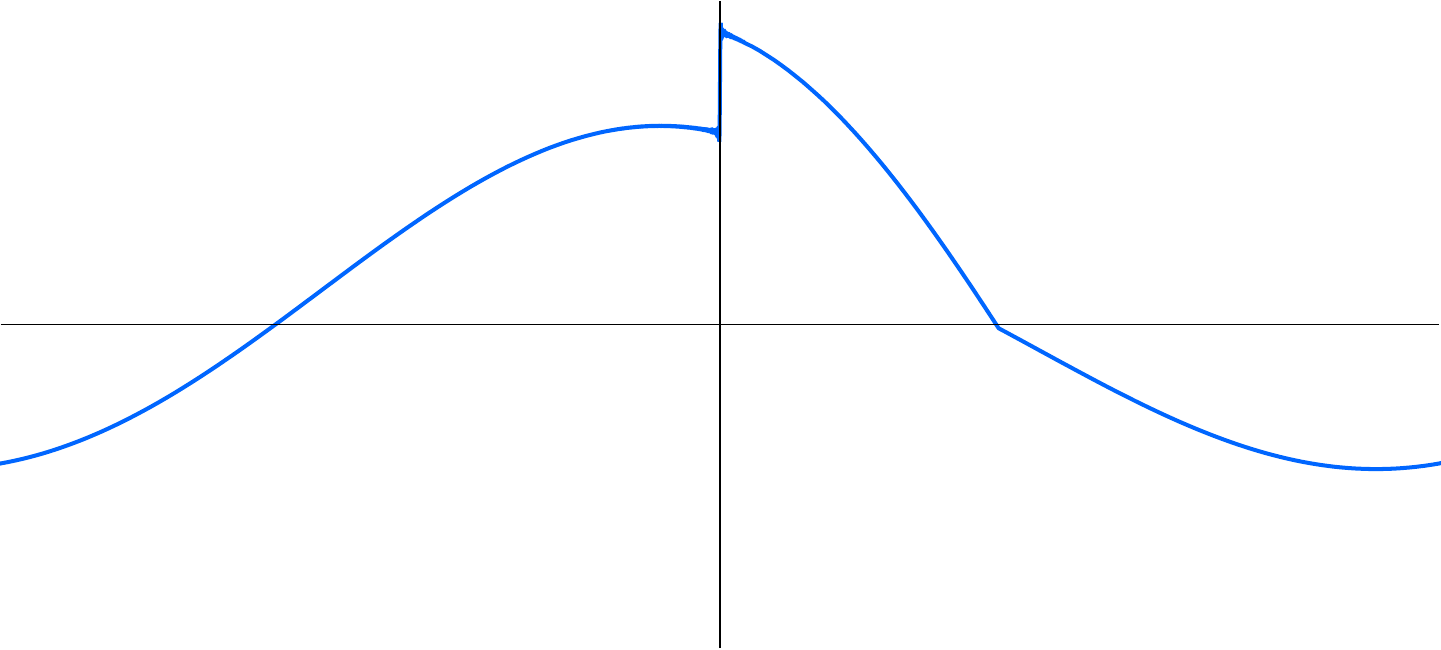}{7.5}{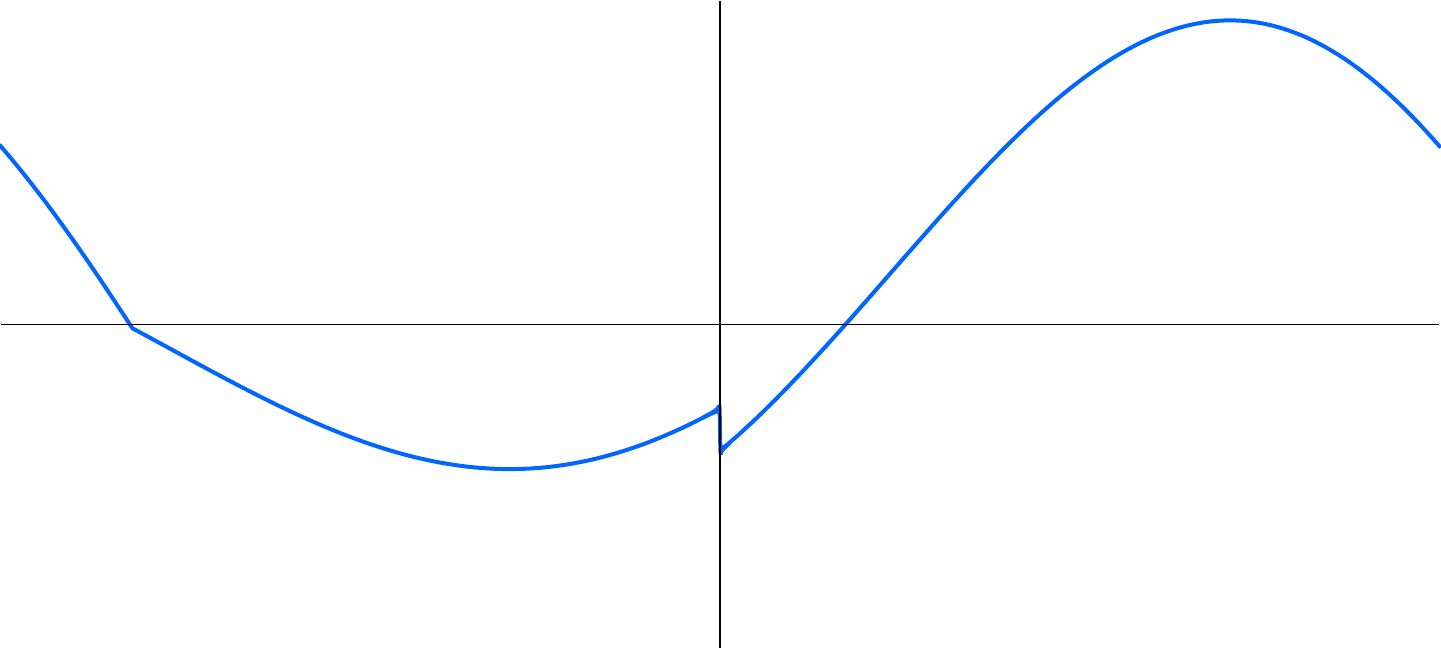}{10.}{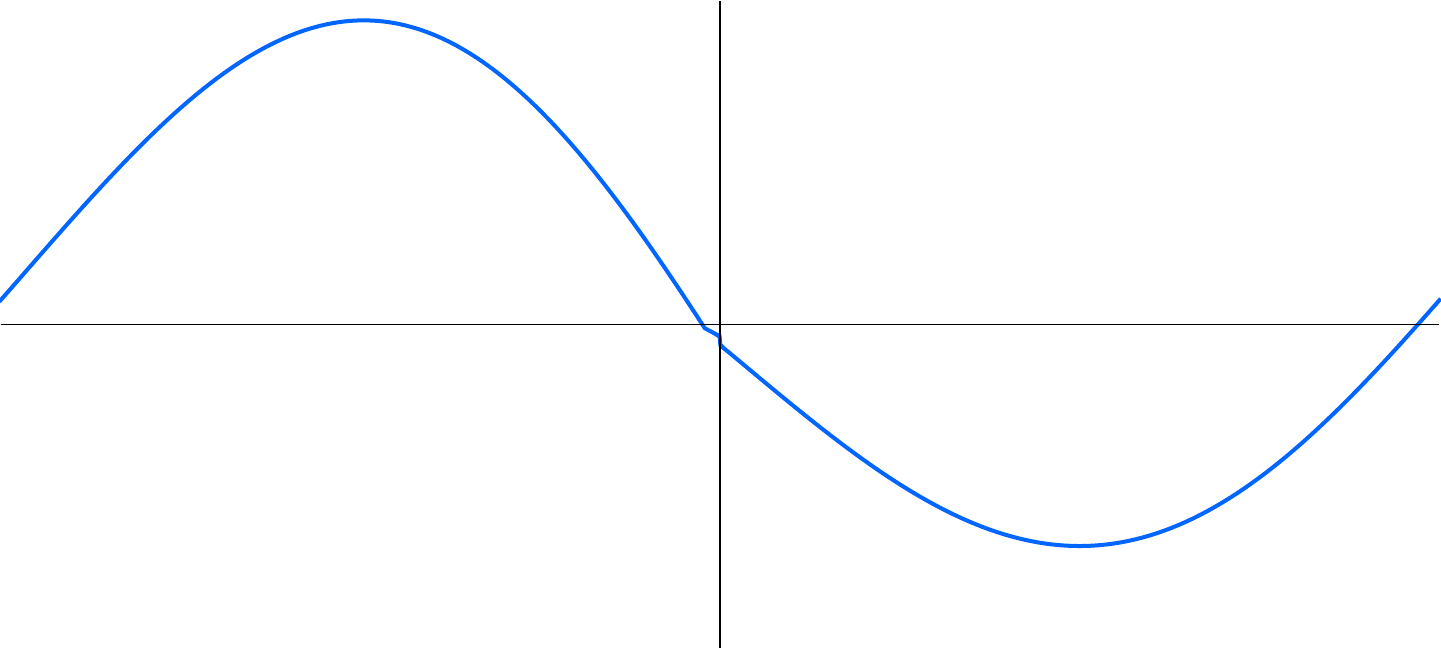}{12.5}}

\Fig{10}{The Unidirectional Dispersive Periodic Lamb Oscillator for $\omega(k)={k^2}$}
{\prefigskip=5pt\figlabl=4.75cm\figwidth=4.5cm
\Fiiil{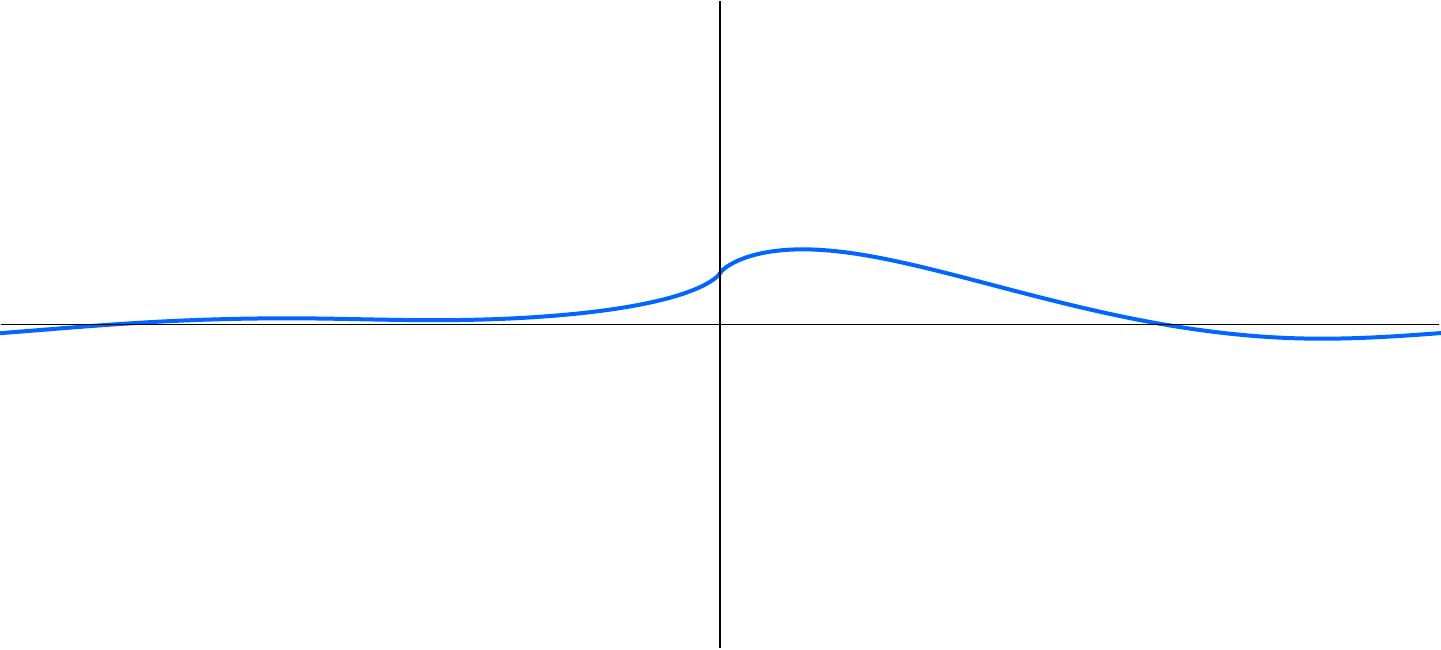}{1.}{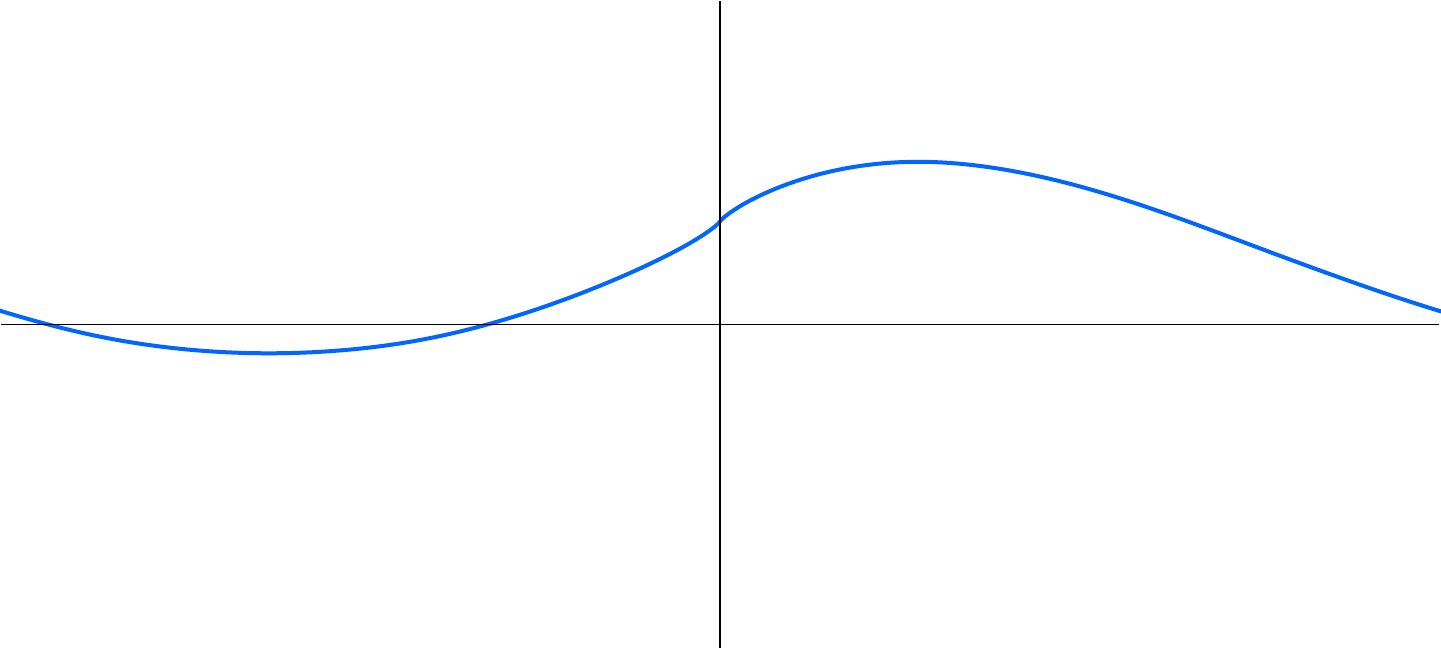}{2.5}{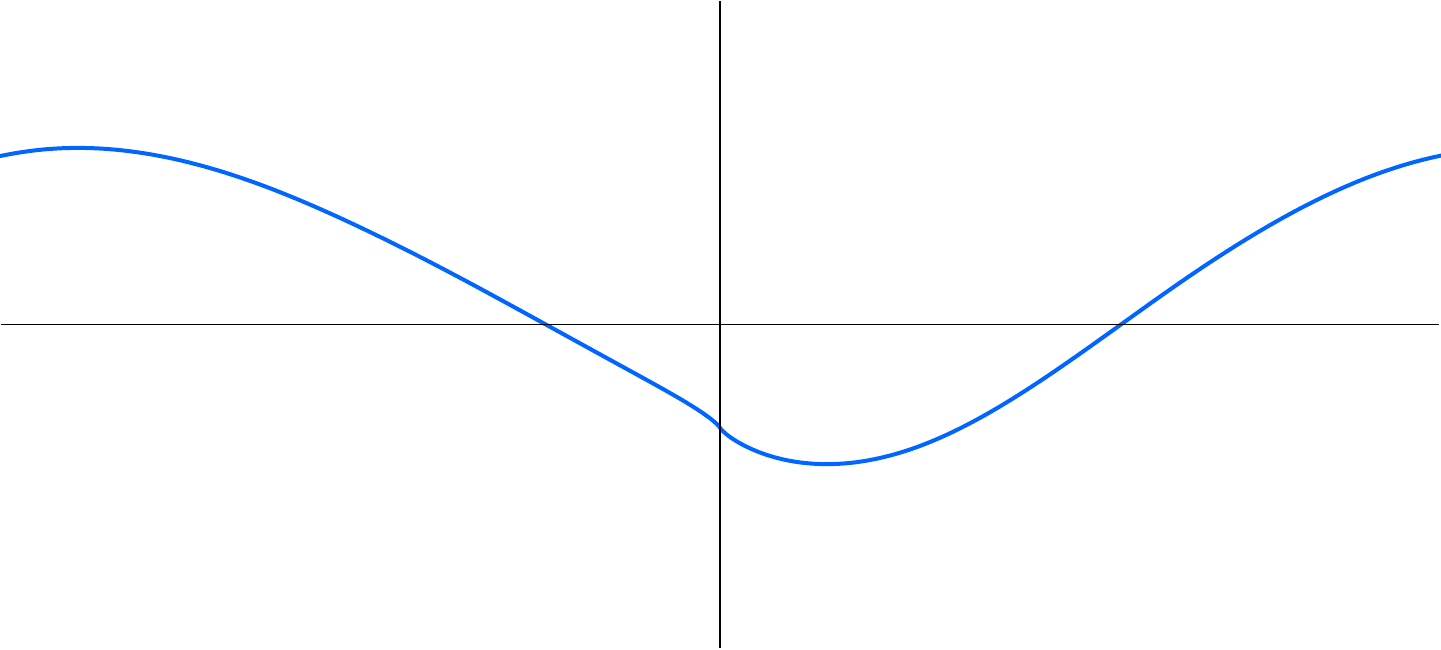}{5.}\vskip-15pt
\Fiiil{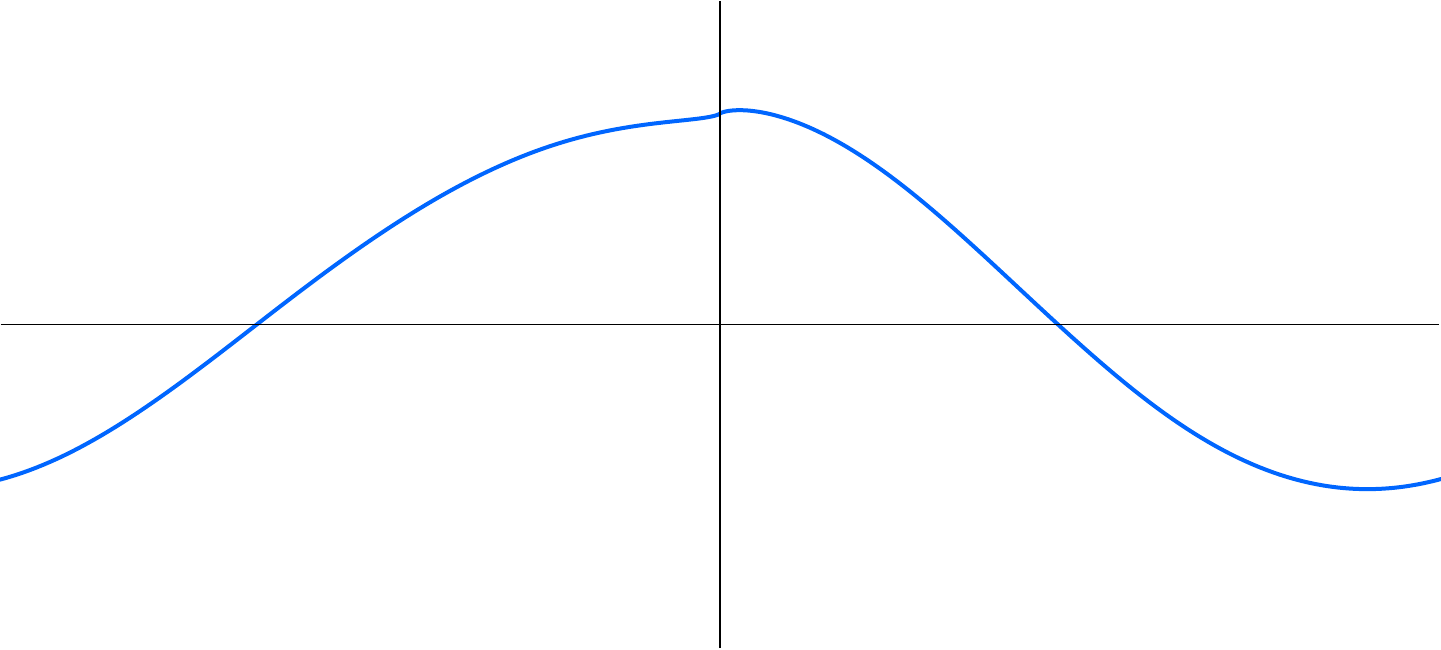}{7.5}{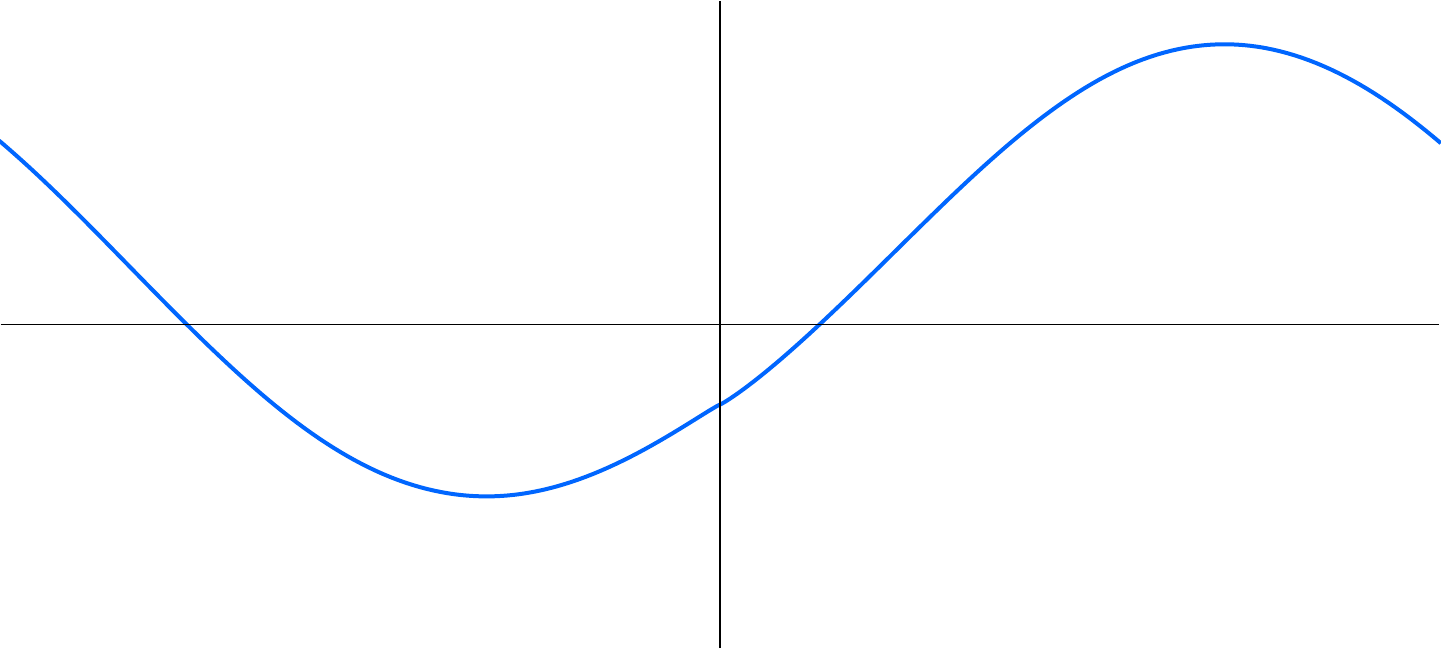}{10.}{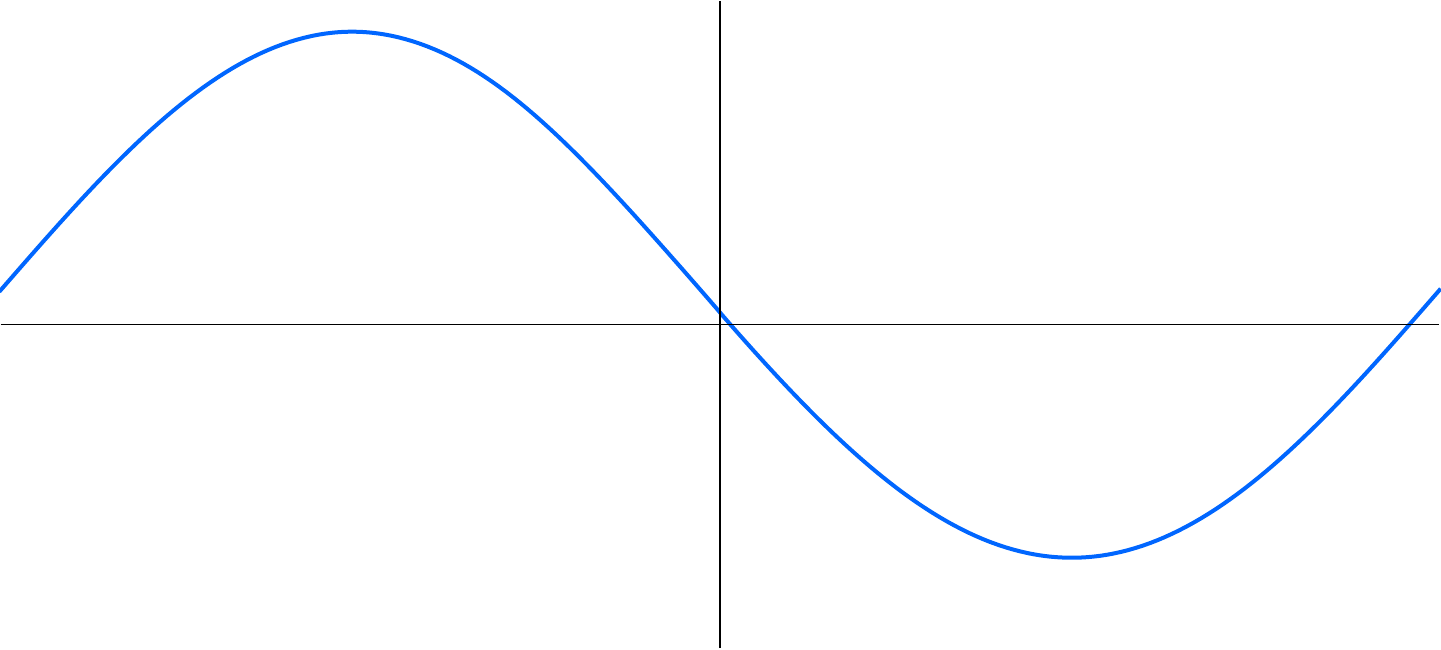}{12.5}\vskip10pt}




As in the bidirectional models, for the higher order dispersion relations, the solution smooths very quickly, as in the following plots.
The amplitude of the oscillations increases with a decaying dispersion relation --- for example, the following figure plots the case $\omega(k)=\sqrt{k}$. Observe the sharp transition region and possible discontinuity for smaller values of $t$, which appear to gradually completely disappear.  One question is whether the solution profiles are fractal anywhere; so far, the plots are inconclusive.  


In the case of asymptotically constant dispersion, \eg $\omega(k)={k^2}/\paz{1+\f3\:k^2}$, we still see significant high frequency oscillations. As in the bidirectional case, the graph thickens and thins at various times, possibly indicating weak convergence to a distributional solution. 

\Fig{12}{The Unidirectional Dispersive Periodic Lamb Oscillator for $\omega(k)=\sqrt{k}$}
{\prefigskip=5pt\figlabl=4.75cm\figwidth=4.5cm
\Fiiil{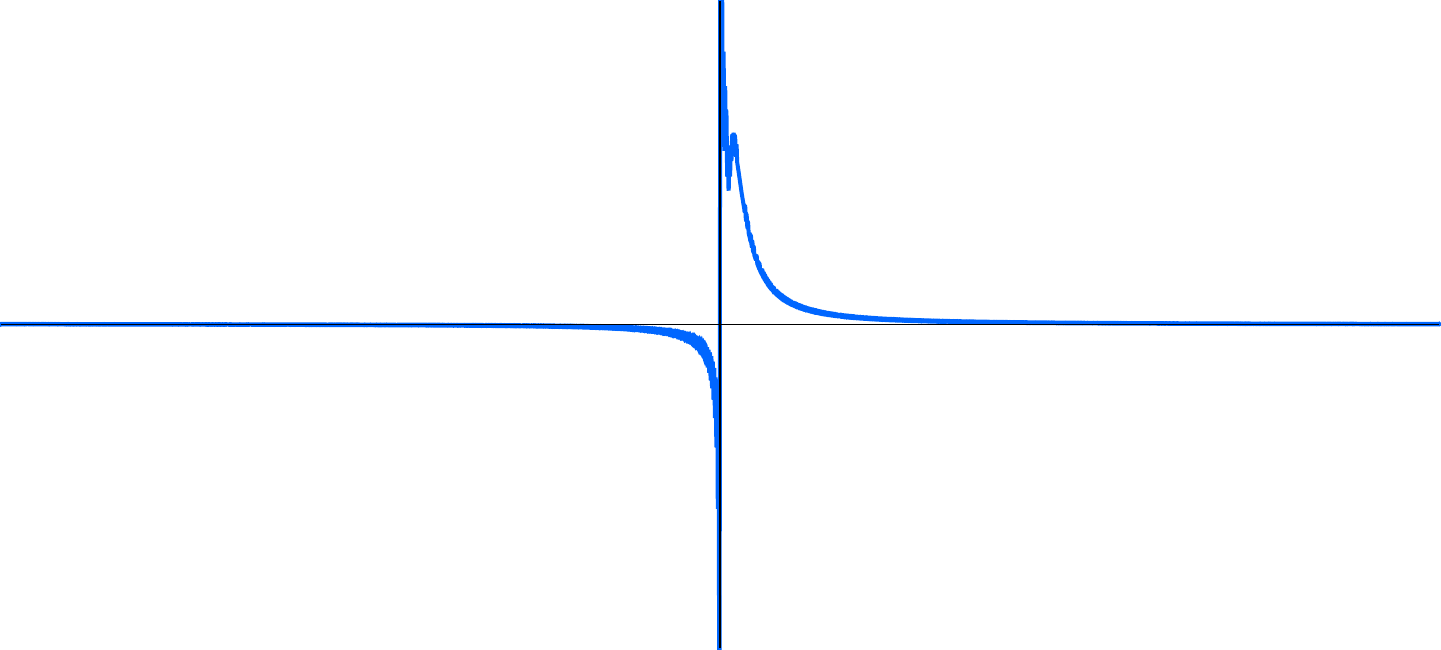}{1}{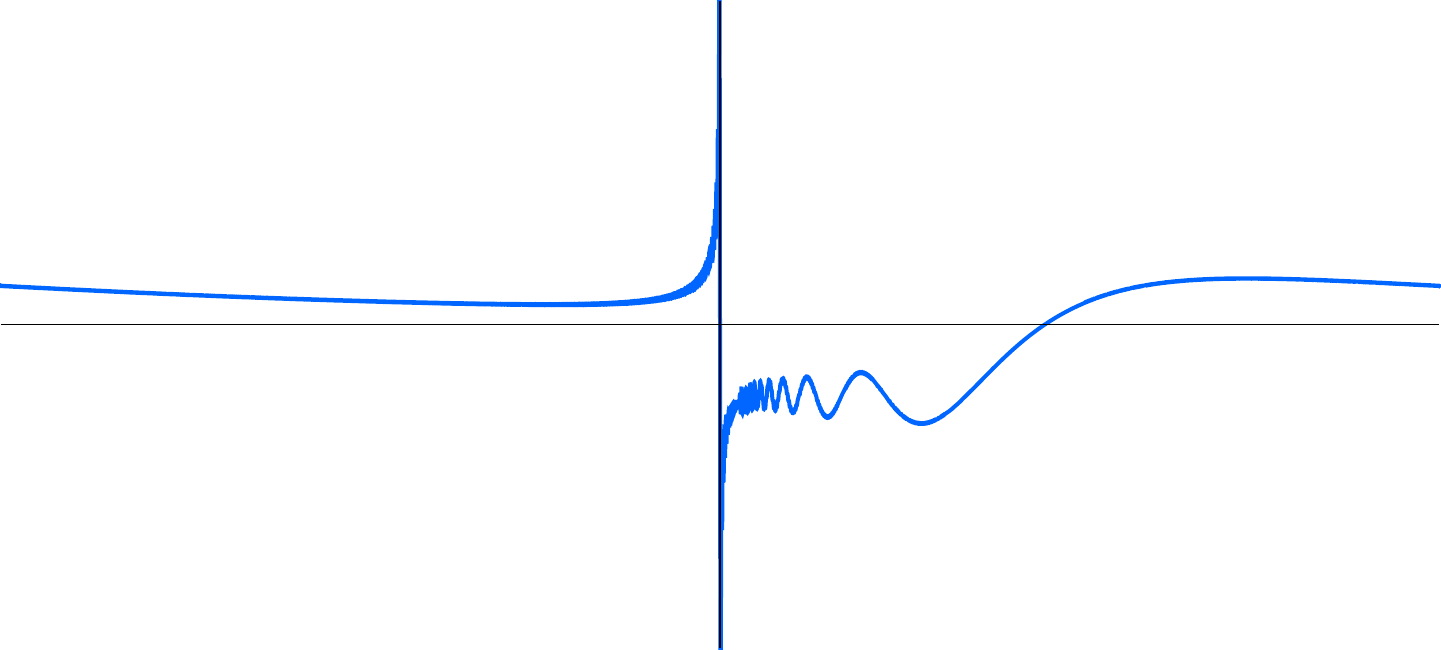}{5.}{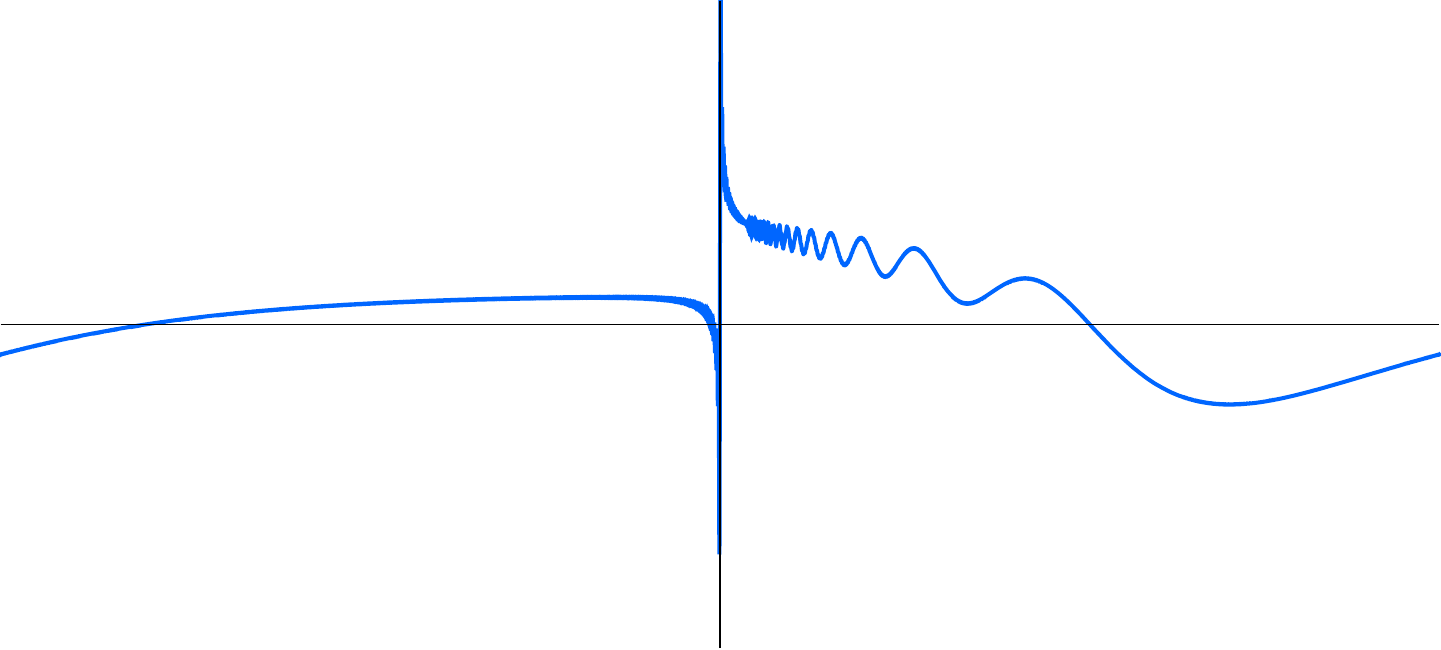}{7.5}\vskip-20pt
\Fiiil{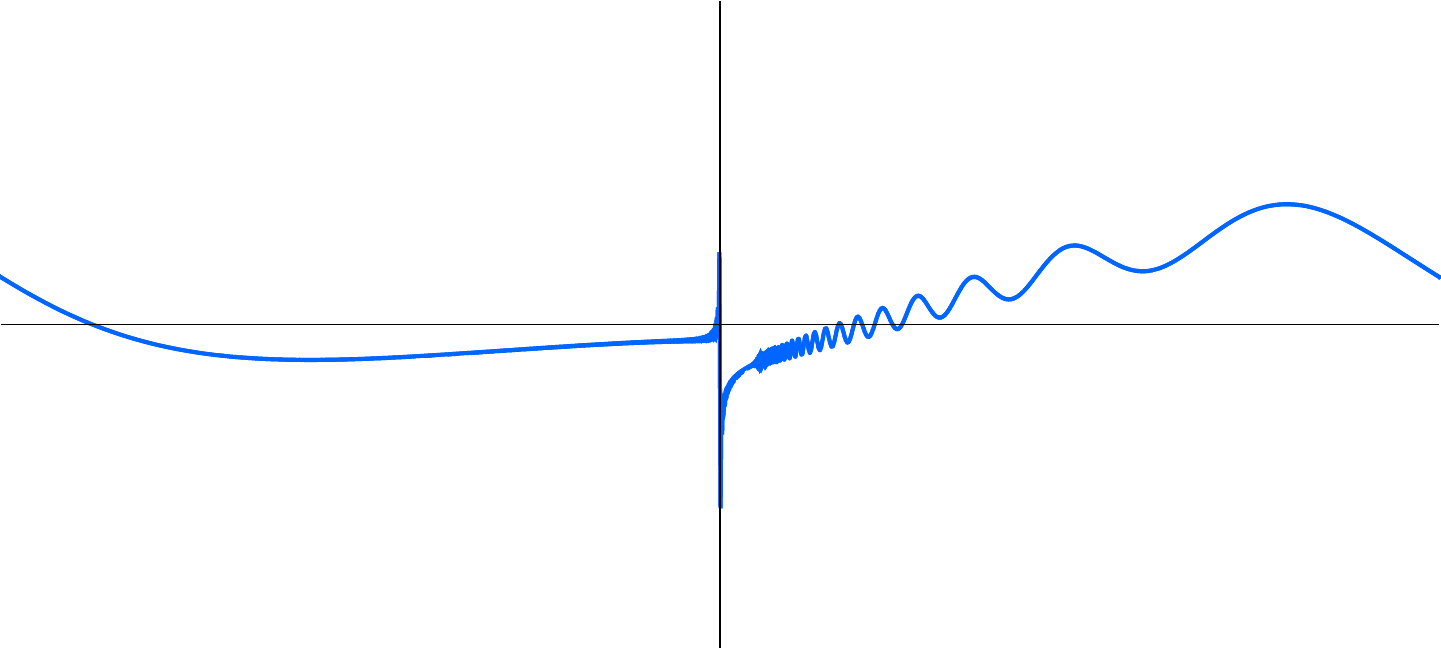}{10.}{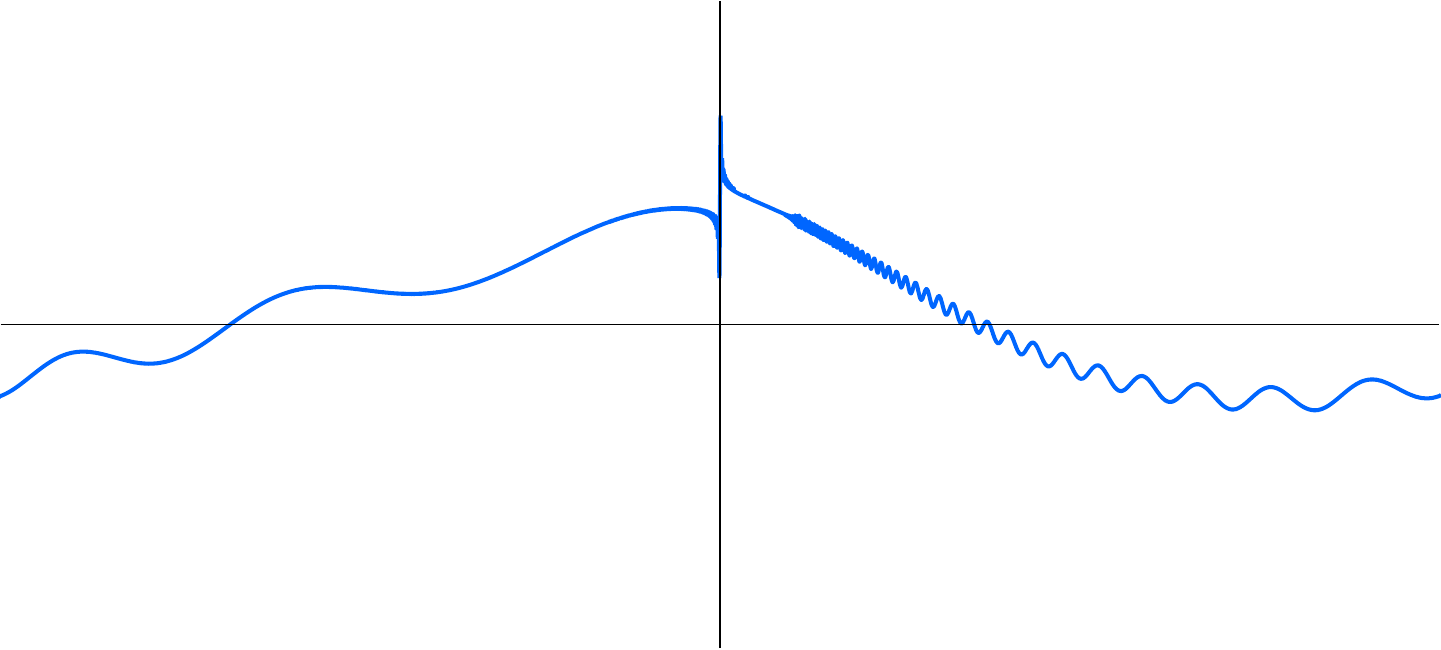}{20.}{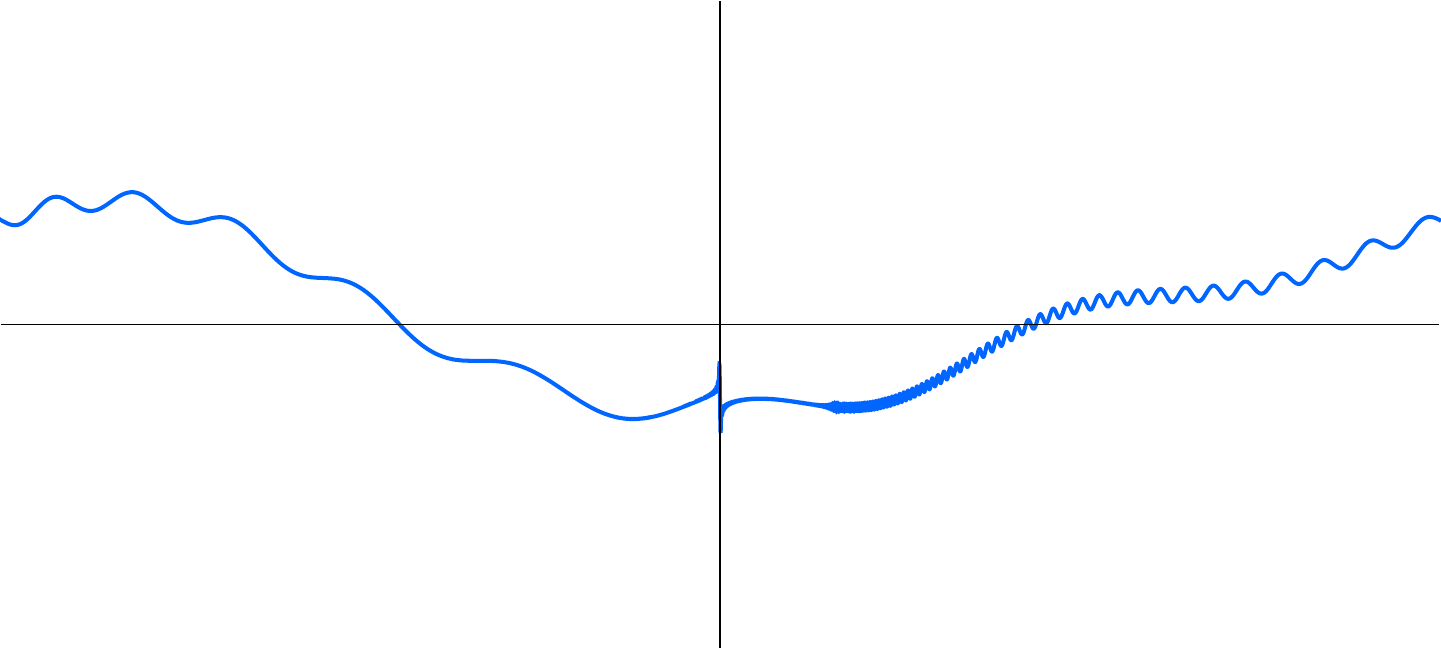}{30.}\vskip-20pt
\Fiiil{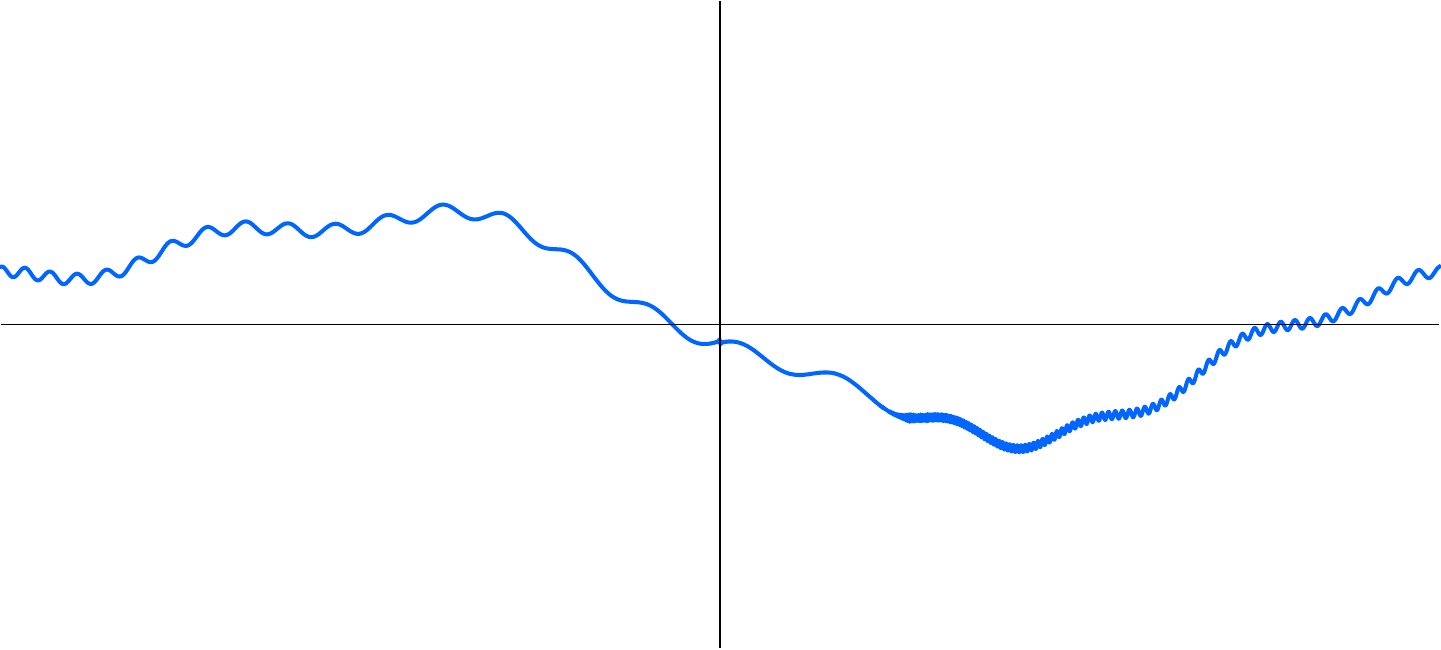}{50.}{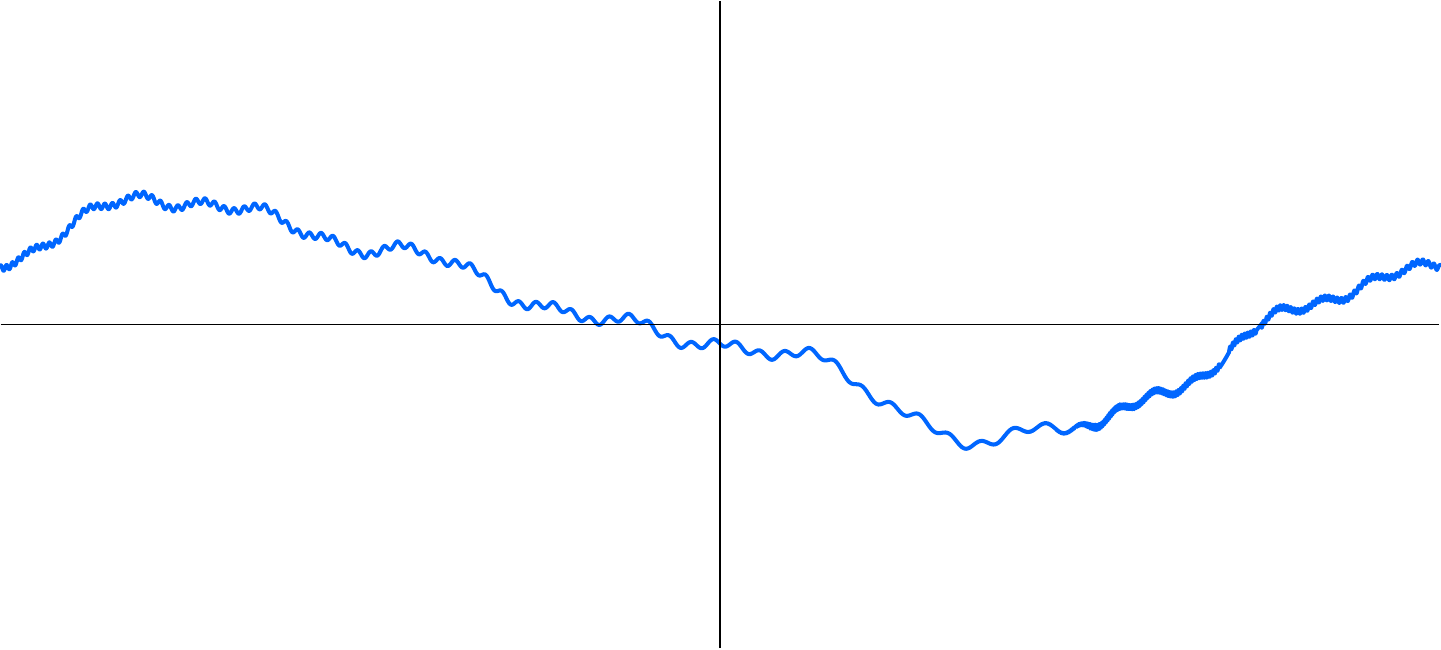}{100.}{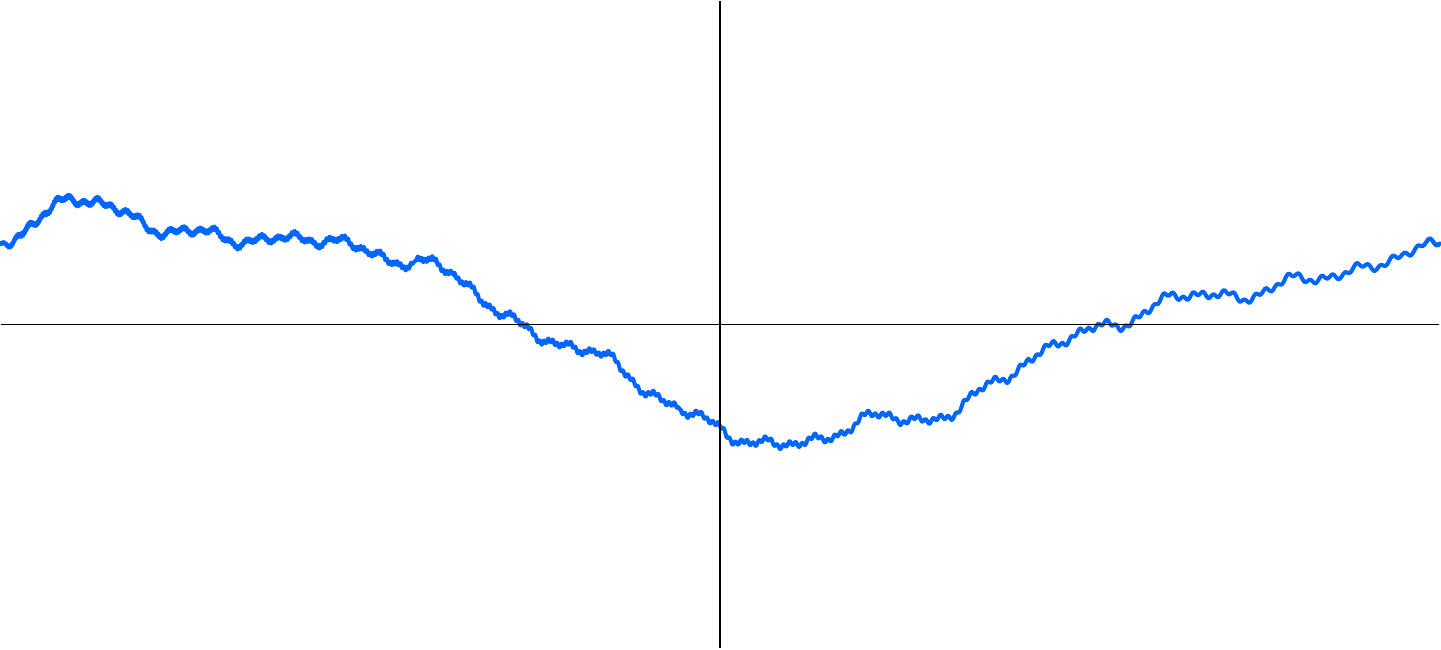}{200.}
}


\Fig{11}{The Unidirectional Dispersive Periodic Lamb Oscillator\break for $\omega(k)={k^2}/\paz{1+\f3\:k^2}$}
{\prefigskip=5pt\figlabl=4.75cm\figwidth=4.5cm 
\Fiiil{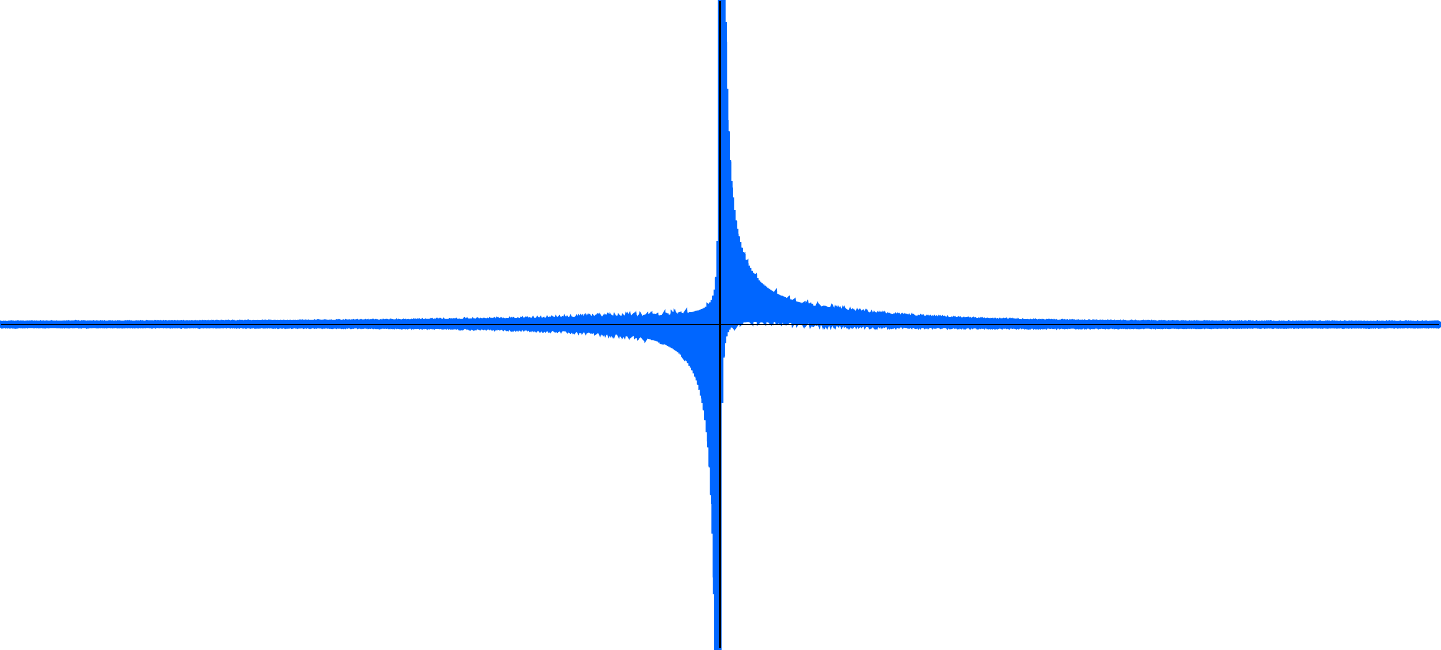}{1.}{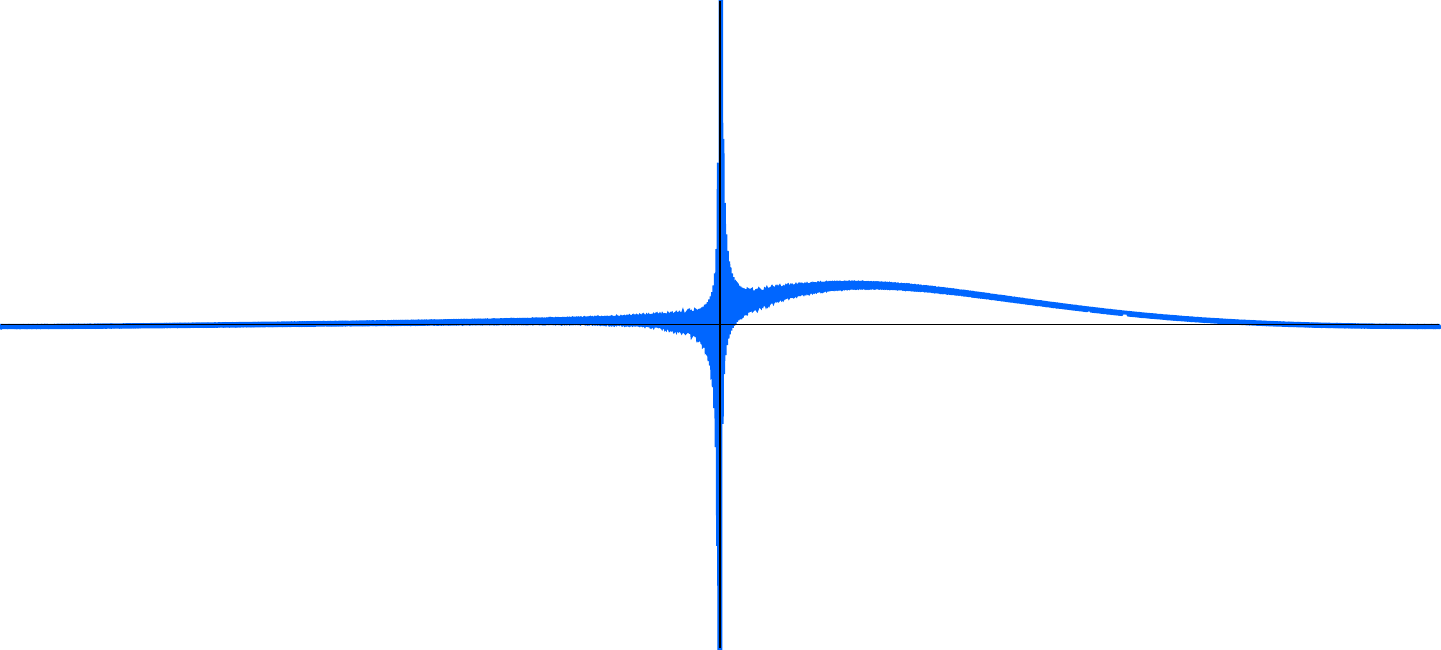}{2.5}{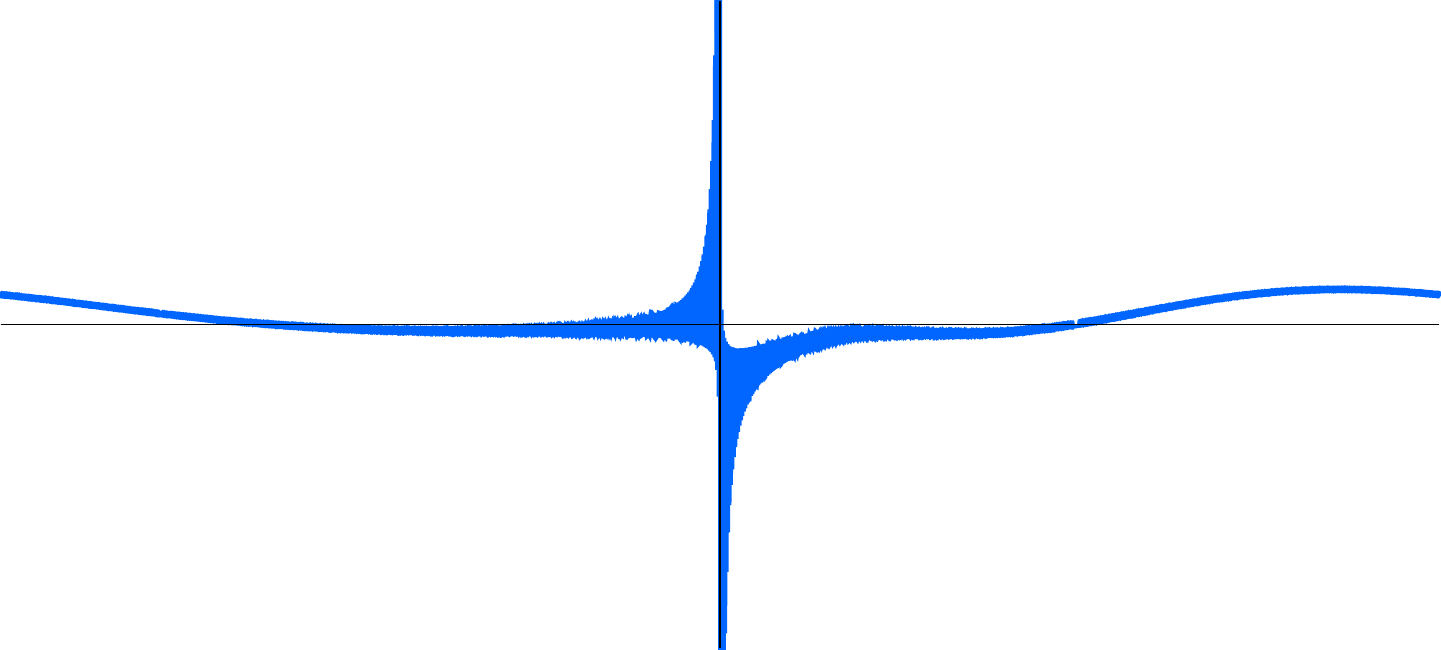}{5.}\vskip-20pt
\Fiiil{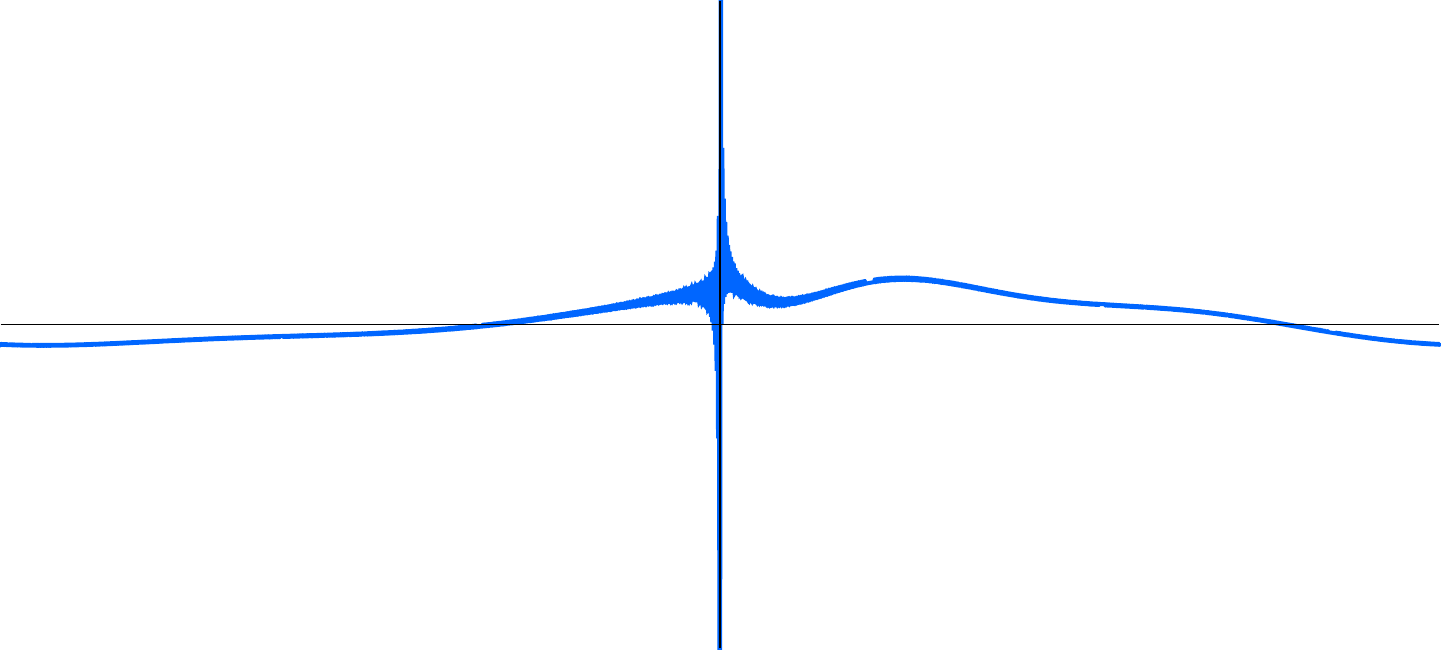}{10.}{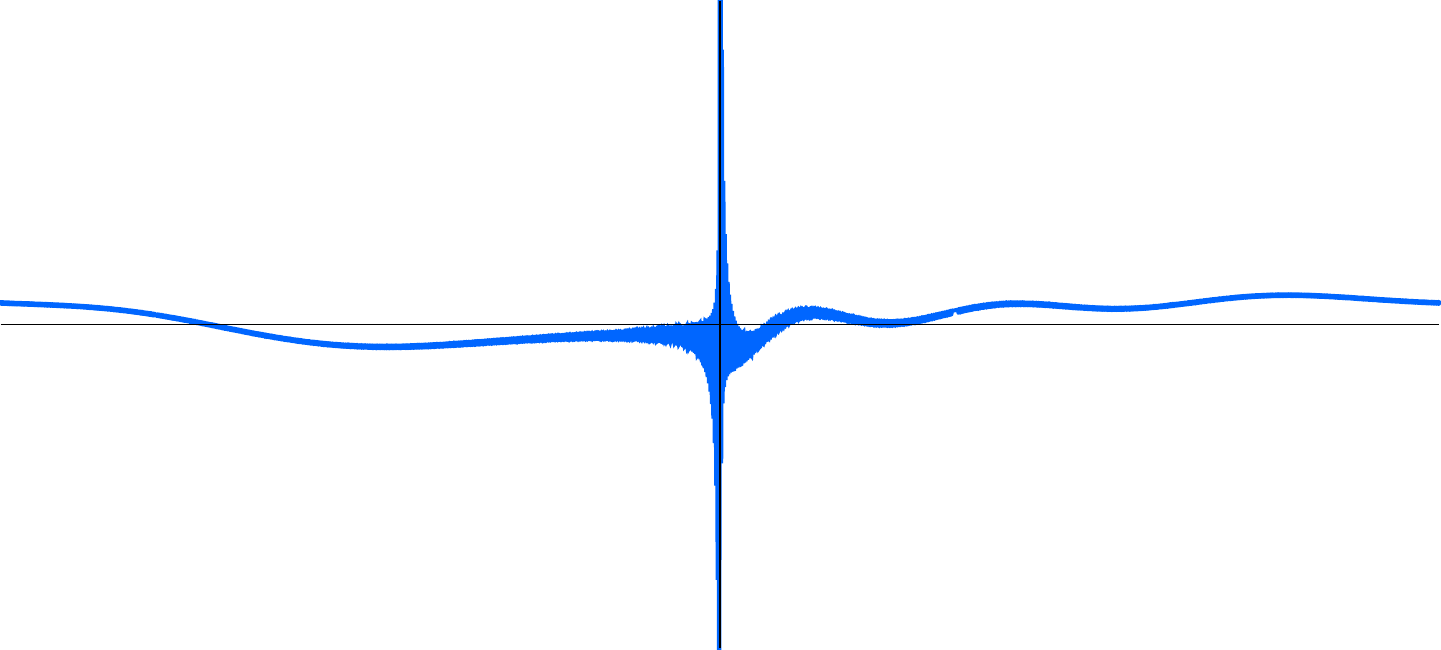}{20.}{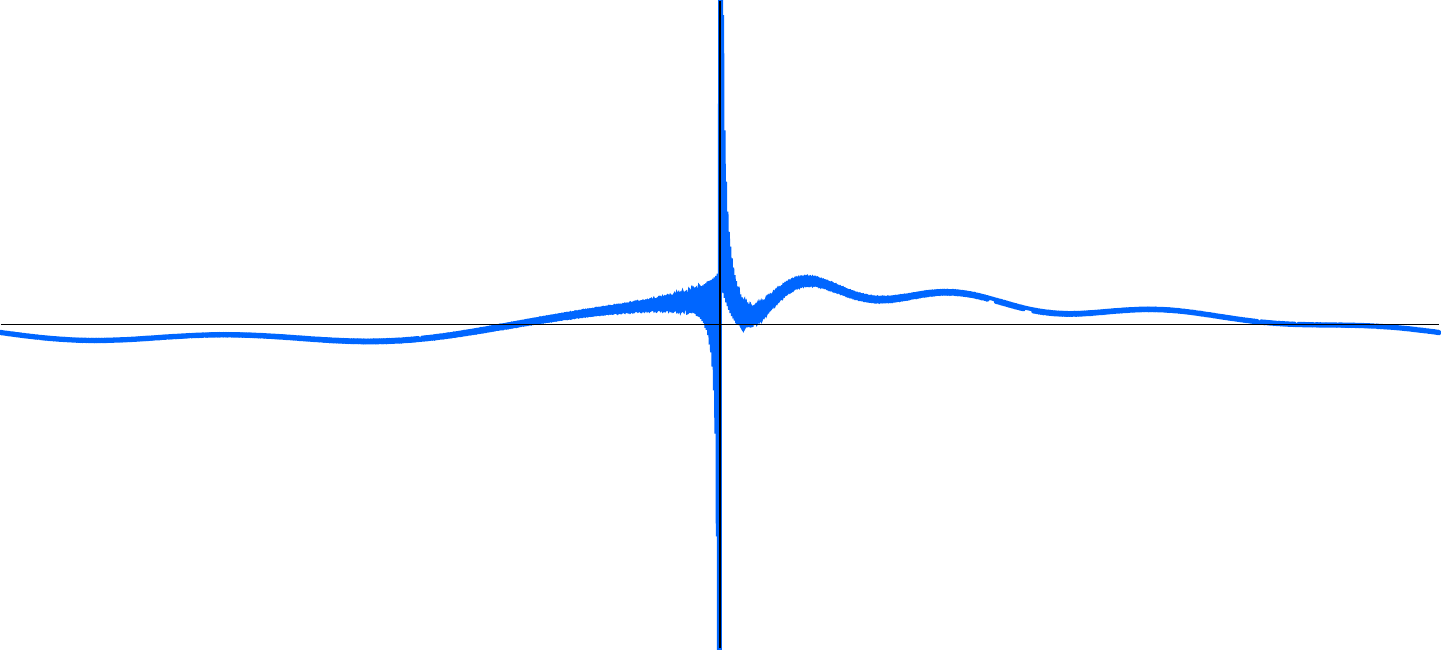}{35.}\vskip-20pt
\Fiiil{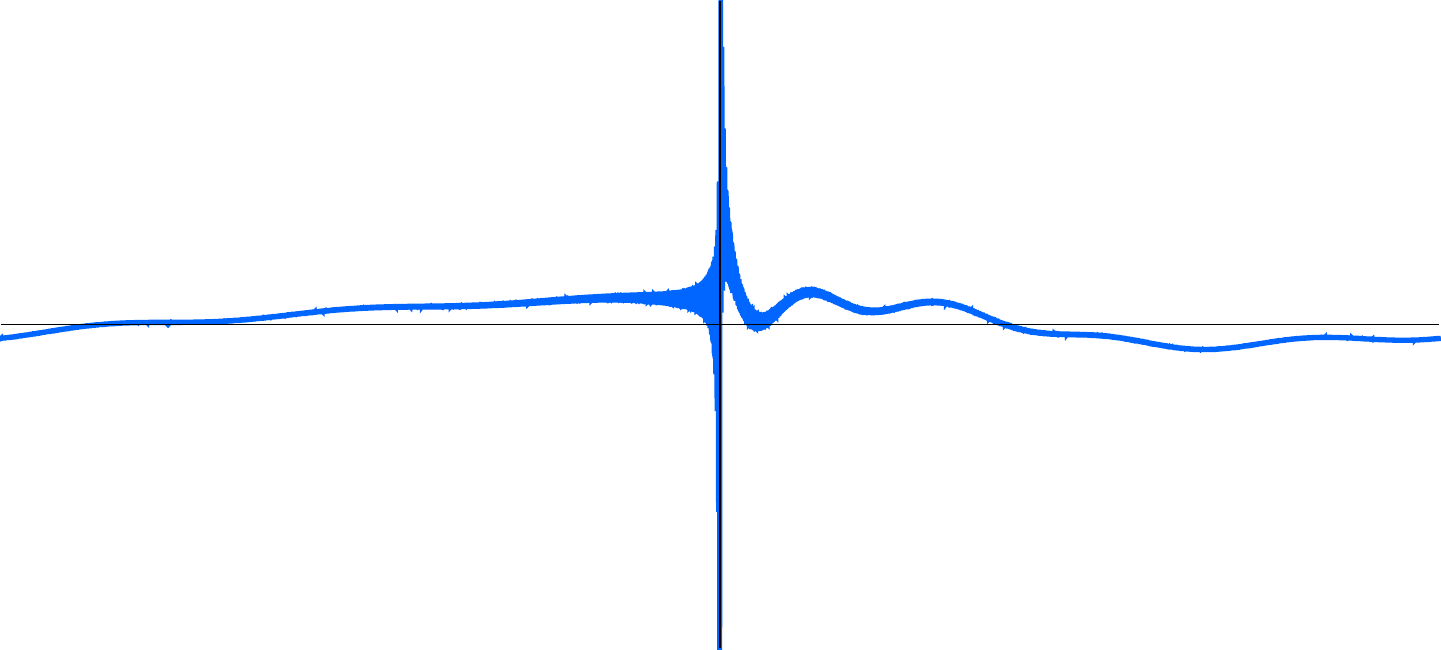}{50.}{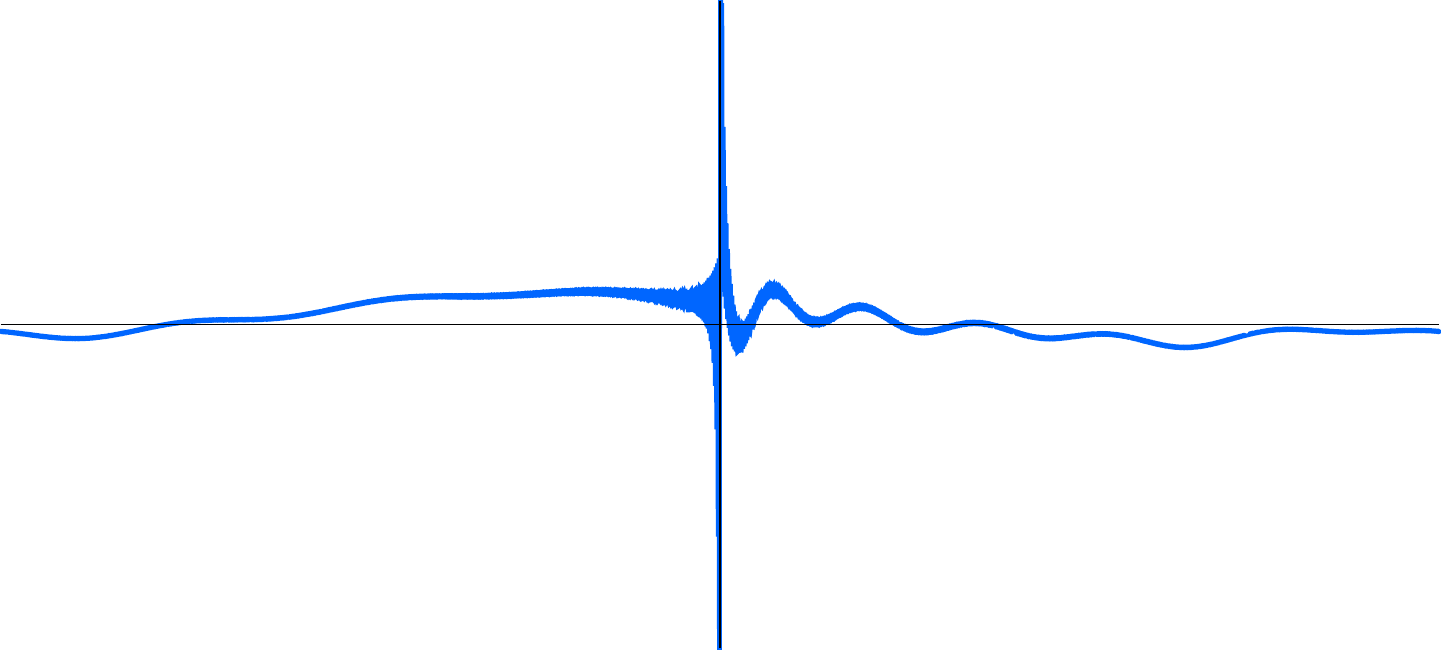}{100.}{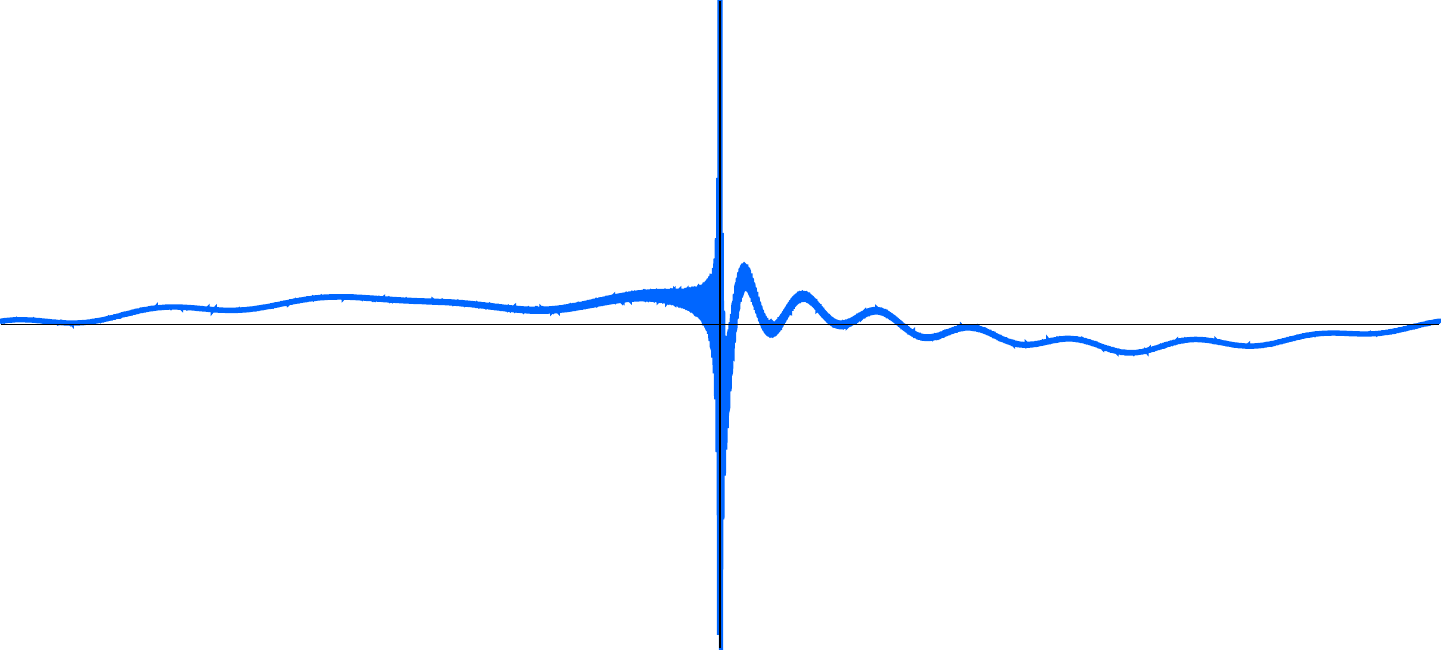}{200.}
}


\Section c Conclusions and Future Research Directions.

The most striking result in this study is that, when subject to periodic boundary conditions, a Lamb oscillator will generate fractal solution profiles in a dispersive medium possessing an asymptotically sub-linear large wave number dispersion relation whereas in an unforced dispersive medium, fractalization and dispersive quantization require that the dispersion relation grow superlinearly. 



Several additional directions of research are indicated. As noted above, it has been shown, \rf{HaBlWe}, that coupling the wave equation or Klein--Gordon equation to a gyroscopic oscillator or more general Chetaev system can induce instabilities through the effects of Rayleigh dissipation. Preliminary numerical experiments with more general dispersive media on periodic domains indicate that the coupling does not appear to destabilize an otherwise stable gyroscopic system. However, a more thorough analysis must be performed to confirm this result.

Furthermore, the effect of the boundary conditions on both unforced and forced dispersive wave equations is not yet clear, although in the unforced regime some preliminary results are available, showing that the fractal behavior of solutions is highly dependent upon the form of the boundary conditions. The corresponding problems in higher dimensional dispersive media also remain largely unexplored.

%




\Section a Acknowledgements.

The authors  wish to thank the referees for their careful reading of the manuscript, correcting some errors, as well as a number of useful suggestions.  We are also grateful to M.~Burak Erdo\utxt gan for alerting us to his recent paper \rf{ErSh}, and for valuable correspondence on the convergence question.

\vskip 25pt


\end{document}